\documentclass[11pt,notitlepage,twoside]{article}
\usepackage{latexsym}
\usepackage{amsmath}
\usepackage{amsfonts}
\usepackage{amssymb}
\usepackage{enumerate}
\usepackage[latin1]{inputenc} 
\renewcommand{\a }{\alpha }

\renewcommand{\d}{\delta }

\newcommand{\D }{\Delta }

\newcommand{\e }{\varepsilon }
\newcommand{\eps }{\varepsilon }

\renewcommand{\l }{\lambda }

\newcommand{\n }{\nabla }

\renewcommand{\O }{\Omega }

\newcommand{\ov}{\overline}

\newcommand{\be}{\begin{equation}}
\newcommand{\ee}{\end{equation}}
\newenvironment{pf}{\noindent{\bf Proof.}\enspace}{
\hfill$\Box$\medskip}
\newenvironment{pfn}[1]{\noindent{\bf Proof of {#1}\enspace}}{
\hfill$\Box$\medskip}
\newcommand{\R}{\mathbb{R}}
 
\newtheorem{thm}{Theorem}[section]
\newtheorem{pro}[thm]{Proposition}
\newtheorem{lem}[thm]{Lemma}

\newtheorem{df}[thm]{Definition}
\numberwithin{equation}{section} 

\title { \Large \textbf{Stable critical point of the Robin function and\\ bubbling
phenomenon for a slightly subcritical elliptic problem
 }}
\author{
{\bf\large Habib Fourti} 
{\bf\large}\vspace{1mm}\\
{\it\small Department of Mathematics and Statistics, College of
Science, King Faisal University, Al-Ahsa 31982, Saudi Arabia},\\
{\it\small University of Monastir, Faculty of Sciences of
Monastir, Tunisia},\\ {\it\small Laboratory LR 13 ES 21, University of Sfax, Faculty of Sciences of Sfax, Tunisia}\\
{\it\small e-mail: hfourti@kfu.edu.sa, habib.fourti@fsm.rnu.tn}
\\
\vspace{1mm}{\bf\large Rabeh Ghoudi}
{\bf\large}\vspace{1mm}\\
{\it\small University of Gabes, Faculty of Sciences of Gabes,
Tunisia}\\
{\it\small e-mail: ghoudi.rabeh@yahoo.fr}\vspace{1mm}}

\begin{document}

\date{}

\maketitle

{\footnotesize \noindent {\bf Abstract.} In this paper, we deal with
 the boundary value problem
$-\Delta u= |u|^{4/(n-2)}u/[\ln (e+|u|)]^\e$ in a bounded smooth
domain $ \Omega$ in $\mathbb{R}^n$, $n\geq 3$ with homogenous
Dirichlet boundary condition. Here $\e>0$.\\ Clapp et al. in Journal
of Diff. Eq. (Vol 275) built a family of solution blowing up if
$n\geq 4$ and $\e$ small enough.
 They conjectured in their paper the existence of sign changing
 solutions which blow up and blow down at the same point.
 Here we give a confirmative answer by proving that our
 slightly subcritical problem has a solution with the shape of sign
 changing bubbles concentrating on a stable critical point of
 the Robin function for $\e$ sufficiently small.

 }

{\footnotesize\section{Introduction and results}
 Let $\Omega$ be a smooth open bounded domain in $\mathbb{R}^n, \ n\geq 3$. This paper deals
with the analysis of solutions of the boundary value problem
\begin{eqnarray*} (P_{\e})\qquad \qquad  \begin{cases}
 -\Delta u= f_{\e}(u)    & \mbox{in} \, \, \O,  \\
 \quad \ \ u= 0     &  \mbox{on}\, \, \partial \O,
\end{cases}
\end{eqnarray*}
where
\begin{equation}\label{nonpower}
f_\e(u):=\frac{|u|^{p-1}u}{[\ln (e+|u|)]^\e}, \ \e\geq 0, \hbox{ and
} p=\frac{n+2}{n-2}.
\end{equation}
Here  $p+1=2n /(n-2)$ is the critical Sobolev exponent for the
embedding of $ H^1_0({\O})$ into $L^{p+1}(\O)$. Observe that
solutions for $(P_\e)$ correspond to critical points in $H^1_0(\O )$
of the energy functional
\begin{align}\label{g0}&I_\e(u)=\frac{1}{2}\int_{\O}|\n
u|^2-\int_{\O}F_\e(u), \hbox{ where }F_\e(s)=\int_0^sf_\e(t)dt.
\end{align}
The nonlinearity $f_\e$ where $\e>0$ has been studied in several
works (see for instance \cite{BFG},
\cite{Pistoianew},\cite{Har},\cite{LLX},\cite{Pardo}). The existence
of solutions for $(P_\e)$ when $\e>0$ is guaranteed by a variational
argument developed in \cite{Har}. Despite the noncompactness of the
embedding of $ H^1_0({\O})$ into $L^{p+1}(\O)$, the Euler functional
$I_\e$ satisfies Palais Smale condition. Thus problem $(P_\e)$ is
considered as a slightly
subcritical one.\\
The critical case corresponds to $\e=0$ and we get
$f_0(u)=|u|^{p-1}u$. In such a situation, the Euler functional $I_0$
does not satisfy Palais Smale condition and thus the corresponding
variational problem happens to be lacking of compactness. The
obtained problem $(P_0)$ is especially meaningful in geometry. In
fact, the Yamabe problem is a version of this problem on manifolds.\\
In the sequel we list some known results concerning the critical
case $\e=0$.
 When $\Omega$ is
starshaped, Pohozaev proved in \cite{Pohozaev} that problem $(P_0)$
has no positive solutions. However, when $\Omega$ is an annulus,
Kazdan and Warner provided in \cite{KW} the existence of a positive
radial solution. By using
 critical points at infinity theory, Bahri and Coron \cite{BC} proved that such a problem has a positive
solution, under the condition that $\Omega$ has nontrivial topology.\\
In this paper, we deal with the slightly subcritical case, i.e. $\e
> 0$. In order to state old and new results, it is useful to recall
some well known definitions.\\
The space $H^1_0(\O)$ is equipped with the norm $\|.\|$ and its
corresponding inner product $\left<.,.\right>$ is defined by
$$ \|u\|^2=\int_{\O}|\n
u|^2;\quad \left<u,v\right>=\int_{\O}\n u \n v,\quad u,v \in
H^1_0(\O).$$ For $a\in \O$ and $\l>0$, let
 \begin{equation}\label{d}
\d_{(a,\l)}(y)=\frac{c_0\l^{(n-2)/2}}{\left(1+\l^2|y-a|^2\right)^{(n-2)/2}},
\hbox{ where } c_0:=(n(n-2))^{(n-2)/4}.
 \end{equation}
The constant $c_0$ is chosen such that $\d_{(a,\l)}$ is the family
of solutions of the following problem
\begin{equation}\label{e}
 -\Delta u=u^{(n+2)/(n-2)},\;\;u>0\; \;\mbox{in}\; \;\R^n.
\end{equation}
Notice that the family $ \d_{(a,\l)}$ achieves the best Sobolev
constant \begin{equation}\label{S_n}
S_n:=\inf\{\|\n
u\|^2_{L^2(\R^n)}\|u\|^{-2}_{L^{2n/(n-2)}(\R^n)}:u\not\equiv 0, \n u
\in (L^{2}(\R^n))^n \,\mbox{and}\;u \in L^{2n/(n-2)}(\R^n)\}.
\end{equation}
 We denote
by $P\d_{(a,\l)}$ the projection of $\d_{(a,\l)}$ onto $ H^1_0(\O)$,
defined by
\begin{eqnarray}\label{q1}
-\Delta P\d_{(a,\l)}=-\D \d_{(a,\l)} \mbox{ in } \O, \quad
P\d_{(a,\l)}=0 \mbox{ on}\;
\partial \O.
\end{eqnarray}
We will denote by $G$ the Green's function and by $H$ its regular
part, that is \begin{equation}\label{defg}G(x,y)=|x-y|^{2-n}-H(x,y)
\quad\mbox{ for }\,\, (x,y)\in \O^2,\end{equation}
 and for $x\in\O$,
$H$ satisfies
\begin{eqnarray*}   \begin{cases}
\Delta  H(x,.)=0   & \mbox{in} \, \, \O,
\\H(x,y)=|x-y|^{2-n},&\mbox{for}\, \,y\in
\partial \O.
\end{cases}
\end{eqnarray*}
We define the Robin function as $R(x)=H(x, x)$, $x\in \O$.\\
The first paper establishing the existence of blowing up solutions
to problem $(P_\e)$ is \cite{Pistoianew} where the authors proved
that any $x_0$ non-degenerate critical point of $R$ generates a
family of single-peak solutions concentrating around $x_0$ as $\e$
goes to $0$. This family of solutions has  the following form
\begin{equation*}
u_\e= P\delta_{( a_{\e} ,\lambda_{\e})}+v_\e,
\end{equation*}
where $a_{\e}\rightarrow x_0, \,\lambda_{\e}\rightarrow\infty$ and
$\|v_\e\|\rightarrow 0$ as $\e$ goes to $0$.\\
 In \cite{Pardo}, the asymptotic
behavior of radially symmetric solutions of $(P_\e)$ was analyzed
when $\Omega$ is a ball.
Very recently, we provided in \cite{BFG} the existence of positive
as well as changing sign solutions that blow up and/or blow down at
different points in $\Omega$. Our solutions have the following
expansion
\begin{equation*}
u_\e= \sum_{i=1}^k\alpha_{i,\e}\gamma_iP\delta_{( a_{i,\e}
,\lambda_{i,\e})}+v
\end{equation*}
where $\gamma_i\in \{-1, 1\}$  for each $1\leq i \leq k$,
$\a_{i,\e}\rightarrow 1$ $\lambda_{i,\e}\rightarrow\infty$, $v_\e\in
E_{(a,\l)}$ and $\|v_\e\|\rightarrow 0$ as $\e$ goes to $0$. Here
\begin{align}\label{vorthogonal}&E_{(a,\l)}:=
\Big\{v\in H^1_0(\O):\bigg<v,P\d_{i}\bigg>=\left<v,\frac{\partial
P\d_{i}}{\partial \l_i}\right>= \left<v,\frac{\partial
P\d_{i}}{\partial (a_i)_j}\right>=0\;\forall \; 1\leq j\leq
n,\,\forall \,1\leq i\leq k\Big\},
\end{align}
where $P\d_i=P\d_{(a_{i,\e},\l_{i,\e})}$ and $( a_{i,\e} )_j$ is the
$j^{th}$ component of $ a_{i,\e} $. Furthermore, the concentration
points $a_{1,\e},\ldots,a_{k,\e}$ are far away from each other and
converge to distinct points in $\O$, which form a non-degenerate
critical point of a function defined explicitly in
terms of the Green function and its regular part.\\
In this paper, we focus on a new kind of solutions. In fact, we
prove the existence of sign-changing bubble towers solutions. These
solutions are constituted by superposition of positive and negative
bubbles with different blow up orders and closed concentration points.\\
Musso and Pistoia have proved in \cite{MP1} the existence of such
solutions for
 the following nonlinear subcritical elliptic problem
\begin{eqnarray}\label{subpb}\begin{cases}
 -\Delta u= |u|^{p-1-\e}u    & \mbox{in} \, \, \O,  \\
 \quad \ \ u= 0     &  \mbox{on}\, \, \partial \O,
\end{cases}
\end{eqnarray}
where $\e>0$. As mentioned in \cite{Pistoianew} and \cite{BFG},
there is an interesting analogy between results obtained for
$(P_\e)$ and those known for the usual elliptic equation
\eqref{subpb}. In this direction, we prove here that $(P_\e)$ shares
also with problem \eqref{subpb}, the existence of tower bubbles
solutions of different blow-up orders, as conjectured in
\cite{Pistoianew}.\\
Let us review some known facts related to tower bubbles solutions.
This phenomenon was firstly discovered by del Pino, Dolbeault and
Musso in \cite{DDM}, where the authors dealt with the slightly
supercritical Brezis-Nirenberg problem. They proved the existence of
positive bubble tower radial solutions obtained by a superposition
of several bubbles centered at the origin point but
with different scaling parameters.\\
Their method strongly relies on the symmetry of the problem. Later,
the construction described above was extended by Ge, Jing and Pacard
in \cite{GJP} to a more general class of domains, namely domains
such that the associated Robin function $R$ admits a non-degenerate
critical point in the case of one tower of bubble. They also proved
the existence of multiple towers of bubbles under a
 non-degeneracy condition on a critical point of a certain functional of point of the domain. 
In \cite{MP1}, the authors studied the same problem and they were
able to remove the
assumption on the domain. 
Indeed, they proved that, in any bounded domain $\O$, the slightly
supercritical Brezis-Nirenberg problem does admit solutions with the
shape of a tower of bubbles. 
 Their
idea was the following: given any domain $\O$, its Robin function is
smooth, positive and unbounded as $x$ approaches the boundary, thus
it has a minimum in $\O$, and
hence at least one critical point.\\
In the same paper \cite{MP1}, as we mentioned before, Musso and
Pistoia have constructed sign-changing bubble tower solutions for
problem $\eqref{subpb}$ which blow up at the minimum of the Robin
function $R$. This work was an extension of that of Pistoia and Weth
\cite{PW} where they proved that if $\Omega$ is symmetric with
respect to the $x_1,\ldots, x_n$ axes, problem \eqref{subpb} has a
sign-changing solution with the shape of a tower of bubbles with
alternate signs, centered at the center of symmetry of the domain.
Finally, we mention the work in \cite{DG}, where the authors studied
the existence of bubble towers changing sign solutions for the
counterpart of \eqref{subpb} when the Laplacian operator is replaced by the biharmonic one.\\
Before introducing our result, we recall the following definition as
stated in \cite{MollePistoia}.
\begin{df}
Let $g:\Omega   \rightarrow \mathbb{R}$ be a $C^1-$function. We say
that $\xi_0$ is a stable critical point of $g$ if $\nabla g(\xi_0)=
0$ and there exists a neighborhood $U\subset\subset \Omega$ of
$\xi_0$ such that $\nabla g(x)\neq 0, \,  \forall x\in \partial U,$
$$\nabla g(x)=
0, \, x\in U \Rightarrow  \, g(x)=g(\xi_0),$$ and
$$\deg (\nabla g, \overline{U},0)\neq 0,$$
where $\deg$ denotes the Brouwer degree.
\end{df}
It is easy to see that, if $\xi_0$ is an isolated minimum or maximum
point of the function $g$, then $\xi_0$ is a stable critical point
of $g$. Moreover, a non-degenerate critical point of $g$ is stable
according to the previous definition.\\

 As we previously mentioned, the aim of the current paper
is, by applying the Liapunov-Schmidt reduction method, proving the
existence of bubble tower sign-changing solutions of $(P_\e)$
concentrating at a stable critical point of the Robin function $R$.
Precisely speaking, our main result can be stated as follows.
\begin{thm}\label{t:1}Assume that $n\geq 3$ and let $\xi_0$ be a stable critical
point of the Robin function $R$. Then, given an integer $k\geq 1$,
there exists $\e_0>0$ such that for each $\e\in (0,\e_0)$, problem
$(P_\e)$ has a solution $u_\e$ of the form
\begin{equation} \label{h:18}
u_\e=
\sum_{i=1}^k(-1)^i\alpha_{i,\e}P\delta_{(\xi_{i,\e},\lambda_{i,\e})}+v_\e,
\end{equation}
where $\xi_{i,\e}=\xi_\e+\frac{1}{\l_{i,\e}}\sigma_{i,\e}\hbox{ for
} i=1,\ldots,k, \ \sigma_{1,\e}=0;\,\sigma_{i,\e}\in \R^n $ for
$i=2,\ldots,k$ and $v_\e\in E_{(\xi,\l)}$.\\ Furthermore, as
$\e\rightarrow 0$, we have
\begin{align*}&\a_{i,\e}\rightarrow 1\hbox{
for } i=1,\ldots,k;\, \xi_\e\rightarrow \xi_0 ; \,
\sigma_{i,\e}\rightarrow 0 \hbox{
for } i=2,\ldots,k
\end{align*}
and the variables $\l_{i,\e}$ and $v_\e$ satisfy
$$
\ \lambda_{i,\e}^{n-2}\left(\frac{
\e}{|\ln\e|}\right)^{2(k-i)+1}=C(k,i,\xi_0)(1+o(1))
\hbox{ for } i=1,\ldots,k;$$
$$\|v_\e\|=
\left\{\begin{array}{ll}O(\e\ln|\ln\e|
),&\hbox{ if }n\leq 6;\\
O\left(\left(\displaystyle\frac{\e}{|\ln\e|}\right)^\frac{n+2}{2(n-2)}\right),
&\hbox{ if }n> 6.\end{array}\right.
$$
Here $C(k,i,\xi_0)$ is positive constant which depends on $k,\, i$
and $\xi_0$.
\end{thm}
We point out that one can remove the assumption $R$ admits a stable
critical point by using the same idea of Musso and Pistoia namely
that the
 Robin function has a minimum in any given domain $\O$, and hence it has at least
 one stable
critical point.\\
Arguing as in \cite{MP2}, one can exhibit an example of a
contractible domain for which such an assumption holds true. Indeed,
following the idea of perturbing domains, Musso and Pistoia
constructed a domain for which the function $R$ has a stable
critical
point.\\
We also mention that the non-degeneracy of the critical point of the
Robin function implies the existence of this kind of solutions,
since by applying implicit function theorem one may prove that  a
non-degenerate critical point is a stable one.
 Thus, our condition is weaker than the
non-degeneracy condition. In fact, when $k=1$ our result improves
the one of Clapp
et al. in \cite{Pistoianew}.\\
 Note that our family of solutions $(u_\e)$ converges weakly to zero and its blow up rate
 satisfies $$\|u_\e\|_{\infty}\sim c(|\ln \e|\e^{-1})^{k-1/2}\hbox{ as }\e\rightarrow
 0.
 $$ Furthermore, the concentration speeds $\l_i$'s verify  $$\l_i \hbox{ is of order }(|\ln
 \e|\e^{-1})^{(2(k-i)+1)/(n-2)},\; \hbox{ for }i=1,\ldots, k \hbox{ and }\frac{\l_i}{\l_{i+1}}\sim c (|\ln
 \e|\e^{-1})^{2/(n-2)},\; \hbox{ for } i=1,\ldots, k-1.$$
 This choice will be justified by exploiting some balancing
 conditions for the parameters of the concentration given in
 Lemmas \ref{D2} and \ref{D3}.
 Indeed, taking $\l_k$ the smallest concentration
 speed by analyzing these balancing conditions, we obtain that $(|\ln
 \e|\e^{-1})^{1/(n-2)}$ and $\l_k$
 needed to be of the same order.
 This will be the subject of Proposition \ref{t:4} where
 we investigate the asymptotic profile of a family
of sign changing bubble tower solutions $u_\e$ which blows up in the
interior of the domain $\Omega$.
Let us point out that, the obtained relations
 between the parameter $\l_i$'s and $\e$ goes along the
 argument developed in \cite{MP1} concerning \eqref{subpb}. In fact, Musso and Pistoia have chosen that the
 smallest concentration speed satisfies $\l_k\sim c \e^{-1/(n-2)}$ instead of the quantity $c(|\ln
 \e|\e^{-1})^{1/(n-2)}$ in our case.\\

To prove our main result, we apply the Liapunov-Schmidt Reduction
method (see for instance \cite{Pistoia} and the references therein).
Thanks to the analysis of the gradient of the energy
 functional performed in \cite{BFG}, we were able to adopt some arguments developed in \cite{MP1} and that
  by choosing suitable concentration rates and concentration points.
  After this choice we have also improved the expansion of the
  gradient with respect to the new variables.\\
  Note that our proof requires also the expansion of the energy
  functional $I_\e$ which is given in Proposition \ref{p23'}. In contrast with the nonlinear term $|u|^{p-1-\e}u$ in problem \eqref{subpb},
   one can
 not write the explicit expression of $F_\e$ the antiderivative
of the nonlinearity $f_\e$ of problem $(P_\e)$ (see \eqref{nonpower}
and \eqref{g0}). We were able to overcome this technical difficulty
and to obtain the asymptotic expansion of $I_\e$, as $\e$ goes to
$0$. From this
  expansion we get a new functional defined on finite space. By
  applying degree theory and taking into account the stability of
  the  critical point of the Robin function we conclude on this new
  functional.\\
  We mention that our new result in Proposition \ref{p23'} is quite
  involving.
   We think that this expansion will be useful to study the existence and to describe the blow up profile of some
positive and changing sign
 solutions, as $\e$ goes to zero.\\

The paper is organized as follows: In Section $2$,
 we collect some basic tools which includes the
asymptotic expansions of the energy functional and its gradient. We
would like to warn the reader that the expansion of the gradient is
quoted from \cite{BFG} and we did not repeat the proofs. However, we
provide the proof of the expansion of the energy functional in case
of close concentration points.
The description of the appropriate solutions is obtained in Section
$3$ and that by analyzing the asymptotic profile of sign changing
bubble tower solutions. In Section $4$ we introduce the precise
profile of our searched solutions and we study its remainder term. 
Then we derive an asymptotic expansion for the reduced energy
functional in Section $5$. Finally, we complete the proof of our
main result in Section $6$. Section $7$ is an appendix where we collect some technical Lemmas used in this work.   \\
Throughout this paper, we use the same $c$ to denote various generic
positive constants independent of $\e$.
\section{The Technical Framework}
Recall that each critical point of the energy functional $I_\e$
(defined by \eqref{g0}) is a solution of $(P_\e)$. To construct
solutions of the form \eqref{h:18}, our argument require the
expansion of $I_\e$ and its gradient for $u=
\sum_{i=1}^k\alpha_i\gamma_iP\delta_{(\xi_i,\lambda_i)}+v$ i.e. $u$
belongs to a neighborhood of potential concentration sets. Here
$\gamma_i\in \{-1,1\}$. We mention that the expansions introduced in
this section will be given in general setting i.e. we will not use,
at this stage, the
parameters information presented in Theorem \ref{t:1}. \\
For $\eta>0$, $k \in \mathbb{N}$ and $
(\gamma_{1},\dots,\gamma_{k})\in\{-1,1\}^k$, let us define
$$
\begin{array}{rcl} V(k,\eta)&=&\Big\{u \in H^1_0(\Omega)\ /\ \exists
\xi_{1},\ldots,\xi_{k} \ \in \O, \ \exists
\lambda_{1},\ldots,\lambda_k>\eta^{-1}, \ \exists
\alpha_1,\ldots,\alpha_k>0 \hbox{ with}
 \\
&& \parallel u-\sum_{i=1}^{k}\alpha_i
\gamma_iP\delta_{(\xi_i,\lambda_i)}\parallel<\eta;\ \mid \a_i-1\mid
<\eta, \ \l_id(\xi_i,\partial \Omega)>\eta^{-1} \ \forall i,\
\varepsilon_{ij}<\eta \ \forall i\neq j \Big\},
\end{array}$$
where \begin{equation}\label{eijdef}
\e_{ij}:=\left(\frac{\l_i}{\l_j}+\frac{\l_j}{\l_i}+\l_i\l_j|\xi_i-\xi_j|\right)^\frac{2-n}{2}
\end{equation}
Note that, the variable $\e_{ij}$ comes from the scalar product
\begin{equation}\label{e_ij}
\displaystyle\int_{\mathbb{R}^n}\nabla\delta_{(\xi_i,\lambda_i)}.
\nabla\delta_{(\xi_j,\lambda_j)}
=\int_{\mathbb{R}^n}\delta_{(\xi_i,\lambda_i)}^{p}
\delta_{(\xi_j,\lambda_j)} =O(\e_{ij})\hbox{ for }i \neq j
,\end{equation} (see \cite{B1} page 4).\\
For simplicity, we denote $P\delta_{(\xi_i,\lambda_i)}$ by
 $P\delta_i$.  Recall that
simple computations show that
\begin{equation}\label{majdd}
|\l_i\frac{\partial P\d_i}{\partial \l_i}|,\,
|\frac{1}{\l_i}\frac{\partial P\d_i}{\partial (\xi_i)_\ell}|\leq
cP\d_i\leq c \d_i \hbox{ and } |\l_i\frac{\partial\d_i}{\partial
\l_i}|,\, |\frac{1}{\l_i}\frac{\partial\d_i}{\partial (\xi_i)_\ell}|
\leq c \d_i\hbox{ for }\ell=0,\ldots,n.
 \end{equation}
 As in \cite{BFG}, we are looking
for a solution of $(P_\e)$ in a small neighbourhood of
$\sum_{i=1}^{k}\alpha_i \gamma_iP\delta_i$.  The authors in
\cite{BFG} investigated the gradient of the functional $I_\e$ in
$V(k,\eta)$. In the sequel, we recall some expansions
 extracted from \cite{BFG}.\\
\begin{pro}\label{p22}
Let $n\geq 3$ and $u=\sum_{j=1}^k\a_j\gamma_jP\d_j+v\in V(k,\eta)$
such that $v\in E_{(\xi,\l)}$. For each $i\in \{1,...,m\}$, we have
the following expansion
\begin{align*}\left<\n
I_\e(u),P\d_i\right>=\gamma_i\a_i\big(1-\a_i^{p-1}\big)S_n^{\frac{n}{2}}+O\Big(\e\ln(\ln\l_i)
+\sum_{j}\frac{1}{(\l_jd_j)^{n-2}}+ \sum_{j\neq
i}\e_{ij}+\|v\|\Big),
\end{align*}
where $S_n$ is the best Sobolev constant defined in \eqref{S_n},
$d_j:=d(\xi_j,\partial\O)$ and $\e_{ij}$ is introduced in
\eqref{eijdef}.\\
\end{pro}
\begin{pro}\label{p23}
Let $n\geq 3$ and $u=\sum_{j=1}^k\a_j\gamma_j P\d_j+v\in V(k,\eta)$
such that $v\in E_{(\xi,\l)}$. For each $i\in \{1,...,k\}$, we have
the following expansion
\begin{align*}
  \left< \nabla I_\e(u),\l_i\frac{\partial P\d_i}{\partial \l_i}\right>&=\gamma_i\Gamma_1\displaystyle
\frac{\a_i^p\e}{\ln \lambda_i}+
 (n-2)\overline{c}_1\frac{\gamma_i\a_i}{2}\left(1-2\a_i^{p-1}\right)\frac{H(\xi_i,\xi_i)}{\l_i^{n-2}}
 \\&+\overline{c}_1\sum_{j\neq i}\gamma_j\a_j\left(1 -\a_j^{p-1}-\a_i^{p-1}\right)
 \biggl(\l_i\frac{\partial \e_{ij}}{\partial
\l_i}+\frac{n-2}{2}\frac{H(\xi_i,\xi_j)}{(\l_i\l_j)^{(n-2)/2}}
\biggr)  \\
&+O\left(\frac{\e}{\ln(\l_i)^2}+\frac{1}{(\l_id_i)^{2n-4}}+\frac{1}{\l_i^{(1-\tau)n}d_i^{2n}}+\sum_{j=1}^k\e^2\ln(\ln\l_j)^2
\right)\\&+O\left( \sum_{j\neq
i}\e_{ij}^{\frac{n}{n-2}}\ln(\e_{ij}^{-1})+\sum_{j\neq
i}\e_{ij}^{2}\ln(\e_{ij}^{-1})^{\frac{2(n-2)}{n}}
+\|v\|^2\right)+\left\{\begin{array}{ll}O\left(\sum_{l=1}^k\frac{\ln(\l_l
d_l)}{(\l_l
d_l)^{n}}\right)  \quad\hbox{ if }n\geq 4;\\
O\left(\sum_{l=1}^k\frac{1}{(\l_l d_l)^{2}}\right) \quad \hbox{ if
}n=3 ,\end{array}\right.
\end{align*}
where
 $\overline{c}_1:=c_0^{\frac{2n}{n-2}}\int_{\mathbb{R}^n}\frac{1}{(1+|x|^2)^{(n+2)/2}}dx$, $\Gamma_1:=\displaystyle c_0\int_{\mathbb{R}^n}
\delta_{(0,1)}^p(y)\ln
(\delta_{(0,1)}(y))\frac{1-|y|^2}{(1+|y|^2)^{n/2}}dy=\frac{(n-2)^2}{4n}S_n^{n/2}$
and $\tau$ is a positive constant small enough.
\end{pro}

\begin{pro}\label{p24} Let $n\geq 3$ and $u=\sum_{l=1}^k\a_l\gamma_lP\d_l+v\in V(k,\eta)$ such that $v\in E_{(\xi,\l)}$.
 For each
 $i\in \{1,...,k\}$ and $j\in \{1,...,n\}$, we have the following expansion
\begin{align*}
&\left<\n I_\e(u),\frac{1}{\l_i}\frac{\partial P\d_i}{\partial
(\xi_i)_j}\right>=
 \gamma_i\left(\a_i^{p}-\frac{\a_i}{2}\right)
 \frac{\overline{c}_1}{\l_i^{n-1}}\frac{\partial H(\xi_i,\xi_i)}{
 \partial (a)_j}\\&
 +\overline{c}_1\sum_{l=1,l\neq i}^k\gamma_l\a_l\Big(1 -\a_l^{p-1}-\a_i^{p-1}\Big)
 \frac{1}{\l_i}\biggr(\frac{\partial \e_{il}}{\partial
(\xi_i)_j}-\frac{1}{(\l_i\l_l)^{(n-2)/2}}\frac{\partial H}{\partial
(a)_j}(\xi_i,\xi_l) \biggr)\\&+O\bigg(\sum_{l=1}^k\e^2\ln(\ln\l_l)^2
 +\sum_{l\neq i}\l_l | \xi_i-\xi_l|\e_{il}^{\frac{n+1}{n-2}
}+\sum_{l\neq i}\e_{il}^{\frac{n}{n-2}}\ln(\e_{il}^{-1})+\sum_{l\neq
i}\e_{il}^{2}\ln(\e_{il}^{-1})^{\frac{2(n-2)}{n}}+\|v\|^2\bigg)
\\&+\left\{\begin{array}{ll}O\left(\sum_{l=1}^k\frac{\ln(\l_l
d_l)}{(\l_l
d_l)^{n}}\right)  \quad\hbox{ if }n\geq 4;\\
O\left(\sum_{l=1}^k\frac{1}{(\l_l d_l)^{2}}\right) \quad \hbox{ if
}n=3 ,\end{array}\right.
\end{align*}
where $\frac{\partial H}{\partial (a)_j}$ denotes the partial
derivative of $H$ with respect to the $j-$th component of the first
variable.
\end{pro}
Now, our aim is to obtain the expansion of the energy functional
$I_\e$ around  $u\in V(k,\eta)$ with close concentration points. We
start by the following useful lemmas.
\begin{lem}\label{lemmaA1}
 Let $\xi\in\Omega$
 and $\lambda_i > 0$ be such that $ \lambda_i d(\xi,
 \partial\Omega)$ is very large. We denote
 \begin{equation}\label{xi}
 \xi_i:=\xi+\frac{\sigma_i}{\l_i}\hbox{ and
 }P\d_i:=P\d_{(\xi_i,\l_i)},
 \end{equation}
 where $\sigma_i$ belongs to a compact set in $\mathbb{R}^n$.\\
We have
\begin{align}
&\int_{\O}\left[f_\e(\a_i\gamma_iP\delta_i)-f_0(\a_i\gamma_iP\delta_i)\right](x-\xi).\n_xP\d_i\nonumber\\&=
 \frac{n-2}{2}\gamma_i\a_i^pS_n^{n/2}\e \ln(\ln\l_i^\frac{n-2}{2}) +
 O\left(\frac{\e}{\ln\l_i}+\frac{\e\ln\ln\l_i}{(\l_id_i)^{n-2}}+\e^2\ln(\ln\l_i)^2\right) .\nonumber
\end{align}
 \end{lem}
\begin{pf}
Taking into account \eqref{xi}, it is easy to see that $\lambda_i
d(\xi,
 \partial\Omega)
\rightarrow \infty $ is equivalent to $\lambda_i d_i := \lambda_i
d(\xi_i,\partial\Omega)\rightarrow \infty$.\\
By using Lemma \ref{lemme1.55} and \eqref{xi}, we have
\begin{align}\label{majpdelta}
(x-\xi).\n_xP\d_i&=(x-\xi_i).\n_xP\d_i+(\xi_i-\xi).\n_xP\d_i=O(\d_i).
\end{align}
Taylor's expansion with respect to $\e$ yields,
\begin{align}\label{epsilon1}
&\int_{\O}\left[f_\e(\a_i\gamma_iP\delta_i)-f_0(\a_i\gamma_iP\delta_i)\right](x-\xi).\n_xP\d_i\nonumber\\
&=-\gamma_i\e
\int_{\O}\a_i^pP\delta_i^p\ln\ln(e+\a_iP\delta_i)(x-\xi).\n_xP\d_i+O\left(
\e^2\int_{\O}P\delta_i^p[\ln\ln(e+\a_iP\delta_i)]^2|(x-\xi)\n_x P\d_i|\right)\nonumber\\
&=-\gamma_i\e
\int_{\O}\a_i^pP\delta_i^p\ln\ln(e+\a_iP\delta_i)(x-\xi).\n_xP\d_i+
O(\e^2(\ln(\ln\l_i))^2),
\end{align}
where we have used \eqref{majpdelta}, $\int_{\O}\d_i^p=O(1)$ and the
fact that
\begin{equation}\label{majlnln1}0\leq\ln\ln(e+\a_iP\d_i)\leq
c\ln(\ln\l_i).\end{equation}
Taking account of Proposition \ref{p21} $(a)$, the mean value
theorem and \eqref{majlnln1} yield
\begin{equation}\label{devg}
P\delta_i^p\ln\ln(e+\a_iP\delta_i)=\delta_i^p\ln\ln(e+\a_i\delta_i)+O(\delta_i^{p-1}\varphi_i\ln\ln\lambda_i).
\end{equation}
Hence, by using \eqref{devg} and the fact that
$\n_x\d_i=-\frac{\partial \d_i}{\partial \xi_i}$, we get
\begin{align}\label{epsilon2}
&\int_{\O}P\delta_i^p\ln\ln(e+\a_iP\delta_i)(x-\xi).\n_xP\d_i\nonumber
\\&=\int_\O\delta_i^p\ln\ln(e+\a_i\delta_i)(x-\xi).\n_x\d_i
+O\left(\ln\ln\lambda_i\int_\O\delta_i^{p-1}\varphi_i|(x-\xi)\n_x P\d_i|+\ln\ln\lambda_i\int_\O\delta_i^{p}|(x-\xi).\n_{x} \varphi_i|\right)\nonumber\\
&=\int_\O\delta_i^p\ln\ln(e+\a_i\delta_i)(x-\xi).\n_x\d_i
+O\left(\ln\ln\lambda_i\left(\int_{B_i^c}\delta_i^{p+1}+\|\varphi_i\|_{L^\infty(B_i)}\int_{B_i}\delta_i^{p}
+\|\n_x\varphi_i\|_{L^\infty(B_i)}\int_{B_i}\delta_i^{p}|x-\xi_i|\right)\right)
\nonumber\\
&=\int_\O\delta_i^p\ln\ln(e+\a_i\delta_i)(x-\xi).\n_x\d_i
+O\left(\frac{\ln\ln\l_i}{(\l_id_i)^n}+\frac{\ln\ln\l_i}{(\l_id_i)^{n-2}}\right)
\nonumber\\
&=\int_\O\delta_i^p\ln\ln(e+\a_i\delta_i)(x-\xi_i).\n_x\d_i-\sum_{\ell=1}^n((\xi_i)_\ell-(\xi)_\ell)
\int_\O\delta_i^p\ln\ln(e+\a_i\delta_i) \frac{\partial
\d_i}{\partial (\xi_i)_\ell}
+O\left(\frac{\ln\ln\l_i}{(\l_id_i)^{n-2}}\right),
\end{align}
where $B_i$ denotes the ball of center $\xi_i$ and radius $d_i/2$.
Note that we have
 \begin{align}\label{l=1}
 \int_\O\delta_i^p\ln\ln(e+\a_i\delta_i)\frac{\partial  \d_i}{\partial (\xi_i)_\ell}=
\int_{B_i}\delta_i^p\ln\ln(e+\a_i\delta_i)\frac{\partial
\d_i}{\partial (\xi_i)_\ell}
+\int_{B_i^c}\delta_i^p\ln\ln(e+\a_i\delta_i)\frac{\partial
\d_i}{\partial (\xi_i)_\ell}
=O\left(\frac{\l_i\ln\ln\l_i}{(\l_id_i)^{n}}\right)
\end{align}
since the function $\delta_i^p\ln\ln(e+\a_i\delta_i)\frac{\partial
\d_i}{\partial (\xi_i)_\ell}$ is
antisymmetric with respect to $(x-\xi_i)_\ell$ in $B_i$.\\
We point out that we have
\begin{equation}\label{relation}
(\xi_i-x).\n_x\d_i=\frac{n-2}{2}\d_i-\l_i\frac{\partial\d_i}{\partial\l_i}.
\end{equation}
Furthermore, the following integral is needed in the proof of
Proposition \ref{p23} and it was computed in \cite{BFG}. We have
\begin{equation}\label{integrale}
\int_\O\delta_i^{p}\ln\ln(e+\a_i\delta_i)\l_i\frac{\partial\d_i}{\partial\l_i}=\frac{\Gamma_1}{\ln
\l_i}+O\left(\frac{1}{(\ln\l_i)^2}+\frac{\ln\ln\l_i}{(\l_id_i)^n}\right)=
O\left(\frac{1}{\ln\l_i}+\frac{\ln\ln\l_i}{(\l_id_i)^n}\right)
\end{equation}
(for more details see the proof of Lemma 2.9 in \cite{BFG}).\\
 Let $\Omega_{\l_i}:=\l_i(\O-\xi_i)$. We recall that
$$\delta_{(\xi_i,\l_i)}=\l_i^{\frac{n-2}{2}}\delta_{(0,1)}(\l_i(.-\xi_i)).$$
An easy computation shows that for any $U>0$ and $\l$ large enough
we have
\begin{align}\label{expansion}
\ln\ln\left(e+\l^\frac{n-2}{2}U\right)&=\ln\ln\left(\l^\frac{n-2}{2}\right)
+O\left(\frac{\ln(e^{1-\frac{n-2}{2}\ln\l}+U)}{\ln\l}\right)\nonumber\\
&=\ln\ln\left(\l^\frac{n-2}{2}\right)+O\left(\frac{1+\ln(U)}{\ln\l}\right).
\end{align}
 Using \eqref{expansion} with taking $U=\a_i\delta_{(0,1)}$,
together with \eqref{relation} and \eqref{integrale}, we obtain
\begin{align}\label{epsilon3}
\displaystyle
&\sum_{\ell=1}^n\int_\O\delta_i^p\ln\ln(e+\a_i\delta_i)(x_\ell-(\xi_i)_\ell)
\frac{\partial  \d_i}{\partial x_\ell}\nonumber\\
=&-\frac{n-2}{2}\int_\O\delta_i^{p+1}\ln\ln(e+\a_i\delta_i)
+\int_\O\delta_i^{p}\ln\ln(e+\a_i\delta_i)\l_i\frac{\partial\d_i}{\partial\l_i}\nonumber\\
=&-\frac{n-2}{2} \int_{\O_{\l_i}}
\delta_{(0,1)}^{p+1}(y)\ln\ln(e+\l_i^{\frac{n-2}{2}}\a_i\delta_{(0,1)}(y))
dy+O\left(\frac{1}{\ln\l_i}+\frac{\ln\ln\l_i}{(\l_id_i)^n}\right)\nonumber\\
=&-\frac{n-2}{2} \ln\ln\left(\l^\frac{n-2}{2}\right)\int_{\O_{\l_i}}
\delta_{(0,1)}^{p+1}(y) \,dy+O\left(\int_{\O_{\l_i}}
\delta_{(0,1)}^{p+1}(y)\frac{1+\ln(\a_i\delta_{(0,1)}(y))}{\ln\l}\,dy\right)
    +O\left(\frac{1}{\ln\l_i}+\frac{\ln\ln\l_i}{(\l_id_i)^n}\right)\nonumber\\
=&-\frac{n-2}{2}\ln\ln\left(\l^\frac{n-2}{2}\right)\int_{\mathbb{R}^n}
\delta_{(0,1)}^{p+1}(y) dy
+O\left(\ln\ln\l_i\int_{B(0,\l_id_i)^c}\d_{(0,1)}^{p+1}\right)+
O\left(\frac{1}{\ln\l_i}+\frac{\ln\ln\l_i}{(\l_id_i)^n}\right)
 \nonumber\\
=&-\frac{(n-2)}{2}\ln\ln\left(\l^\frac{n-2}{2}\right)S_n^{n/2}
+O\left(\frac{1}{\ln\l_i}+\frac{\ln\ln\l_i}{(\l_id_i)^n}\right).
\end{align}
Thus, the desired result follows from \eqref{epsilon1},
\eqref{epsilon2}, \eqref{l=1} and \eqref{epsilon3}.
\end{pf}
\begin{lem}\label{lem25}Let $\xi\in\Omega$
 and $\lambda_i,\lambda_j > 0$ be such that $\min(\lambda_i,\lambda_j) d(\xi,
 \partial\Omega)$ is very large. We denote
\begin{equation}\label{xixj}
 \xi_l:=\xi+\frac{\sigma_l}{\l_l}\hbox{ for }l=i,j
\end{equation}
where $\sigma_l$ is in compact set of $\mathbb{R}^n$. Moreover, we
suppose that $\e_{ij}$ is very small. Then, there hold
\begin{align}\label{int1}
\int_{\O}&P\delta_i^p(x-\xi).\n_xP\d_i=
-\frac{n-2}{2}S_n^{\frac{n}{2}}+n\overline{c}_1\frac{H(\xi_i,\xi_i)}{\l_i^{n-2}}
+\left\{\begin{array}{ll} O\left(\frac{1}{(\l_id_i)^2}\right),\quad\mbox{ if }\, n=3;\\
O\left(\frac{ \ln(\l_id_i)}{(\l_id_i)^n}\right),\quad\mbox{ if }\,
n\geq 4,
\end{array}\right.
\end{align}
\begin{align} \label{int2}\int_{\O}P\delta_j^p(x-\xi).\n_xP\d_i=& -\frac{ n-2
}{2}\overline{c}_1\left(\e_{ij}-
 \frac{H(\xi_i,\xi_j)}{(\l_i\lambda_j)^{(n-2)/2}}\right)-\overline{c}_1\left(\l_j\frac{\partial\e_{ij}}{\partial\l_j}+\frac{n-2}{2}
 \frac{H(\xi_i,\xi_j)}{(\l_i\lambda_j)^{(n-2)/2}}\right)\notag\\
 &+ \overline{c}_1\sum_{\ell=1}^n((\xi_j)_\ell-(\xi)_\ell) \left(\frac{\partial\e_{ij}}{\partial
(\xi_j)_\ell}-
 \frac{1}{(\l_i\lambda_j)^{(n-2)/2}}\frac{\partial H}{\partial (b)_\ell}(\xi_i,\xi_j)\right)\notag\\
 &+\left\{
                             \begin{array}{ll}
                              O\left(\e^2_{ij}\ln(\e_{ij}^{-1})^\frac{2}{3}+\sum_{l=i,j}\frac{1}{(\l_l d_l)^2} \right), & \hbox{if }n=3; \\
                               O\left(\e_{ij}^\frac{n}{n-2}
\ln (\e_{ij}^{-1})+\sum_{l=i,j} \frac{\ln (\l_l d_l)}{(\l_l
d_l)^n}\right), & \hbox{if }n\geq 4,
                             \end{array}
                           \right.
\end{align}
\begin{align}\label{int3}
p\int_{\O}P\delta_i^{p-1}P\delta_j(x-\xi).\n_xP\d_i=&-\frac{ n+2
}{2}\overline{c}_1\left(\e_{ij}-
 \frac{H(\xi_i,\xi_j)}{(\l_i\lambda_j)^{(n-2)/2}}\right)+\overline{c}_1\left(\l_i\frac{\partial\e_{ij}}{\partial\l_i}+\frac{n-2}{2}
 \frac{H(\xi_i,\xi_j)}{(\l_i\lambda_j)^{(n-2)/2}}\right)\notag\\
 &- \overline{c}_1\sum_{\ell=1}^n((\xi_i)_\ell-(\xi)_\ell) \left(\frac{\partial\e_{ij}}{\partial
(\xi_i)_\ell}-
 \frac{1}{(\l_i\lambda_j)^{(n-2)/2}}\frac{\partial H}{\partial (a)_\ell}(\xi_i,\xi_j)\right)\notag\\&+\left\{
                             \begin{array}{ll}
                              O\left(\e^2_{ij}\ln(\e_{ij}^{-1})^\frac{2}{3}+\sum_{l=i,j}\frac{1}{(\l_l d_l)^2} \right), & \hbox{if }n=3; \\
                               O\left(\e_{ij}^\frac{n}{n-2}
\ln (\e_{ij}^{-1})+\sum_{l=i,j} \frac{\ln (\l_l d_l)}{(\l_l
d_l)^n}\right), & \hbox{if }n\geq 4,
                             \end{array}
                           \right.
\end{align}
where $\frac{\partial H}{\partial (a)_\ell}$ and $\frac{\partial
H}{\partial (b)_\ell}$ denote the partial derivatives of $H$ with
respect to the $\ell-$th components of the first variable and the
second variable.
\end{lem}
\begin{pf}
Taking into account \eqref{xixj}, it is easy to see that
$min(\lambda_i,\l_j )d(\xi,\partial \O) \rightarrow \infty $ is
equivalent to $\lambda_l d_l := \lambda_l
d(\xi_l,\partial\Omega)\rightarrow \infty$ for $l=i,j$.\\
We remark that
\begin{equation}\label{majphi}
(x-\xi).\n_x\varphi_i=O(\d_i)
\end{equation}
since \eqref{majpdelta} holds and we have $(x-\xi).\n_x\d_i=O(\d_i)$
which derive from \eqref{relation}, \eqref{xixj} and the fact that  $\n_x\d_i=-\frac{\partial\d_i}{\partial \xi_i}$. \\
Let us compute the first integral. Integrating by parts and using
\eqref{q1}, \eqref{majpdelta}, \eqref{majphi} and Holder's
inequality, we get
\begin{align}\label{int11}
\int_{\O}P\delta_i^p(x-\xi).\n_xP\d_i=&\sum_{\ell=1}^n\int_{\O}
\delta_i^p(x_\ell-(\xi)_\ell)\frac{\partial P\d_i}{\partial
x_\ell}-p
\int_{\O}\delta_i^{p-1}\varphi_i(x-\xi).\n_x  \d_i\notag\\
&+O\left(
\int_{B_i}\delta_i^{p-1}\varphi_i^2+\int_{B_i}\delta_i^{p-1}\varphi_i
|(x - \xi)\n_x\varphi_i|+
\int_{B_i^c}\delta_i^{p+1}\right) \notag\\
=&-\sum_{\ell=1}^n\int_{\O}\frac{\partial}{\partial x_\ell}
\left((x_\ell-(\xi)_\ell)\delta_i^p\right) P\d_i-p
\int_{\O}\delta_i^{p-1}\varphi_i(x-\xi_i).\n_x
\d_i-p\int_{\O}\delta_i^{p-1}\varphi_i(\xi_i-\xi).\n_x
\d_i\notag\\&+O\left(
\int_{B_i}\delta_i^{p-1}\varphi_i^2+\int_{B_i}\delta_i^{p-1}\varphi_i
|(x - \xi)\n_x\varphi_i|+
\frac{1}{(\l_id_i)^n}\right)\notag\\=&-n\int_{\O}\delta_i^pP\d_i
-p\int_{\O}\delta_i^{p-1}P\d_i(x-\xi_i).\n_x
\d_i-p\int_{\O}\delta_i^{p-1}P\d_i(\xi_i-\xi).\n_x
\d_i\nonumber\\&-p \int_{\O}\delta_i^{p-1}\varphi_i(x-\xi_i).\n_x
\d_i-p\int_{\O}\delta_i^{p-1}\varphi_i(\xi_i-\xi).\n_x \d_i+\left\{
                             \begin{array}{ll}
                              O\left(\frac{1}{(\l_i d_i)^2} \right), & \hbox{if }n=3; \\
                               O\left(\frac{\ln(\l_i d_i)}
 {(\lambda_i d_i)^{n}}\right), & \hbox{if }n\geq 4.
                             \end{array}
                           \right.
\end{align}
Using \eqref{int11}, \eqref{relation}, \eqref{xixj},
\cite[Lemma2.4]{BFG} and the fact that
$\n_x\d_i=-\frac{\partial\d_i}{\partial \xi_i}$, we obtain
\begin{align*}
\int_{\O}P\delta_i^p(x-\xi).\n_xP\d_i
=&-\frac{n-2}{2}\int_{\O}\delta_i^pP\d_i-
p\int_{\O}\delta_i^{p-1}\l_i\frac{\partial\d_i}{\partial\l_i} P\d_i
+p\int_{\O}\delta_i^{p-1}P\d_i(\xi_i-\xi).\frac{\partial
\d_i}{\partial\xi_i}\nonumber\\
&+\frac{n+2}{2}\int_\O\d_i^p\varphi_i-p\int_\O\d_i^{p-1}\varphi_i\l_i\frac{\partial\d_i}
{\partial\l_i}+p\int_{\O}\delta_i^{p-1}\varphi_i(\xi_i-\xi).\frac{\partial
\d_i}{\partial\xi_i}\notag\\&+\left\{
                             \begin{array}{ll}
                              O\left(\frac{1}{(\l_i d_i)^2} \right), & \hbox{if }n=3; \\
                               O\left(\frac{\ln(\l_i d_i)}
 {(\lambda_i d_i)^{n}}\right), & \hbox{if }n\geq 4,
                             \end{array}
                           \right.
                           \notag\\=&-\frac{n-2}{2}\int_{\O}\delta_i^{p+1}+n\int_\O\d_i^p\varphi_i
-2p\int_{\mathbb{R}^n}\delta_i^{p}\l_i\frac{\partial\d_i}{\partial\l_i}
+2p\int_{B_i}\delta_i^{p}(\xi_i-\xi).\frac{\partial
\d_i}{\partial\xi_i}\nonumber\\
& +\left\{
                             \begin{array}{ll}
                              O\left(\frac{1}{(\l_i d_i)^2} \right), & \hbox{if }n=3; \\
                               O\left(\frac{\ln(\l_i d_i)}
 {(\lambda_i d_i)^{n}}\right), & \hbox{if }n\geq 4,
                             \end{array}
                           \right.\notag\\
=&-\frac{n-2}{2}S_n^{n/2}+n\overline{c}_1\frac{H(\xi_i,\xi_i)}{\l_i^{n-2}}+\left\{
                             \begin{array}{ll}
                              O\left(\frac{1}{(\l_i d_i)^2} \right), & \hbox{if }n=3; \\
                               O\left(\frac{\ln(\l_i d_i)}
 {(\lambda_i d_i)^{n}}\right), & \hbox{if }n\geq 4,
                             \end{array}
                           \right.
\end{align*}
since the two last integrals are equal to $0$. Thus \eqref{int1} follows.\\
We focus now on Eq. \eqref{int2}. Using an integration by parts,
\eqref{q1}, \eqref{relation}, \eqref{xixj}, the fact that
$\n_x\d_i=-\frac{\partial\d_i}{\partial \xi_i}$ and \cite[Lemma
2.4]{BFG}, we have
\begin{align*}
\int_{\O}&P\delta_j^p(x-\xi).\n_xP\d_i\notag\\=&\sum_{\ell=1}^n\int_{\O}\delta_j^p(x_\ell-(\xi)_\ell))
\frac{\partial P\d_i}{\partial x_\ell} +O\left(\int_{\O}\delta_j^{p-1}\varphi_j \d_i \right)\notag\\
=&-\sum_{\ell=1}^n\int_{\O}\frac{\partial}{\partial
x_\ell}(\delta_j^p(x_\ell-(\xi)_\ell))  P\d_i
+O\left(\int_{\O}\delta_j^{p-1}\varphi_j \d_i \right)\notag\\
=&-n\int_{\O}\d_j^pP\d_i-p\int_{\O}\d_j^{p-1}P\d_i(x-\xi_j).\n_x\d_i
-p\sum_{\ell=1}^n((\xi_j)_\ell-(\xi)_\ell)\int_{\O}\d_j^{p-1}\frac{\partial
\delta_j}{\partial
x_\ell}P\d_i+O\left(\int_{\O}\delta_j^{p-1}\varphi_j \d_i
\right)\notag\\=&-\frac{ n-2 }{2}\int_{\O}\delta_j^pP\d_i-
p\int_{\O}\delta_j^{p-1}\l_j\frac{\partial \delta_j }{\partial \l_j}
P\d_i+p\sum_{\ell=1}^n\frac{(\sigma_j)_\ell}{\l_\ell}\int_{\O}\delta_j^{p-1}
\frac{1}{\l_j}\frac{\partial \delta_j }{\partial
(\xi_j)_\ell}P\d_i+O\left(\int_{\O}\delta_j^{p-1}\varphi_j \d_i
\right)\notag\\=& -\frac{ n-2 }{2}\overline{c}_1\left(\e_{ij}-
 \frac{H(\xi_i,\xi_j)}{(\l_i\lambda_j)^{(n-2)/2}}\right)-\overline{c}_1\left(\l_j\frac{\partial\e_{ij}}{\partial\l_j}+
 \frac{n-2}{2}
 \frac{H(\xi_i,\xi_j)}{(\l_i\lambda_j)^{(n-2)/2}}\right)\notag\\
 &+ \overline{c}_1\sum_{\ell=1}^n\frac{(\sigma_j)_\ell}{\l_\ell} \left(\frac{\partial\e_{ij}}{\partial
(\xi_j)_\ell}-
 \frac{1}{(\l_i\lambda_j)^{(n-2)/2}}\frac{\partial H}{\partial (b)_\ell}(\xi_i,\xi_j)\right)+\left\{
                             \begin{array}{ll}
                              O\left(\e^2_{ij}\ln(\e_{ij}^{-1})^\frac{2}{3}+\sum_{l=i,j}\frac{1}{(\l_l d_l)^2} \right), & \hbox{if }n=3; \\
                               O\left(\e_{ij}^\frac{n}{n-2}
\ln (\e_{ij}^{-1})+\sum_{l=i,j} \frac{\ln (\l_l d_l)}{(\l_l
d_l)^n}\right), & \hbox{if }n\geq 4.
                             \end{array}
                           \right.
\end{align*}
Finally, we prove \eqref{int3}. Integrating by parts, we have
 \begin{align}\label{int30}
p\int_{\O}P\delta_i^{p-1}P\delta_j(x-\xi).\n_xP\d_i&=\sum_{\ell=1}^n\int_{\O}
P\delta_j(x_\ell-(\xi)_\ell)\frac{\partial (P\d_i^p)}{\partial
x_k}\nonumber\\&=- \sum_{\ell=1}^n\int_{\O} P\d_i^p\frac{\partial
}{\partial
x_\ell}((x_\ell-(\xi)_\ell)P\delta_j)\notag\\&=-n\int_{\O} P\d_i^p
P\delta_j
-\sum_{\ell=1}^n\int_{\O}P\delta_i^p(x_\ell-(\xi)_\ell)\frac{\partial
P\d_j}{\partial x_\ell}\notag\\&=-n\int_{\O} P\d_i^p
P\delta_j-\int_{\O}P\delta_i^p(x-\xi).\n _xP\d_j.
\end{align}
We know that
\begin{align}\label{int31}
\int_{\O} P\d_i^p P\delta_j= \overline{c}_1\left(\e_{ij}-
 \frac{H(a_i,a_j)}{(\l_i\lambda_j)^{(n-2)/2}}\right)+\left\{
                             \begin{array}{ll}
                              O\left(\e^2_{ij}\ln(\e_{ij}^{-1})^\frac{2}{3}+\frac{1}{(\l_j d_j)^2} \right), & \hbox{if }n=3; \\
                               O\left(\e_{ij}^\frac{n}{n-2}
\ln (\e_{ij}^{-1})+\sum_{l=i,j} \frac{\ln (\l_l d_l)}{(\l_l
d_l)^n}\right), & \hbox{if }n\geq 4,
                             \end{array}
                           \right.
\end{align}
see for instance \cite[Lemma 2.5]{BFG}. The result of \eqref{int3}
follows from
\eqref{int30}, \eqref{int2} and \eqref{int31}.\\
 Thus the proof of the lemma is thereby completed.
\end{pf}

Our result, concerning the expansion of $I_\e$, can be stated as
follows.
\begin{pro}\label{p23'}
Let $n\geq 3$ and $u=\sum_{j=1}^k\a_j\gamma_j P\d_j+v\in V(k,\eta)$
such that $v\in E_{(\xi,\l)}$. Assume that
\begin{equation*}
 \xi_i=\xi+\frac{\sigma_i}{\l_i}\hbox{ for }i=1,\ldots,k,
\end{equation*}
where $\xi$ is some point in $\Omega$ and $\sigma_i$ is in compact
set of $\mathbb{R}^n$ for each $i= 1,...,k$.\\
We have the following expansion
\begin{align*}
 I_\e(u)
 =&\sum_{i=1}^k\Big(\frac{\a_i^2}{2}-\frac{\a_i^{p+1}}{p+1}\Big)S_n^{\frac{n}{2}}-\ov c_1\sum_{i=1}^k
 \Big(\frac{\a_i^2}{2}-\a_i^{p+1}\Big)\frac{H(\xi_i,\xi_i)}{\l_i^{n-2}}
 +\frac{\e (n-2)S_n^{\frac{n}{2}}}{2n}\sum_{i=1}^k\a_i^{p+1}
 \ln(\ln\l_i^\frac{n-2}{2})\notag\\&
 -\ov c_1\sum_{i=1}^k\sum_{j\neq i}\gamma_i\gamma_j
 \Big(\frac{n+2}{2n}\a_i^{p}\a_j+\frac{n-2}{2n}\a_j^{p}\a_i-\frac{\a_i\a_j}{2}\Big)
 \Big(\e_{ij}-\frac{H(\xi_i,\xi_j)}{(\l_i\l_j)^{(n-2)/2}}\Big)\notag\\
 &+\frac{\overline{c}_1}{n}\sum_{i=1}^k \sum_{j\neq i} \gamma_i\gamma_j\left(\a_i^p\a_j\l_i
 \frac{\partial\e_{ij}}{\partial\l_i}-\a_j^p\a_i\l_j\frac{\partial\e_{ij}}{\partial\l_j}
 +\frac{n-2}{2}\big(\a_i^p\a_j-\a_j^p\a_i\big)
 \frac{H(\xi_i,\xi_j)}{(\l_i\lambda_j)^{(n-2)/2}}\right)\notag\\
 &+\frac{\overline{c}_1}{n}\sum_{i=1}^k \sum_{j\neq i} \a_i\a_j^p\gamma_i\gamma_j
 \sum_{\ell=1}^n((\xi_j)_\ell-(\xi)_\ell) \left(\frac{\partial\e_{ij}}{\partial
(\xi_j)_\ell}-
 \frac{1}{(\l_i\lambda_j)^{(n-2)/2}}\frac{\partial H}{\partial (b)_\ell}(\xi_i,\xi_j)\right)\notag\\
 &-\frac{\overline{c}_1}{n}\sum_{i=1}^k \sum_{j\neq i}\a_i ^p\a_j\gamma_i\gamma_j
 \sum_{\ell=1}^n((\xi_i)_\ell-(\xi)_\ell) \left(\frac{\partial\e_{ij}}{\partial
(\xi_i)_\ell}-
 \frac{1}{(\l_i\lambda_j)^{(n-2)/2}}\frac{\partial H}{\partial (a)_\ell}(\xi_i,\xi_j)\right)+\frac{1}{2}Q(v,v)-f(v)
 \\&+o(\|v\|^2)+O(\e^2\ln(\ln \l_{\max})^2)+\left\{ \begin{array}{ll}
                              O\left(\sum_{i=1}^k\frac{1}{(\l_i d_i)^{2}}+\sum_{j\neq i}\e_{ij}^2\ln(
\e_{ij}^{-1})^\frac{2}{3} \right), & \hbox{if }n=3; \\
                               O\left(\sum_{i=1}^k\frac{\ln(\l_i d_i)}{(\l_i
d_i)^{n}}+\sum_{j\neq i}\e_{ij}^\frac{n}{n-2}\ln
(\e_{ij}^{-1})\right), & \hbox{if }n\geq4.
                             \end{array}
                           \right.
\end{align*}
where $S_n$ is the best Sobolev constant, $\ov c_1$ is introduced in
Proposition \ref{p23}
and
\begin{eqnarray*}\begin{cases}
 Q(v,v)=\displaystyle\int_{\Omega}|\nabla
v |^2 -\frac{n+2}{n-2} \sum_{i=1}^k \displaystyle\int_{\Omega}
\d_i^{p-1}v^2,\\
 f(v) =\displaystyle\int_{\Omega}f_\e(\sum_{i=1}^k\a_i\gamma_i P\d_i)v.\end{cases}
 \end{eqnarray*}
 \end{pro}
\begin{pf}
Recall  that $\underline{u}:= \sum_{j=1}^k\a_j\gamma_jP\delta_j$.
Integrating by parts and using \eqref{nonpower} and \eqref{g0}, we
have
\begin{align}\label{dev}
I_{\varepsilon}( u)=&\frac{1}{2}<u,u>-\int_{\O}F_\e(u)\nonumber\\
=&\frac{1}{2}<\underline{u},\underline{u}>
+\frac{1}{2}<v,v>-\int_{\O}\left[F_\e(u)-F_\e(\underline{u})-F'_\e(\underline{u})v-
\frac{1}{2}F_\e''(\underline{u})v^2\right]\nonumber\\
&-\int_{\O}F_\e(\underline{u})-\int_{\O}F'_\e(\underline{u})v-
\frac{1}{2}\int_{\O}F_\e''(\underline{u})v^2\nonumber\\=&I_{\varepsilon}(\underline{u})
-\int_{\O}F'_\e(\underline{u})v+\frac{1}{2}<v,v>-
\frac{1}{2}\int_{\O}F_0''(\underline{u})v^2\nonumber\\
&-
\frac{1}{2}\int_{\O}(F_\e''(\underline{u})-F_0''(\underline{u}))v^2-\int_{\O}\left[F_\e(u)
-F_\e(\underline{u})-F'_\e(\underline{u})v-
\frac{1}{2}F_\e''(\underline{u})v^2\right]\nonumber\\
=&I_{\varepsilon}( \underline{u})
-\displaystyle\int_{\Omega}f_\e(\underline{u})v+\frac{1}{2}\left\{\displaystyle\int_{\Omega}|\nabla
v |^2 -\frac{n+2}{n-2} \sum_{i=1}^k \displaystyle\int_{\Omega}
\d_i^{p-1}v^2\right\}+R_{\eps,\a,a,\l}(v),
\end{align}
where \begin{align}R_{\eps,\a,\l,a}:=&-
\frac{1}{2}\int_{\O}(F_\e''(\underline{u})-F_0''(\underline{u}))v^2-\int_{\O}\left[F_\e(u)-F_\e(\underline{u})-F'_\e(\underline{u})v-
\frac{1}{2}F_\e''(\underline{u})v^2\right]\nonumber\\
&-\frac{n+2}{2(n-2)}\displaystyle\int_{\Omega}\left(|\underline{u}|^{p-1}-\sum_{i=1}^k
\d_i^{p-1}\right)v^2\nonumber\\=&o(\|v\|^2).
\end{align}
Using \eqref{g0} and integrating by parts, we have
\begin{align}I_\e(\underline{u})&=\frac{1}{2}\int_{\O}|\n\underline{u}|^2-\int_{\O}F_\e(\underline{u})\nonumber
\\&=\frac{1}{2}\int_{\O}|\n
\underline{u}|^2+\frac{1}{n}\sum_{\ell=1}^n\int_{\O}f_\e(\underline{u})(x_\ell-(\xi)_\ell)
\frac{\partial\underline{u}}{\partial
x_\ell}\notag\\&=\frac{1}{2}\int_{\O}|\n
\underline{u}|^2+\frac{1}{n}\sum_{i=1}^k\a_i\gamma_i\int_{\O}f_\e(\underline{u})\,
(x-\xi). \n_x P\d_i.
\end{align}
On one hand, using \cite[Lemma 2.2]{BFG} and \cite[Lemma 2.4]{BFG},
we get
\begin{align}
\int_{\O}|\n \underline{u}|^2=&\sum_{i=1}^k\a_i^2\|P\d_i\|^2
+\sum_{i\neq
j}\a_i\a_j\gamma_i\gamma_j\left<P\d_i,P\d_j\right>\nonumber\\
=&\displaystyle\sum_{i=1}^k\a_i^2\left(S_n^\frac{n}{2}-\ov
c_1\frac{H(\xi_i,\xi_i)}{\l_i^{n-2}}\right)+\sum_{i\neq
j}\a_i\a_j\gamma_i\gamma_j\ov
c_1\left(\e_{ij}-\frac{H(\xi_i,\xi_j)}{(\l_i\l_j)^\frac{n-2}{2}}\right)\nonumber\\&
+O\left(\sum_{i=1}^k\frac{\ln(\l_i d_i)}{(\l_i d_i)^{n}}+\sum_{j\neq
i}\e_{ij}^\frac{n}{n-2}\ln (\e_{ij}^{-1})\right).
\end{align}
On the other hand, we write
\begin{align}
\int_{\O}& f_\e(\underline{u})(x-\xi).\n_x
P\d_i\notag\\&=\int_{\O}\biggr[f_\e(\underline{u})-f_\e(\a_i\gamma_iP\delta_i)-f_\e(\sum_{j\neq
i}\a_j\gamma_jP\delta_j)-f'_\e(\a_i\gamma_iP\delta_i)\sum_{j\neq
i}\a_j\gamma_jP\delta_j\biggr](x-\xi).\n_x P\d_i\nonumber\\
&+\int_{\O}f_\e(\a_i\gamma_iP\delta_i)(x-\xi).\n_x
P\d_i+\int_{\O}f_\e(\sum_{j\neq i}\a_j\gamma_jP\delta_j)(x-\xi).\n_x
P\d_i\nonumber\\&+\int_{\O}f'_\e(\a_i\gamma_iP\delta_i)(\sum_{j\neq
i}\a_j\gamma_jP\delta_j)(x-\xi).\n_x P\d_i\nonumber\\
&=:B_1+B_2+B_3+B_4.
\end{align}
We start with the last integral. Using \eqref{f'}, \eqref{majdd},
\eqref{majpdelta}, \eqref{majlnln1} and \eqref{e_ij}, we get
\begin{align}
B_4&:=\displaystyle\int_{\O}f'_\e(\a_i\gamma_iP\delta_i)(\sum_{j\neq
i}\a_j\gamma_jP\delta_j)(x-\xi).\n_x P\d_i\nonumber\\
&=\int_{\O}f'_0(\a_i\gamma_iP\delta_i)(\sum_{j\neq
i}\a_j\gamma_jP\delta_j)(x-\xi).\n_x
P\d_i+\int_{\O}\left[f'_\e(\a_i\gamma_iP\delta_i)-f'_0(\a_i\gamma_iP\delta_i)\right]
(\sum_{j\neq
i}\a_j\gamma_jP\delta_j)(x-\xi).\n_x P\d_i\nonumber\\
&=p\int_{\O}(\a_iP\delta_i)^{p-1}(\sum_{j\neq
i}\a_j\gamma_jP\delta_j)(x-\xi).\n_x P\d_i+O\left(\e\ln(\ln
\l_i)\sum_{j\neq
i}\int_{\Omega}\delta_i^p\delta_j\right)\nonumber\\
&=\sum_{j\neq
i}\gamma_j\a_i^{p-1}\a_jp\int_{\O}P\delta_i^{p-1}P\delta_j(x-\xi).\n_x
P\d_i+O\left(\e\ln(\ln \l_i)\sum_{j\neq i}\e_{ij}\right).
\end{align}
Using the first claim of Lemma \ref{lemmaB1}, \eqref{majdd},
\eqref{majpdelta} and the fact that $\ln\ln(e+|\sum_{j\neq
i}\a_j\gamma_jP\d_j |)=O( \ln\ln\l_{\max})$ where
$\l_{\max}:=\max(\l_1,\ldots,\l_m)$, we have
\begin{align}
B_3&:=\displaystyle\int_{\O}f_\e(\sum_{j\neq
i}\a_j\gamma_jP\delta_j)(x-\xi).\n_x P\d_i\nonumber\\
&=\int_{\O}f_0(\sum_{j\neq i}\a_j\gamma_jP\delta_j)(x-\xi).\n_x
P\d_i+\int_{\O}\left[f_\e(\sum_{j\neq
i}\a_j\gamma_jP\delta_j)-f_0(\sum_{j\neq i}\a_j\gamma_jP\delta_j)\right](x-\xi).\n_x P\d_i\nonumber\\
&=\int_{\O}f_0(\sum_{j\neq
i}\a_j\gamma_jP\delta_j)(x-\xi).\n_x P\d_i+O\left(\e\ln(\ln\l_{\max})\sum_{j\neq i}\int_{\Omega}\delta_j^p\delta_i\right)\nonumber\\
&=\sum_{j\neq i}\a_j^p\gamma_j\int_{\O}P\delta_j^p(x-\xi).\n_x
P\d_i+O\left(\sum_{l\neq j;l,j\neq
i}\int_{\O}P\delta_j^{p-1}\inf(P\delta_j,P\delta_l)\delta_i+\e\ln(\ln
\l_{\max})\sum_{j\neq i}\e_{ij}\right)\nonumber\\
&=\sum_{j\neq i}\a_j^p\gamma_j\int_{\O}P\delta_j^p(x-\xi).\n_x
P\d_i+O\left(\sum_{l\neq j}\e_{lj}^\frac{n}{n-2}\ln \e_{lj}^{-1}
+\e^2\ln(\ln \l_{\max})^2+\sum_{j\neq i}\e_{ij}^2\right).
\end{align}
Now, Lemma \ref{lemmaA1} implies that
\begin{align}
B_2&:=\displaystyle\int_{\O}f_\e(\a_i\gamma_iP\delta_i)(x-\xi).\n_x P\d_i\nonumber\\
&=\int_{\O}\left[f_\e(\a_i\gamma_iP\delta_i)-f_0(\a_i\gamma_iP\delta_i)\right](x-\xi).\n_x
P\d_i+
\int_{\O}f_0(\a_i\gamma_iP\delta_i)(x-\xi).\n_x P\d_i\nonumber\\
&=
 \frac{n-2}{2}\gamma_i\a_i^pS_n^{n/2}\e \ln(\ln\l_i^\frac{n-2}{2}) +\gamma_i\a_i^p\int_{\O}P\d_i^p(x-\xi).\n_x
 P\d_i+
 O\left(\frac{\e}{\ln\l_i}+\frac{\e\ln\ln\l_i}{(\l_id_i)^{n-2}}+\e^2\ln(\ln\l_i)^2\right) .
\end{align}
In the sequel, we compute $B_1$. Let $\O_1:=\{x:\ |\sum_{j\neq
i}\a_j\gamma_jP\delta_j(x)|\leq \frac{1}{2}\a_iP\delta_i(x)\}$.
\begin{align}B_1&:=\int_{\O}\left[f_\e(\underline{u})-f_\e(\a_i\gamma_iP\delta_i)-f_\e(\sum_{j\neq
i}\a_j\gamma_jP\delta_j)-f'_\e(\a_i\gamma_iP\delta_i)\sum_{j\neq
i}\a_j\gamma_jP\delta_j\right](x-\xi).\n_x P\d_i\nonumber \\
&=\int_{\O_1}\ldots +\int_{\O\setminus\O_1}\ldots
\nonumber \\
&=B_{11}+B_{12}.
\end{align}
Observe that, in $\O_1$, it holds $1/2 \a_iP\d_i\leq
|\a_i\gamma_iP\delta_i+\theta\sum_{j\neq
i}\a_j\gamma_jP\delta_j|\leq 3/2 \a_iP\d_i$ for each $\theta\in
(0,1)$. By using the mean value theorem and \eqref{B''}, we have
\begin{align}
|B_{11}|
&\leq\displaystyle\int_{\O_1}\left|f_\e(\underline{u})-f_\e(\a_i\gamma_iP\delta_i)
-f'_\e(\a_i\gamma_iP\delta_i)\sum_{j\neq
i}\a_j\gamma_jP\delta_j\right|\delta_i+\int_{\O_1}|f_\e(\sum_{j\neq
i}\a_j\gamma_jP\delta_j)|\delta_i\nonumber\\
&\leq
c\displaystyle\int_{\O_1}\left|f''_\e(\a_i\gamma_iP\delta_i+\theta\sum_{j\neq
i}\a_j\gamma_jP\delta_j)\right|\left(\sum_{j\neq
i}\a_jP\delta_j\right)^2\delta_i+\int_{\O_1}(\sum_{j\neq
i}\a_jP\delta_j)^p\delta_i, \hbox{ for some }\theta\in(0,1)\nonumber\\
&\leq c\displaystyle\int_{\O_1}\left(\sum_{j\neq
i}\a_jP\delta_j\right)^2\delta_i^{p-1}+c\sum_{j\neq
i}\int_{\O_1}(\delta_j \delta_i)^\frac{n}{n-2}\nonumber\\
&\leq  c \left\{ \begin{array}{ll}
                              O\left(\sum_{j\neq i}\e_{ij}^2\ln(
\e_{ij}^{-1})^\frac{2}{3} \right), & \hbox{if }n=3; \\
                               O\left(\sum_{j\neq i}\e_{ij}^\frac{n}{n-2}\ln (\e_{ij}^{-1})\right), & \hbox{if }n\geq4.
                             \end{array}
                           \right.
\end{align}
Using again the mean value theorem, \eqref{majpdelta}, \eqref{B.1}
and the fact that $|\a_j-1|<\eta $ for each $j$, we obtain
\begin{align}\label{d2f2}
|B_{12}|
&\leq\displaystyle\int_{\O\setminus\O_1}\left(|f_\e(\a_i\gamma_iP\delta_i)|+|f'_\e(\a_i\gamma_iP\delta_i)
||\sum_{j\neq
i}\a_j\gamma_jP\delta_j|\right)\delta_i+\int_{\O\setminus\O_1}|f_\e(\underline{u})-f_\e(\sum_{j\neq
i}\a_j\gamma_jP\delta_j)|\delta_i\nonumber\\
&\leq
c\displaystyle\int_{\O\setminus\O_1}\a_i^pP\delta_i^{p}\delta_i+\a_i^{p-1}P\delta_i^{p-1}|\sum_{j\neq
i}\a_j\gamma_jP\delta_j|\delta_i+\int_{\O\setminus\O_1}|f'_\e(\sum_{j\neq
i}\a_j\gamma_jP\delta_j+\theta\a_i\gamma_iP\delta_i)|\a_iP\delta_i\delta_i, \hbox{ for some }\theta\in(0,1)\nonumber\\
&\leq c \sum_{j\neq i}\e_{ij}^\frac{n}{n-2}\ln
(\e_{ij}^{-1})+c\int_{\O\setminus\O_1}|\sum_{j\neq
i}\a_j\gamma_jP\delta_j|^{p-1}P\delta_i\delta_i\nonumber\\
&\leq c \left\{ \begin{array}{ll}
                              O\left(\sum_{j\neq i}\e_{ij}^2\ln(
\e_{ij}^{-1})^\frac{2}{3} \right), & \hbox{if }n=3; \\
                               O\left(\sum_{j\neq i}\e_{ij}^\frac{n}{n-2}\ln (\e_{ij}^{-1})\right), & \hbox{if }n\geq4.
                             \end{array}
                           \right.
\end{align}
 Combining \eqref{dev}-\eqref{d2f2} and Lemma \ref{lem25}, the
result of Proposition \ref{p23'} follows.
\end{pf}
\section{Asymptotic behavior}
In this section, we investigate the asymptotic profile of a family
of sign changing bubble tower solutions $u_\e$ blowing up in the
interior of the domain $\Omega$. Our aim is to look for suitable
conditions on the parameters $\l_i$ and $\xi_i$ to construct such
solution in the next sections.\\ For simplicity, we will assume that
$u_\e$ is a solution of $(P_\e)$ having the following form
\begin{equation}\label{h:1}
u_\e=\gamma_1 P\d_{(\xi_{1 ,\e},\l_{1,\e})}+\gamma_2 P\d_{(\xi_{2
,\e},\l_{2,\e})}+v_\e
\end{equation} where
\begin{equation}\label{gamma}\gamma_1=-\gamma_2,\, |\gamma_1|=1,
\end{equation} \begin{equation}\label{di}
d(\xi_i,\partial \Omega)>c\hbox{ for }i=1,2
\end{equation} and $v_\e$
satisfies
\begin{equation}\label{vpart}\|v_\e\|=\min_{\l_1,\l_2,\xi_1,\xi_2}
 \|u_\e-\gamma_1 P\d_{(\xi_{1 ,\e},\l_{1,\e})}-\gamma_2 P\d_{(\xi_{2 ,\e},\l_{2,\e})}\|.
 \end{equation}
Therefore we get $\left<v_\e,\frac{\partial P\d_{i}}{\partial
\l_i}\right>= \left<v_\e,\frac{\partial P\d_{i}}{\partial
(\xi_i)_j}\right>=0\;\forall \; 1\leq j\leq n,\,\forall \,1\leq
i\leq 2$ and by multiplying the equation $-\Delta u_\e=f_\e(u_\e)$
by $v_\e$ and integrating by parts the following estimate holds
\begin{align}\label{vestimate}
\|v_\e\| \leq c\sum_{j=1}^2 \e\ln(\ln \l_j)+c \left\{
\begin{array}{ccc}
\displaystyle\sum_{j=1}^2\frac{1}{\l_j^{n-2}}+\e_{12}\ln (\e_{12}^{-1})&\mbox{if}& n < 6,\\
\displaystyle\sum_{j=1}^2\frac{\ln\l_j}{\l_j^{4}}+\e_{12}^{\frac{n+2}{2(n-2)}}\ln (\e_{12}^{-1})&\mbox{if}& n = 6,\\
\displaystyle\sum_{j=1}^2\frac{1}{\l_j^\frac{n+2}{2}}+\e_{12}^{\frac{n+2}{2(n-2)}}\ln
(\e_{12}^{-1})&\mbox{if}& n > 6.
\end{array}
\right.
\end{align}
The proof of \eqref{vestimate} is similar to that of Eq. $(3.7)$ in
\cite{BFG}. \\Our main result in this section is stated as follows.
\begin{pro}\label{t:4}
Let $n\geq 6$ and let $(u_\varepsilon)$ be a family of sign-changing
solutions of $(P_\varepsilon)$ having the expansion \eqref{h:1}
 and satisfying \eqref{gamma}-\eqref{vpart}. Assume that
\begin{align}\label{hh:1}\frac{\l_{1,\e}}{\l_{2,\e}}\rightarrow \infty\;\;\mbox{as}\;\;\varepsilon\rightarrow 0
\, \hbox{ such that } \ln\l_{1,\e}\leq c'\ln\l_{2,\e} \, \hbox{ for
some  } c'>1.\end{align} Then the concentration points
$\xi_{1,\varepsilon}$, $\xi_{2,\varepsilon}$ and the concentration
speeds $\l_{1,\varepsilon}$, $\l_{2,\varepsilon}$ satisfy
\begin{align*}\l_{2,\varepsilon}
|\xi_{1,\varepsilon}-\xi_{2,\varepsilon}|\to 0 \mbox{ as }
\varepsilon \to 0,\, \mbox{ and }\, \end{align*}
\begin{equation}\label{h:14} \frac{\Gamma_2}{3}\varepsilon\frac{\l_{2,\e}^{n-2}}{\ln \l_{2,\e}}
 \to \frac{1}{\overline{\Lambda}^2},\quad
\xi_{\varepsilon,i}\to  \xi_0 \in \Omega \mbox{ for } \,\,
i=1,2,\,\,\mbox{and}\,\,\,\frac{\l_{1,\varepsilon}}{\l_{2,\varepsilon}^3}\to
\overline{\Lambda}^\frac{4}{n-2},\end{equation} where $\Gamma_2$
is a positive constant defined in \eqref{t3}.\\
Furthermore, $\xi_0 $ is a critical point of the Robin function $R$
 and $\overline \Lambda$ satisfies \begin{equation}\label{Lambda}
R(\xi_0)\overline\Lambda^2 = 4.\end{equation}
\end{pro}
We point out that the condition $\ln\l_{1,\e}\leq c'\ln\l_{2,\e}$ is
very helpful in our argument. In fact, thanks to this condition and
 $\l_{1,\e}/\l_{2,\e}\rightarrow \infty$ which implies $\ln\l_{2,\e}\leq \ln\l_{1,\e}$, we derive that the quantities
 $\ln \l_{i,\e}$'s are of the same order,
 and this will
help us to analyze the balancing conditions introduced in Lemma
\ref{D2}.\\
For simplicity we shall write $P\delta_i$ for $P\d_{(\xi_{i
,\e},\l_{i,\e})}$, $\xi_i $ for $\xi_{i ,\e}$ and $\l_i$ for
$\l_{i,\e}$. Throughout this section, we restrict ourselves to the
case
$n \geq 6$.\\
Next, we provide two balancing conditions, that need to be satisfied
by
 the parameters of concentration, required in the proof of Proposition \ref{t:4}. The first one
 concerns
 the concentration speeds.
Multiplying the equation $-\Delta u_\e=f_\e(u_\e)$ by
$\l_i\frac{\partial P\d_{i}}{\partial \l_i}$ and integrating by
parts with taking into account \eqref{vestimate} and \eqref{di}, we
obtain
\begin{lem}\label{D2}
Let $n\geq 6$ and $u_\e=\gamma_1 P\d_1+\gamma_2 P\d_2+v_\e$ a
solution of $(P_\e)$. For each $i,l\in \{1,2\}$ such that $l\neq i$,
we have the following expansion
\begin{align*}
 &\gamma_i\Gamma_1\displaystyle
\frac{\a_i^p\e}{\ln \lambda_i}-
 (n-2)\overline{c}_1\frac{\gamma_i}{2}\frac{H(\xi_i,\xi_i)}{\l_i^{n-2}}
-\overline{c}_1\gamma_l\biggl(\l_i\frac{\partial \e_{il}}{\partial
\l_i}+\frac{n-2}{2}\frac{H(\xi_i,\xi_l)}{(\l_i\l_l)^{(n-2)/2}}
\biggr)  \\
&=O\left(\frac{\e}{\ln(\l_i)^2}+\frac{1}{\l_i^{2n-4}}+\frac{1}{\l_i^{(1-\tau)n}}+\e^2\ln(\ln\l_1)^2+\e^2\ln(\ln\l_2)^2
+ \frac{\ln\l_1 }{\l_1^{n}}+ \frac{\ln\l_2}{\l_2^{n}}+
\e_{12}^{\frac{n}{n-2}}\ln(\e_{12}^{-1})\right)
\end{align*}
where $\Gamma_1$ and $\overline{c}_1$ are the same constant defined
in Proposition \ref{p23} and $\tau$ is a positive constant small
enough.
\end{lem}
The proof of Lemma \eqref{D2} that we omit here looks like the
expansion of $\left< \nabla I_\e(u),\l_i\frac{\partial
P\d_i}{\partial \l_i}\right>$ given in Proposition \ref{p23}.
Indeed, being a solution of $(P_\e)$, $u_\e$ is a critical point of
$I_\e$ and therefore we have
$\left< \nabla I_\e(u_\e),\l_i\frac{\partial P\d_i}{\partial \l_i}\right>=0$.\\
Our second balancing condition focus on the concentration points and
it is obtained in similar way. Namely, we have
\begin{lem}\label{D3} Let $n\geq 6$ and $u_\e=\gamma_1 P\d_1+\gamma_2 P\d_2+v_\e$ a solution of
$(P_\e)$. For each
 $i,l\in \{1,2\}$ such that $l\neq i$ and $j\in \{1,...,n\}$, we have the following expansion
\begin{align*}
& \frac{\gamma_i}{2}
 \frac{\overline{c}_1}{\l_i^{n-1}}\frac{\partial H(\xi_i,\xi_i)}{
 \partial (a)_j}
 -\overline{c}_1\gamma_l
 \frac{1}{\l_i}\biggr(\frac{\partial \e_{il}}{\partial
(\xi_i)_j}-\frac{1}{(\l_i\l_l)^{(n-2)/2}}\frac{\partial H}{\partial
(a)_j}(\xi_i,\xi_l) \biggr)\\
&=O\bigg(\e^2\ln(\ln\l_1)^2+\e^2\ln(\ln\l_2)^2
 +\frac{\ln\l_1 }{\l_1^{n}}+\frac{\ln\l_2 }{\l_2^{n}}
 +\l_l | \xi_i-\xi_l|\e_{12}^{\frac{n+1}{n-2}
}+\e_{12}^{\frac{n}{n-2}}\ln(\e_{12}^{-1})
\bigg),
\end{align*}
where $\frac{\partial H}{\partial (a)_j}$ denotes the partial
derivative of $H$ with respect to the $j-$th component of the first
variable.
\end{lem}
 \begin{pfn}{Proposition
\ref{t:4}}: We recall that we have $\lambda_1/\lambda_2\rightarrow
+\infty$ as $\varepsilon \rightarrow 0$. From Proposition \ref{D2},
it is easy to obtain  the following estimate
 \begin{equation}\label{t0}
 \frac{\varepsilon}{\ln \l_i}=O\biggl(  \frac{1}{\lambda_2 ^{n-2}}
 +\varepsilon_{12}\biggr), \hbox{ for }i=1,2.
 \end{equation}
Let $K$ be a compact set in $\O$. It is easy to see that
\begin{equation*} 0<c\leq H(x,y)\leq C,\quad \forall x,y\in K
.\end{equation*} Note that the concentration points $\xi_1$ and
$\xi_2$ satisfy \eqref{di}. Thus, there exists a compact set $K$ in
$\O$ such that $\xi_1 ,\, \xi_2\in K$ for any $\e$ and therefore
\begin{equation}\label{minmajH}
 0<c\leq H(\xi_i,\xi_j)\leq C \quad\hbox{ for } i,j=1,2.
\end{equation}
We start by the following lemma.
\begin{lem}\label{l:1a}Under the assumptions of  Proposition
\ref{t:4}, we have
\begin{align*}\lambda_2|\xi_1-\xi_2|\rightarrow 0 \mbox{ as }\varepsilon\rightarrow
0.\end{align*}
\end{lem}
\begin{pf}
Note that from \eqref{hh:1}, we have $\ln\l_2< \ln\l_1\leq
c'\ln\l_2$ and up to subsequence we get
\begin{equation}\label{gam}\ln\l_1=\kappa\ln \l_2\, (1+o(1))
\end{equation}
for some real number $\kappa\geq 1$. In the sequel, we distinguish two cases: $\kappa>1$ and $\kappa=1$. \\
Let us start by the case $\kappa>1$. In a first step, we prove that
$\lambda_2|\xi_1-\xi_2|$ is bounded. In fact, arguing by
contradiction, we assume that $\lambda_2|\xi_1-\xi_2|\rightarrow
+\infty$ as $\varepsilon\rightarrow 0$. Then
$\lambda_1|\xi_1-\xi_2|\rightarrow +\infty$ as
$\varepsilon\rightarrow 0$ since $\l_1$ and $\l_2$ satisfy
\eqref{hh:1}. Thus, we get
\begin{align}
\label{t1}
&\varepsilon_{12}=\frac{1}{(\lambda_1\lambda_2|\xi_1-\xi_2|^2)^{(n-2)/2}}
+o(\varepsilon_{12}),\\&\label{t2} \lambda_i\frac{\partial
\varepsilon_{12}}{\partial \lambda_i} =-\frac{n-2}{2}
\frac{1}{(\lambda_1\lambda_2|\xi_1-\xi_2|^2)^{(n-2)/2}}
+o(\varepsilon_{12}).
 \end{align}
Using \eqref{t0} and \eqref{t2}, Lemmas \ref{D2} and \ref{D3} become
 \begin{equation}\label{t3}
\frac{H(\xi_i,\xi_i)}{\lambda_i ^{n-2}}+\frac{G(\xi_1,\xi_2)}
{(\lambda_1\lambda_2)^{(n-2)/2}}-\Gamma_2\frac{\varepsilon}{\ln
\l_i}=o\biggl(\varepsilon_{12}+ \frac{1}{\lambda_2 ^{n-2}} \biggr),
\hbox{ here }\Gamma_2=\frac{2\Gamma_1}{\overline{c}_1(n-2)}>0,
 \end{equation}
 \begin{equation}\label{t3'} \frac{1}{\lambda_i^{n-1}}\frac{\partial
H(\xi_i,\xi_i)}{\partial
a}+\frac{2}{\lambda_i}\frac{1}{(\lambda_1\lambda_2)^\frac{n-2}{2}}\frac{\partial
G(\xi_i,\xi_l)}{\partial a}=o\biggr(\varepsilon_{12}^\frac{n-1}{n-2}
+\ \frac{1}{\lambda_2 ^{n-1}} \biggr).\end{equation}
 Multiplying \eqref{t3} by $\ln\l_i$
for $i=1,2$ and we subtract, we get
$$-\frac{H(\xi_2,\xi_2)}{\lambda_2 ^{n-2}}+(\kappa-1)\frac{G(\xi_1,\xi_2)}
{(\lambda_1\lambda_2)^{(n-2)/2}}=o\biggl(\varepsilon_{12}+
\frac{1}{\lambda_2 ^{n-2}} \biggr),$$ where we have used
\eqref{hh:1}, \eqref{gam} and \eqref{minmajH}. This equation, with
taking into account \eqref{hh:1}, \eqref{defg} and \eqref{minmajH},
implies
\begin{equation}\label{s1}c\varepsilon_{12}\leq\frac{1}{\lambda_2^{n-2}}\leq C\varepsilon_{12},\end{equation}
i.e $\frac{1}{\lambda_2^{n-2}}$ and $\varepsilon_{12}$ are
comparable.
 Using \eqref{defg}, \eqref{t1}, \eqref{s1} together with \eqref{hh:1} and the fact that  $| \partial H(\xi_1,\xi_2)
/\partial \xi_i|\leq c$, we obtain
 \begin{equation}\label{s2} \frac{1}{(\lambda_1\lambda_2)^{(n-2)/2}} \bigg|
\frac{1}{\lambda_i}\frac{\partial G}{\partial \xi_i}(\xi_1,\xi_2)
\bigg|\geq c\varepsilon_{12}^\frac{n-1}{n-2}\;\;\;\mbox{for }
i=1,2.\end{equation}
 Clearly, \eqref{t3'} for $i=1$, \eqref{s1}
and \eqref{s2} give a contradiction and therefore $\lambda_2|\xi_1-\xi_2|$ is bounded.\\
On the second step, we claim that $\lambda_2|\xi_1-\xi_2|=o(1)$. In
fact, arguing by contradiction, we assume that
$\lambda_2|\xi_1-\xi_2|\nrightarrow 0$ as $\varepsilon\rightarrow0$.
Thus using the first step, there exist two positive constants $c'_1$
and $c'_2$ such that $c'_1<\lambda_2|\xi_1-\xi_2|<c'_2$. In this
case, we have
\begin{align}\label{zz1}
 & \varepsilon_{12}\mbox{ and }
(\lambda_2/\lambda_1)^{(n-2)/2} \mbox{ are of the same
order},\\
 & \varepsilon_{12}=\frac{1}{(\frac{\lambda_1}{\lambda_2}
+\lambda_1\lambda_2|\xi_1-\xi_2|^2)^{(n-2)/2}}+o(\varepsilon_{12})
\label{zz2},\\
 & \lambda_1\frac{\partial
\varepsilon_{12}}{\partial\lambda_1}=-\frac{n-2}{2}
\frac{1}{(\frac{\lambda_1}{\lambda_2}+\lambda_1\lambda_2|\xi_1-\xi_2|^2)^{(n-2)/2}}
+o(\varepsilon_{12})\notag,\\
 & \lambda_2\frac{\partial
\varepsilon_{12}}{\partial\lambda_2}=
\frac{(n-2)\lambda_1/\lambda_2}{(\frac{\lambda_1}{\lambda_2}
+\lambda_1\lambda_2|\xi_1-\xi_2|^2)^\frac{n}{2}}
-\frac{(n-2)/2}{(\frac{\lambda_1}{\lambda_2}+\lambda_1\lambda_2|\xi_1-\xi_2|^2)^\frac{n-2}{2}}
+o(\varepsilon_{12})\notag.
\end{align}
Thus, using Lemma \ref{D2} through \eqref{hh:1}, \eqref{t0} and
\eqref{minmajH}, we obtain
\begin{equation} \label{w1'}
 \frac{1}{(\frac{\lambda_1}{\lambda_2}+
\lambda_1\lambda_2|\xi_1-\xi_2|^2)^{(n-2)/2}}-\Gamma_2
 \frac{\varepsilon}{\ln \l_1}=o\biggr(\frac{1}{\lambda_2^{n-2}}+
\varepsilon_{12}\biggr),\end{equation}
\begin{align}
 \label{w2'} \frac{H(\xi_2,\xi_2)}
{\lambda_2^{n-2}} & +\frac{1}{(\frac{\lambda_1}{\lambda_2}+
\lambda_1\lambda_2|\xi_1-\xi_2|^2)^\frac{n-2}{2}}\notag\\
 & -\frac{\lambda_1}{\lambda_2}\frac{2}{(\frac{\lambda_1}{\lambda_2}
+\lambda_1\lambda_2|\xi_1-\xi_2|^2)^\frac{n}{2}}-\Gamma_2
 \frac{\varepsilon}{\ln \l_2}
=o\biggr(\frac{1}{\lambda_2^{n-2}}+\varepsilon_{12}\biggr).
\end{align}
Multiplying \eqref{w2'} by $\ln\l_2$ and \eqref{w1'} by $\ln\l_1$
and we subtract, we obtain
\begin{equation}\label{w3'}\frac{H(\xi_2,\xi_2)}
{\lambda_2^{n-2}}+
\frac{1-\kappa}{(\frac{\lambda_1}{\lambda_2}+\lambda_1\lambda_2|\xi_1-\xi_2|^2)^{(n-2)/2}}
-\frac{\lambda_1}{\lambda_2}
\frac{2}{(\frac{\lambda_1}{\lambda_2}+\lambda_1\lambda_2|\xi_1-\xi_2|^2)^{n/2}}
=o\biggr(\frac{1}{\lambda_2^{n-2}}+\varepsilon_{12}\biggr),
\end{equation}
since we have used \eqref{t0} and \eqref{gam}. \\
Recall that $\kappa>1$ and we have $c\leq H(\xi_2,\xi_2)\leq C$ from
\eqref{minmajH}. So \eqref{w3'} and \eqref{zz2} imply that
$\varepsilon_{12}$ and $1/\lambda_2^{n-2}$ are of the same order.
One one hand, using Lemma \ref{D3} for $i=2$, we obtain
\begin{align}\label{w4'}\frac{1}{\lambda_2^{n-1}}\frac{\partial H(\xi_2,\xi_2)}{\partial
a}-\frac{1}{\lambda_2}\frac{\partial \varepsilon_{12}}{\partial
\xi_2}=o\biggr(\e_{12}^\frac{n-1}{n-2}\biggr).
\end{align}
On the other hand, since we have assumed that
$c_1'\leq\lambda_2|\xi_1-\xi_2|\leq c_2'$, an easy computations show
that
\begin{align}\label{w5'}&\biggr|\frac{1}{\lambda_2}\frac{\partial \varepsilon_{12}}{\partial
\xi_2}\biggr|\geq
c\varepsilon_{12};\quad\frac{1}{\lambda_2^{n-1}}\frac{\partial
H(\xi_2,\xi_2)}{\partial
a}=O\biggr(\frac{1}{\lambda_2^{n-1}}\biggr)=o(\varepsilon_{12}).
\end{align}
We note that \eqref{w4'} and \eqref{w5'} give a contradiction.\\
 Thus, we get $\lambda_2|\xi_1-\xi_2| \to 0$
as $\varepsilon \to 0$ which completes the proof for the case
$\kappa>1$.\\
When $\kappa=1$, we argue similarly. We follow the same steps, but
without multiplying the equations obtained from Lemma \ref{D2} by
the coefficient $\ln\l_i$. In fact, to get rid of the case
$\lambda_2|\xi_1-\xi_2|\rightarrow +\infty$ and to prove that it is
bounded, we argue by contradiction. Subtracting the equations
\eqref{t3} for $i=1, 2$ and using \eqref{gam}(with $\kappa=1$), we
obtain that $1/\l_2^{n-2}=o(\e_{12})$ instead of $1/\l_2^{n-2}$ and
$\e_{12}$ are of the same order, found in the first case. This will
help us to conclude in the first step. Concerning the second step,
namely excluding the case $c'_1<\lambda_2|\xi_1-\xi_2|<c'_2$, we
subtract the equations \eqref{w1'} and \eqref{w2'} we get, as in the
case $\kappa>1$, $1/\l_2^{n-2}$ and $\e_{12}$ are of the same order
and we conclude similarly.
\end{pf}

Notice that Lemma \ref{l:1a} implies, as $\varepsilon \rightarrow
0$, that
\begin{align}\label{w1}&\varepsilon_{12}=
\biggr(\frac{\lambda_2}{\lambda_1}\biggr)^\frac{n-2}{2}+o(\varepsilon_{12}),
\\&\lambda_1\frac{\partial
\varepsilon_{12}}{\partial\lambda_1}=-\frac{n-2}{2}
\biggr(\frac{\lambda_2}{\lambda_1}\biggr)^\frac{n-2}{2}+o(\varepsilon_{12})
\quad \hbox{ and }\quad\lambda_2\frac{\partial
\varepsilon_{12}}{\partial\lambda_2}=
\frac{n-2}{2}\biggr(\frac{\lambda_2}{\lambda_1}\biggr)^\frac{n-2}{2}+
o(\varepsilon_{12})\label{w3}.
\end{align}
Using Lemma \ref{D2}, \eqref{w3}, \eqref{hh:1} and \eqref{t0}, we
obtain
\begin{align}&\label{w5}\biggr(\frac{\lambda_2}{\lambda_1}\biggr)^{(n-2)/2}
-\Gamma_2 \frac{\varepsilon}{\ln \l_1}
=o\biggr(\frac{1}{\lambda_2^{n-2}}+
\varepsilon_{12}\biggr),\\&\label{w6}\frac{H(\xi_2,\xi_2)}
{\lambda_2^{n-2}}-\biggr(\frac{\lambda_2}{\lambda_1}\biggr)
^{(n-2)/2}-\Gamma_2
 \frac{\varepsilon}{\ln \l_2}=o\biggr(\frac{1}{\lambda_2^{n-2}}+\varepsilon_{12}\biggr).
\end{align}
We start by proving that the constant $\kappa $ introduced in
\eqref{gam}
need to be greater than $1$. We argue by contradiction.\\
Assume that $\kappa=1$. Subtracting Equations \eqref{w5} and
\eqref{w6} and using \eqref{t0} and \eqref{gam}, we get rid of the
terms $\e/\ln\l_i$. More precisely, we obtain
\begin{equation}\label{leps}
2\biggr(\frac{\lambda_2}{\lambda_1}\biggr)^{(n-2)/2}
-\frac{H(\xi_2,\xi_2)} {\lambda_2^{n-2}}
=o\biggr(\frac{1}{\lambda_2^{n-2}}+ \varepsilon_{12}\biggr)
 \end{equation} Taking into account
\eqref{w1}, \eqref{leps} and \eqref{minmajH}, we obtain
$(\lambda_2/\lambda_1)^{(n-2)/2}$ and $1/\lambda_2^{n-2}$ are of the
same order, which implies that $\lambda_1$ and $\lambda_2^3$ are of
the same order. This contradicts \eqref{gam} with $\kappa=1$. Thus we have  $\kappa>1$.\\
Subtracting now the products of \eqref{w5} with $\ln\l_1$ and Eq.
\eqref{w6} multiplied by $\ln\l_2$ and taking into account
\eqref{w1} and \eqref{w3} together with \eqref{t0} and \eqref{gam},
we derive as previously that
\begin{equation}\label{sameorder}\varepsilon_{12}\hbox{
and }1/\lambda_2^{n-2}\hbox{ are of the same order},\end{equation}
which implies that $\lambda_1$ and $\lambda_2^3$ are of the same
order. We deduce the existence of a positive constant $k_0$ such
that
\begin{equation}\label{w10'}
 \frac{\lambda_1}{\lambda_2^3}=k_0(1+o(1)).
\end{equation}
This result asserts that $\kappa=3$ and we get
\begin{equation}\label{lneq}\ln\l_1=3\ln\l_2\; (1+o(1)).
\end{equation}
 Furthermore, Lemma \ref{D3} for $i=2$ and \eqref{hh:1} imply
\begin{equation}\label{w10} \frac{1}{\lambda_2^{n-1}}\frac{\partial
H(\xi_2,\xi_2)}{\partial
a}-\frac{1}{\l_2}\frac{\partial\e_{12}}{\partial
\xi_2}=o\biggr(\varepsilon_{12}^\frac{n-1}{n-2} +\
\frac{1}{\lambda_2 ^{n-1}} \biggr).\end{equation}
Under the previous assumptions, we obtain the following result.
\begin{lem}\label{e12o}
$$\frac{1}{\l_1}|\frac{\partial\e_{12}}{\partial
(\xi_1)_\ell}|=o\biggr(\frac{\l_2}{\l_1} \frac{1}{\lambda_2 ^{n-1}}
\biggr), \hbox{ for each }1\leq \ell\leq n.$$
\end{lem}
\begin{pf}
For sake of simplicity, we may assume that $\ell =1$. The same
argument holds true for $\ell=2,\ldots,n$.\\ Let
$P\psi_1^1:=\frac{1}{\l_1}\frac{\partial P\delta_1}{\partial
(\xi_1)_1}$ and $\underline{u}:= \sum_{j=1}^2\gamma_jP\delta_j$. On
one hand, we have
\begin{equation}\label{usol}
\displaystyle\left< \n I_\e(u_\e),P\psi_{1}^1\right>=0
\end{equation}
since $u_\e$ is a solution of $(P_\e)$. On the other hand, arguing
as in the proof of \cite[Proposition 2.11]{BFG}, we get
\begin{align}\label{d2d}
\displaystyle\left< \n
I_\e(u),P\psi_{1}^1\right>=&<u,P\psi_1^1>-\int_{\O}f_\e(u)P\psi_1^1\nonumber\\
=&<u,P\psi_1^1>-\int_{\O}\left[f_\e(u)-f_\e(\underline{u})-f'_\e(\underline{u})v\right]P\psi_1^1-\int_{\O}
f_\e(\underline{u})P\psi_1^1-\int_{\O}f'_\e(\underline{u})vP\psi_1^1\nonumber\\
=&<\underline{u},P\psi_1^1>-A-B-C
\end{align}
since $v\in E_{(\xi,\lambda)}$. 
Note that \cite[Proposition 2.11]{BFG} concerns the expansion of
$\left< \n I_\e(u),\l_i\partial P\delta_i/\partial \l_i\right>$.
Thanks to the fact that $P\psi_1^1$ satisfies $|P\psi_1^1|\leq
c\delta_1$ as well, some remaining terms in \eqref{d2d} may be
computed as in
\cite{BFG}.\\
 By the mean value theorem, there
exists $\theta=\theta(x)\in(0, 1)$ such that
$$A=\int_{\O}\left[f_\e'(\underline{u}+\theta v)-f'_\e(\underline{u})\right]vP\psi_1^1.$$
As computed in \cite{BFG}, we have
 $$
  |A|\leq c\left\{\begin{array}{ll}\|v\|^2\quad \quad\ \quad \quad \quad\hbox{ if }n\leq 6;\\
 \| v \| ^2 +   ( \| v \|^p + \e \|v\| ) \e_{12}^{1/2} \ln(\e_{12}^{-1})^{(n-2)/(2n)}\hbox{ if }n>6
 .\end{array}\right.
$$
  Thus we obtain
   \begin{equation}\label{Aestimate}
  A=o\biggr(\frac{\l_2}{\l_1} \frac{1}{\lambda_2 ^{n-1}}
\biggr) \hbox{ for }n\geq 6,
  \end{equation}
where we have used \eqref{vestimate}, \eqref{t0}, \eqref{w10'} and
the fact that
$\varepsilon_{12}$ and $1/\lambda_2^{n-2}$ are of the same order.\\
We also compute $C$ as in \cite{BFG}. Note that
$\ln\ln(e+|\underline{u}|)=O( \ln\ln\l_{1})=O( \ln\ln\l_{2})$ since
we have \eqref{lneq}. Using \eqref{f'} and the fact that
$$f_0'(\underline{u})=p[P\delta_1]^{p-1}+O([P\d_2]^{p-1}{\bf{1}}_{\{P\d_1\leq
P\d_2\}}+[P\delta_1]^{p-2}P\d_2{\bf{1}}_{\{P\d_2\leq P\d_1\}}),$$ we
get, in view of \eqref{vestimate}, \eqref{t0}, \eqref{sameorder} and
\eqref{w10'} the following
\begin{align}\label{Cdev}
C&:=\displaystyle\int_{\O}f'_\e(\underline{u})vP\psi_1^1\nonumber\\
&=\int_{\O}f'_0(\underline{u})vP\psi_1^1+\int_{\O}[f'_\e(\underline{u})-f'_0(\underline{u})]vP\psi_1^1\nonumber\\
&=p\int_{\O}[P\delta_1]^{p-1}P\psi_1^1 v+O \left(\int_{P\d_1\leq
P\d_2}[P\delta_2]^{p-1}P\delta_1 |v| + \int_{P\d_2\leq
P\d_1}[P\delta_1]^{p-1}P\delta_2| v|
+\|v\|\e\ln\ln\l_{2}\right)\nonumber\\
 &=C_1+ O\left(\e^2\ln(\ln\l_{2})^2+\|v\|^2+
\e_{12}^{\frac{n+2}{n-2}}\ln(\e_{12}^{-1})^\frac{n+2}{n}+\e_{12}^{2}\ln(\e_{12}^{-1})^{\frac{2(n-2)}{n}}
\right)\nonumber\\
&=C_1+ o\biggr(\frac{\l_2}{\l_1} \frac{1}{\lambda_2 ^{n-1}} \biggr).
\end{align}
As in \cite{BFG}, we have
\begin{align}\label{d2f}
C_1
&=\left\{\begin{array}{ll}
                                                      O\left(\|v\|\frac{\ln(\l_1d_1)}{(\l_1d_1)^4}\right), & \hbox{if }n=6; \\
                                                      O\left(\frac{\|v\|}{(\l_1d_1)^{(n+2)/2}}\right), & \hbox{if
}n>6,
                                                     \end{array}\right.\nonumber\\
                                                     &=o\biggr(\frac{\l_2}{\l_1} \frac{1}{\lambda_2 ^{n-1}} \biggr)
\end{align}
since $v$ satisfies \eqref{vestimate} and we have \eqref{t0},
\eqref{w10'} and \eqref{sameorder} i.e. $\varepsilon_{12}$ and
$1/\lambda_2^{n-2}$ are of the same order.
From \eqref{Cdev} and \eqref{d2f}, we get
\begin{equation}\label{cetimate}
C=o\biggr(\frac{\l_2}{\l_1} \frac{1}{\lambda_2 ^{n-1}} \biggr).
\end{equation}
Now, we compute $B$. Using \eqref{di}, \eqref{hh:1} and Proposition
\ref{p21}, we write
$$\underline{u}= \gamma_1\d_1+\gamma_2\d_2+O(\frac{1}{\l_2^\frac{n-2}{2}})$$
and, by \eqref{B'}, we get
$$f_\e(\underline{u})=f_\e\left(\gamma_1\d_1+\gamma_2\d_2\right)+O\left(\frac{1}{\l_2^\frac{n-2}{2}}
\left[|\gamma_1\d_1+\gamma_2\d_2|^{p-1}+(\frac{1}{\l_2^\frac{n-2}{2}})^{p-1}\right]
\right).$$ Thus
\begin{align}\label{Bestimate}
B:=&\int_{\O}f_\e(\underline{u})P\psi_1^1\nonumber\\
=&\int_{\O}f_\e(\underline{u})\frac{1}{\l_1}
\frac{\partial\d_1}{\partial
(\xi_1)_1}+\int_{\O}f_\e(\underline{u})\frac{1}{\l_1}
\frac{\partial\varphi_1}{\partial (\xi_1)_1}\nonumber\\
=&\int_{\O}\left(f_\e(\gamma_1\d_1+\gamma_2\d_2)+O\left(\frac{1}{\l_2^\frac{n-2}{2}}
\left[(\gamma_1\d_1+\gamma_2\d_2)^{p-1}+(\frac{1}{\l_2^\frac{n-2}{2}})^{p-1}\right]\right)\right)\frac{1}{\l_1}
\frac{\partial\d_1}{\partial (\xi_1)_1}
\nonumber\\
&+
O\left(\frac{1}{\l_1^\frac{n}{2}}\int_{\O}f_\e(\underline{u})\right)\nonumber\\
=&\int_{\O}f_\e(\gamma_1\d_1+\gamma_2\d_2)\frac{1}{\l_1}
\frac{\partial\d_1}{\partial (\xi_1)_1}+O\left(
\int_{\O}\left[\frac{1}{\l_2^\frac{n-2}{2}}(\d_1^{p-1}+\d_2^{p-1})+\frac{1}{\l_2^\frac{n+2}{2}}\right]|\frac{1}{\l_1}
\frac{\partial\d_1}{\partial
(\xi_1)_1}|\right)+O\left(\frac{1}{\l_2^\frac{3n}{2}}\right)
\end{align}
where we have used Proposition \ref{p21}, \eqref{w10'} and the fact
that $\int_{\O}f_\e(\underline{u})=o(1)$.\\ Note that, using
\eqref{w10'}, we have $$o\biggr(\frac{\l_2}{\l_1} \frac{1}{\lambda_2
^{n-1}} \biggr)= o\biggr( \frac{1}{\lambda_2 ^{n+1}} \biggr) .$$ In
the sequel, we will write the estimates of the remaining terms of
$B$ in
their $\l_2-$orders.\\
Taking into account \eqref{Bestimate}, \eqref{w10'}, the fact that
$$\int_{\O}\frac{1}{\l_1}
|\frac{\partial\d_1}{\partial (\xi_1)_1}|\leq
c\int_{\O}|x-\xi_1|\d_1^\frac{n}{n-2}\leq
\frac{c}{\l_1^\frac{n}{2}}\int_{\O}\frac{dx}{|x-\xi_1|^{n-1}}\leq\frac{c}{\l_1^\frac{n}{2}}$$
and
$$\int_{\O}\d_1^{p-1}\frac{1}{\l_1}
|\frac{\partial\d_1}{\partial (\xi_1)_1}|\leq
c\int_{\O}|x-\xi_1|\d_1^\frac{n+4}{n-2}\leq
\frac{c}{\l_1^\frac{n-2}{2}},$$ we obtain
\begin{align}\label{Bdef}
B=&\int_{\O}f_\e(\gamma_1\d_1+\gamma_2\d_2)\frac{1}{\l_1}
\frac{\partial\d_1}{\partial
(\xi_1)_1}+o\left(\frac{1}{\l_2^{n+1}}\right)\nonumber\\
=:&B'+o\left(\frac{1}{\l_2^{n+1}}\right).
\end{align}
Now, we introduce a new function, that is
$$\widetilde{\d}_2:=\d_{(\xi_1,\l_2)}.$$
We write
$\gamma_1\d_1+\gamma_2\d_2=\gamma_1\d_1+\gamma_2\widetilde{\d}_2+\gamma_2(\d_2-\widetilde{\d}_2)$
and thus, by the mean value theorem there exists some $\theta \in
(0,1)$ such that
\begin{align}\label{mvt}f_\e(\gamma_1\d_1+\gamma_2\d_2)=&f_\e(\gamma_1\d_1+\gamma_2\widetilde{\d}_2)+
f'_\e(\gamma_1\d_1+\gamma_2\widetilde{\d}_2)\gamma_2(\d_2-\widetilde{\d}_2)\nonumber\\
&+\frac{1}{2}f_{\e}''(\gamma_1\d_1+\gamma_2\widetilde{\d}_2+\theta
\gamma_2
 (\d_2- \widetilde{\d}_2))(\d_2- \widetilde{\d}_2)^2.
\end{align}
An easy computation, through the result of Lemma \ref{l:1a}, gives
us the following
\begin{equation}\label{lemA}
|\d_2-\widetilde{\d}_2|\leq c \l_2 |\xi_1-\xi_2|\widetilde{\d}_2.
\end{equation}
By using \eqref{mvt} and \eqref{lemA}, the integral $B'$ becomes
\begin{align}\label{B'def}
B'=&\int_{\O}f_\e(\gamma_1\d_1+\gamma_2\widetilde{\d}_2)\frac{1}{\l_1}\frac{\partial\d_1}{\partial
(\xi_1)_1}+\int_{\O}
f'_\e(\gamma_1\d_1+\gamma_2\widetilde{\d}_2)\gamma_2(\d_2-\widetilde{\d}_2)\frac{1}{\l_1}\frac{\partial\d_1}{\partial
(\xi_1)_1}\nonumber\\
&+\frac{1}{2}\int_{\O}f_{\e}''(\gamma_1\d_1+\gamma_2\widetilde{\d}_2+\theta
\gamma_2
 (\d_2- \widetilde{\d}_2))(\d_2- \widetilde{\d}_2)^2\frac{1}{\l_1}\frac{\partial\d_1}{\partial
(\xi_1)_1}\nonumber\\
=:&B'_1+B'_2+B'_3.
\end{align}
Recall that we have $d_1=d(\xi_1,\partial \O)>c$ from assumption
\eqref{di}. We get
\begin{align}\label{B1'def}
B'_1&:=\int_{\O}f_\e(\gamma_1\d_1+\gamma_2\widetilde{\d}_2)\frac{1}{\l_1}\frac{\partial\d_1}{\partial
(\xi_1)_1}\nonumber\\
&=\int_{B(\xi_1,d_1/2)}\ldots +\int_{\O\setminus B(\xi_1,d_1/2)}\ldots\nonumber\\
&=O\left(\int_{
B(\xi_1,d_1/2)^c}\d_1^{p+1}+\widetilde{\d}_2^p\frac{\d_1}{\l_1|x-\xi_1|}\right)
\end{align}
since the function
$f_\e(\gamma_1\d_1+\gamma_2\widetilde{\d}_2)\frac{1}{\l_1}\frac{\partial\d_1}{\partial
(\xi_1)_1}$ is antisymmetric with respect to $(x-\xi_1)_1$ in
$B(\xi_1,d_1/2)$ and we have used that
$$|\frac{1}{\l_1}\frac{\partial\d_1}{\partial
(\xi_1)_1}(x)|\leq c \frac{\d_1(x)}{\l_1|x-\xi_1|}.$$ A simple
computations lead to
$$\int_{
B(\xi_1,d_1/2)^c}\d_1^{p+1}= O(\frac{1}{\l_1^n})\quad \hbox{ and }
\int_{B(\xi_1,d_1/2)^c}\widetilde{\d}_2^p\frac{\d_1}{\l_1|x-\xi_1|}=
O(\frac{1}{\l_1^\frac{n}{2}\l_2^\frac{n+2}{2}}).$$ Hence,
\eqref{B1'def} and \eqref{w10'} assert that
\begin{equation}\label{B1'}
B'_1=o\left(\frac{1}{\l_2^{n+1}}\right).
\end{equation}
Now, expanding as follows
$$f'_\e(\gamma_1\d_1+\gamma_2\widetilde{\d}_2)=f'_\e(\gamma_1\d_1)+f'_\e(\gamma_2\widetilde{\d}_2)
+O\left(\inf(\d_1,\widetilde{\d}_2)^{p-1}\right),$$ the second
integral $B'_2$ is written
\begin{align}\label{B'2def}
B'_2=&\int_{\O}\gamma_2f'_\e(\gamma_1\d_1)[\d_2-\widetilde{\d}_2]\frac{1}{\l_1}\frac{\partial\d_1}{\partial
(\xi_1)_1}
+\int_{\O}\gamma_2f'_e(\gamma_2\widetilde{\d}_2)[\d_2-\widetilde{\d}_2]\frac{1}{\l_1}\frac{\partial\d_1}{\partial
(\xi_1)_1}\nonumber\\
&+O\left(\int_{\O}\inf(\d_1,\widetilde{\d}_2)^{p-1}|\d_2-\widetilde{\d}_2|\frac{1}{\l_1}|\frac{\partial\d_1}{\partial
(\xi_1)_1}|\right)\nonumber\\
=:&B'_{21}+B'_{22}+B'_{23}.
\end{align}
Observe that
\begin{equation}\label{equivalence}\widetilde{\d}_2(x)\leq
\d_1(x)\Leftrightarrow\l_1\l_2|x-\xi_1|^2\leq 1\Leftrightarrow x\in
B(\xi_1,\frac{1}{\l_2^2}).
\end{equation}
Taking into account Lemma \ref{l:1a} ,\eqref{w10'}, \eqref{lemA},
\eqref{equivalence} and the fact that
$|\frac{1}{\l_1}\frac{\partial\d_1}{\partial (\xi_1)_1}(x)|\leq c
|x-\xi_1|\d_1^\frac{n}{n-2}(x)$, we get
\begin{align}\label{B'23}
B'_{23}=&\int_{\widetilde{\d}_2\leq
\d_1}\inf(\d_1,\widetilde{\d}_2)^{p-1}|\d_2-\widetilde{\d}_2|\frac{1}{\l_1}|\frac{\partial\d_1}{\partial
(\xi_1)_1}|+\int_{\d_1\leq\widetilde{\d}_2}\inf(\d_1,\widetilde{\d}_2)^{p-1}|\d_2-\widetilde{\d}_2|\frac{1}{\l_1}|\frac{\partial\d_1}{\partial
(\xi_1)_1}|\nonumber\\
\leq &c\l_2|\xi_1-\xi_2|\left(\int_{\widetilde{\d}_2\leq
\d_1}\widetilde{\d}_2^{p}(x)|x-\xi_1|\d_1^\frac{n}{n-2}\,dx+\int_{\d_1\leq\widetilde{\d}_2}
\widetilde{\d}_2(x)|x-\xi_1|\d_1^\frac{n+4}{n-2}(x)\,dx\right)\nonumber\\
\leq
&c\l_2|\xi_1-\xi_2|\left(\int_{B(\xi_1,1/\l_2^2)}\l_2^\frac{n+2}{2}\frac{dx}{\l_1^\frac{n}{2}|x-\xi_1|^{n-1}}+
\int_{B(\xi_1,1/\l_2^2)^c}
\l_2^\frac{n-2}{2}\frac{dx}{\l_1^\frac{n+4}{2}|x-\xi_1|^{n+3}}\right)\nonumber\\
\leq
&c\l_2|\xi_1-\xi_2|\left(\frac{\l_2^\frac{n+2}{2}}{\l_1^\frac{n}{2}}\int_{0}^{1/\l_2^2}\;dr+\frac{\l_2^\frac{n-2}{2}}{\l_1^\frac{n+4}{2}}
\int_{\frac{1}{\l_2^2}}^\infty \frac{dr}{r^4}\right)\nonumber\\
\leq
&c\l_2|\xi_1-\xi_2|\left(\left(\frac{\l_2}{\l_1}\right)^\frac{n}{2}\frac{1}{\l_2}+\frac{1}{\l_2^{n+7}}(\l_2^2)^3\right)\nonumber\\
=&o\left(\frac{1}{\l_2^{n+1}}\right).
\end{align}
In view of \eqref{f'}, \eqref{lemA}, \eqref{t0}, \eqref{sameorder},
\eqref{w10'}, Lemma \ref{l:1a} and the fact that $\ln\ln (\d_1)\leq
c\ln\ln \l_2$, we have
\begin{align}
B'_{21}:=&\int_{\O}\gamma_2f'_\e(\gamma_1\d_1)[\d_2-\widetilde{\d}_2]\frac{1}{\l_1}\frac{\partial\d_1}{\partial
(\xi_1)_1}\nonumber\\
=&\int_{\O}\gamma_2f'_0(\gamma_1\d_1)[\d_2-\widetilde{\d}_2]\frac{1}{\l_1}\frac{\partial\d_1}{\partial
(\xi_1)_1}+O\left(\e\ln\ln \l_2
\int_{\O}\d_1^{p-1}\l_2|\xi_1-\xi_2|\widetilde{\d}_2\frac{1}{\l_1}\frac{\partial\d_1}{\partial
(\xi_1)_1}\right)\nonumber\\
=&\gamma_2p\int_{\mathbb{R}^n}\d_1^{p-1}\d_2\frac{1}{\l_1}\frac{\partial\d_1}{\partial
(\xi_1)_1}+O\left(\int_{\O^c}\d_1^{p+1}+\e\ln\ln \l_2
\l_2|\xi_1-\xi_2|
\l_2^\frac{n-2}{2}\int_{\O}\d_1^{p}\right)\nonumber\\
=&\frac{\gamma_2\overline{c}_1}{\l_1}\frac{\partial\e_{12}}{\partial
(\xi_1)_1}+O(\l_2|\xi_1-\xi_2|\e_{12}^\frac{n+1}{n-2})+O(\frac{1}{\l_1^n})+o\left((\frac{\l_2}{\l_1})^\frac{n-2}{2}
\e\ln\ln\l_2\right)\nonumber\\
=&\frac{\gamma_2\overline{c}_1}{\l_1}\frac{\partial\e_{12}}{\partial
(\xi_1)_1}+o\left(\frac{1}{\l_2^{n+1}}\right)
\end{align}
where we have used the fact that
$\d_1^{p-1}\widetilde{\d}_2\frac{1}{\l_1}\frac{\partial\d_1}{\partial
(\xi_1)_1}$ is antisymmetric with respect to $(x-\xi_1)_1$ in
$\mathbb{R}^n$ and the following result
\begin{equation}\label{F11}p\int_{\mathbb{R}^n}\d_1^{p-1}\d_2\frac{1}{\l_1}\frac{\partial\d_1}{\partial
(\xi_1)_1}=\int_{\mathbb{R}^n}\d_2^p\frac{1}{\l_1}\frac{\partial\d_1}{\partial
(\xi_1)_1}=\frac{\overline{c}_1}{\l_1}\frac{\partial\e_{12}}{\partial
(\xi_1)_1}+O(\l_2|\xi_1-\xi_2|\e_{12}^\frac{n+1}{n-2})\end{equation}
given in \cite[Estimate F11]{B1}. Here $\overline{c}_1$ is the positive constant introduced in Proposition \ref{p23}. \\
Using again \eqref{f'} and \eqref{lemA} and Lemma \ref{l:1a}, we
obtain
\begin{align}\label{B'22}
B'_{22}:=&
\int_{\O}\gamma_2f'_\e(\gamma_2\widetilde{\d}_2)[\d_2-\widetilde{\d}_2]\frac{1}{\l_1}\frac{\partial\d_1}{\partial
(\xi_1)_1}\nonumber\\
=&\int_{\O}\gamma_2f'_0(\gamma_2\widetilde{\d}_2)[\d_2-\widetilde{\d}_2]\frac{1}{\l_1}\frac{\partial\d_1}{\partial
(\xi_1)_1}+O\left(\l_2|\xi_1-\xi_2|\e\ln\ln\l_2\int_{\O}\widetilde{\d}_2^{p}\frac{1}{\l_1}\frac{\partial\d_1}{\partial
(\xi_1)_1}\right)\nonumber\\
=&\int_{\O}\gamma_2\widetilde{\d}_2^{p-1}[\d_2-\widetilde{\d}_2]\frac{1}{\l_1}\frac{\partial\d_1}{\partial
(\xi_1)_1}+o\left(\e\ln\ln \l_2
\left(\frac{\l_2}{\l_1}\right)^\frac{n}{2}\right)\nonumber\\
=:&B''_{22}+o\left(\frac{1}{\l_2^{n+1}}\right)
\end{align}
where we have used \eqref{t0}, \eqref{w10'} and the fact that
\begin{align}\label{e12n}
\int_{\O}\widetilde{\d}_2^{p}\frac{1}{\l_1}|\frac{\partial\d_1}{\partial
(\xi_1)_1}|&\leq c
\int_{\O}\widetilde{\d}_2^{p}|x-\xi_1|\d_1^\frac{n}{n-2}\nonumber\\
&\leq
\frac{c}{\l_1^\frac{n}{2}}\int_{\O}\frac{\widetilde{\d}_2^{p}}{|x-\xi_1|^{n-1}}
\nonumber\\
&\leq
\frac{c}{\l_1^\frac{n}{2}}\frac{1}{\l_2^\frac{n-2}{2}}\l_2^{n-1}\int_{\mathbb{R}^n}\frac{dy}{(1+|y|^2)^\frac{n+2}{2}|y|^{n-1}}
\nonumber\\
&\leq c\left(\frac{\l_2}{\l_1}\right)^\frac{n}{2}\int_{0}^\infty
\frac{dr}{(1+|r|^2)^\frac{n+2}{2}}
\nonumber\\
&\leq c\left(\frac{\l_2}{\l_1}\right)^\frac{n}{2}.
\end{align}
 Observe that, using
\eqref{lemA}, we get
$$\d_2^p=\widetilde{\d}_2^p+p\widetilde{\d}_2^{p-1}(\d_2-\widetilde{\d}_2)+
O\left(\l_2^2|\xi_1-\xi_2|^2\widetilde{\d}_2^{p}\right).$$ Thus,
using \eqref{F11}, Lemma \ref{l:1a} and the antisymmetry of the
function
$\widetilde{\d}_2^{p}\frac{1}{\l_1}\frac{\partial\d_1}{\partial
(\xi_1)_1}$ with respect to $(x-\xi_1)_1$ in $B(\xi_1, d_1/2)$, we
obtain
\begin{align}\label{B"22}
B''_{22}:=&\gamma_2\int_{\O}\widetilde{\d}_2^{p-1}[\d_2-\widetilde{\d}_2]\frac{1}{\l_1}\frac{\partial\d_1}{\partial
(\xi_1)_1}\nonumber\\
=&\frac{\gamma_2}{p}\int_{\O}\d_2^{p}\frac{1}{\l_1}\frac{\partial\d_1}{\partial
(\xi_1)_1}-\frac{\gamma_2}{p}\int_{\O}\widetilde{\d}_2^{p}\frac{1}{\l_1}\frac{\partial\d_1}{\partial
(\xi_1)_1}+O\left(\l_2^2|\xi_1-\xi_2|^2\int_{\O}\widetilde{\d}_2^{p}\frac{1}{\l_1}|\frac{\partial\d_1}{\partial
(\xi_1)_1}|\right)\nonumber\\
=&\frac{\gamma_2}{p}\int_{\mathbb{R}^n}\d_2^{p}\frac{1}{\l_1}\frac{\partial\d_1}{\partial
(\xi_1)_1}-\frac{\gamma_2}{p}\int_{B(\xi_1,
d_1/2)}\widetilde{\d}_2^{p}\frac{1}{\l_1}\frac{\partial\d_1}{\partial
(\xi_1)_1}+O\left(\int_{B(\xi_1,
d_1/2)^c}\widetilde{\d}_2^p\d_1+\l_2^2|\xi_1-\xi_2|^2\int_{\O}\widetilde{\d}_2^{p}|x-\xi_1|\d_1^\frac{n}{n-2}\right)\nonumber\\
=&\frac{\gamma_2\overline{c}_1}{p}\frac{1}{\l_1}\frac{\partial\e_{12}}{\partial
(\xi_1)_1}+O\left(\l_2|\xi_1-\xi_2|\e_{12}^\frac{n+1}{n-2}+\frac{1}{\l_1^{n+1}}+
\l_2^2|\xi_1-\xi_2|^2\frac{1}{\l_1^\frac{n}{2}}\int_{\O}\widetilde{\d}_2^{p}\frac{1}{|x-\xi_1|^{n-1}}\right)\nonumber\\
=&\frac{\gamma_2\overline{c}_1}{p}\frac{1}{\l_1}\frac{\partial\e_{12}}{\partial
(\xi_1)_1}+o\left(\frac{1}{\l_2^{n+1}}+\l_2|\xi_1-\xi_2|\left(\frac{\l_2}{\l_1}\right)^\frac{n}{2}\right)\nonumber\\
=&\frac{\gamma_2\overline{c}_1}{p}\frac{1}{\l_1}\frac{\partial\e_{12}}{\partial
(\xi_1)_1}(1+o(1))+o\left(\frac{1}{\l_2^{n+1}}\right)
\end{align}
where we have used, in the last equality, \eqref{w1} i.e. $\e_{12}$
and $(\l_2/\l_1)^\frac{n-2}{2}$ are of the same order and the fact
that
$$\frac{1}{\l_1}|\frac{\partial\e_{12}}{\partial (\xi_1)_1}|=(n-2)\l_2|(\xi_1-\xi_2)_1|\e_{12}^\frac{n}{n-2}.$$
In view of \eqref{B'2def}, \eqref{B'23}-\eqref{B'22} and
\eqref{B"22}, we get
\begin{equation}\label{B'2}
B'_2=\gamma_2\overline{c}_1(1+\frac{1}{p})\frac{1}{\l_1}\frac{\partial\e_{12}}{\partial
(\xi_1)_1}(1+o(1))+o\left(\frac{1}{\l_2^{n+1}}\right).
\end{equation}
To estimate $B'_3$, we split the integral as follows:
\begin{align}\label{B3}
B'_3&:=\frac{1}{2}\int_{\O}f_{\e}''(\gamma_1\d_1+\gamma_2\widetilde{\d}_2+\theta
\gamma_2
 (\d_2- \widetilde{\d}_2))(\d_2- \widetilde{\d}_2)^2\frac{1}{\l_1}\frac{\partial\d_1}{\partial
(\xi_1)_1}\nonumber\\
&=\int_{\O_1}\ldots +\int_{\O_1^c}\ldots\nonumber\\
 &=:B'_{31}+B'_{32}.
\end{align}
where $\O_1:=\{x:\ |\gamma_1\delta_1(x)+\gamma_2\delta_2(x)|\leq
M\l_2|\xi_1-\xi_2|\widetilde{\d}_2(x)\}$. Here M is a positive
constant greater than $4c$
 where $c$ is the constant introduced in \eqref{lemA}. Using this fact, it is easy to see
 that in $\O_1^c$ we have
 $$|\gamma_1\d_1+\gamma_2\widetilde{\d}_2+\theta
\gamma_2
 (\d_2- \widetilde{\d}_2)|\geq \frac{M}{2}|\xi_1-\xi_2|\widetilde{\d}_2 \hbox{ for each }\theta\in (0,1).$$
 Thus in view of \eqref{lemA}, \eqref{B''}, Lemma \ref{l:1a} and the fact that $p-2\leq 0$ for $n\geq 6$, we get
 \begin{align}\label{B32}
B'_{32}&=O\left(\int_{\O_1^c}\left(\l_2|\xi_1-\xi_2|\widetilde{\d}_2\right)^{p-2}(\d_2-
\widetilde{\d}_2)^2\frac{1}{\l_1}|\frac{\partial\d_1}{\partial
(\xi_1)_1}|\right)\nonumber\\
&=O\left((\l_2|\xi_1-\xi_2|)^{p-1}\l_2|\xi_1-\xi_2|\int_{\O}\widetilde{\d}_2^{p}\frac{1}{\l_1}|\frac{\partial\d_1}{\partial
(\xi_1)_1}|\right)\nonumber\\
&=o(\frac{1}{\l_1}|\frac{\partial \e_{12}}{\partial(\xi_1)_1}|)
 \end{align}
 where we have used \eqref{e12n} and the fact that $\e_{12}$ and
 $(\frac{\l_2}{\l_1})^\frac{n}{2}$ are of the same order.\\
Observe that, in $\O_1$, it holds
$|\gamma_1\d_1+\gamma_2\widetilde{\d}_2|\leq
 c \l_2|\xi_1-\xi_2|\widetilde{\d}_2$, by using \eqref{lemA}.
 Therefore, the expansion \eqref{mvt}, \eqref{B.1}, \eqref{lemA} and \eqref{e12n}
 imply
 \begin{align}\label{B31}
B'_{31}=&\int_{\O_1}\left(f_\e(\gamma_1\d_1+\gamma_2\d_2)-f_\e(\gamma_1\d_1+\gamma_2\widetilde{\d}_2)-
f'_\e(\gamma_1\d_1+\gamma_2\widetilde{\d}_2)\gamma_2(\d_2-\widetilde{\d}_2)\right)\frac{1}{\l_1}\frac{\partial\d_1}{\partial
(\xi_1)_1}\nonumber\\
&=O\left(\int_{\O}\left(\l_2|\xi_1-\xi_2|\widetilde{\d}_2\right)^{p}\frac{1}{\l_1}|\frac{\partial\d_1}{\partial
(\xi_1)_1}|\right)\nonumber\\
&=O\left((\l_2|\xi_1-\xi_2|)^{p-1}\l_2|\xi_1-\xi_2|\int_{\O}\widetilde{\d}_2^{p}\frac{1}{\l_1}|\frac{\partial\d_1}{\partial
(\xi_1)_1}|\right)\nonumber\\
&=o\left(\frac{1}{\l_1}|\frac{\partial
\e_{12}}{\partial(\xi_1)_1}|\right).
\end{align}
Hence, \eqref{B3}, \eqref{B32} and \eqref{B31} assert
\begin{align}\label{B3'}
B'_{3}=o(\frac{1}{\l_1}|\frac{\partial
\e_{12}}{\partial(\xi_1)_1}|).
\end{align}
Lastly, we compute $<\underline{u},P\psi_1^1>$. In view of
Proposition \ref{p21}, \eqref{majdd}, \eqref{w10'}, \eqref{F11} and
Holder's inequality, we obtain
\begin{align}\label{scalar}
<\underline{u},P\psi_1^1>=&<\gamma_1\d_1+\gamma_2\d_2,\frac{1}{\l_1}\frac{\partial
\d_1}{\partial
(\xi_1)_1}>+<\underline{u},\frac{1}{\l_1}\frac{\partial
\varphi_1}{\partial
(\xi_1)_1}>\nonumber\\
=&\int_{\O}\gamma_2\d_2^p\frac{1}{\l_1}\frac{\partial \d_1}{\partial
(\xi_1)_1} +\int_{\O}\gamma_1\d_1^p\frac{1}{\l_1}\frac{\partial
\d_1}{\partial
(\xi_1)_1} +O\left(\frac{1}{\l_1^\frac{n}{2}}\int_{\O}\d_1^{p}+\d_2^{p}\right)\nonumber\\
=&\gamma_2\int_{\R^n}\delta_2^p\frac{1}{\lambda_1}\frac{\partial
\delta_1}{\partial (\xi_1)_1}
+\int_{B(\xi_1,d_1/2)}\gamma_1\d_1^p\frac{1}{\l_1}\frac{\partial
\d_1}{\partial (\xi_1)_1}
+O\left(\int_{\Omega^c}\delta_2^p\d_1 +\int_{B(\xi_1,d_1/2)^c}\d_1^{p+1}\right)+o(\frac{1}{\l_2^\frac{3n}{2}})\nonumber\\
=&\gamma_2\frac{\ov c_1}{\l_1}\frac{\partial\e_{12}}{\partial
(\xi_1)_1}+
O\left(\l_2|\xi_1-\xi_2|\e_{12}^\frac{n+1}{n-2}+\frac{1}{\l_1^\frac{n-2}{2}\l_2^\frac{n+2}{2}}+\frac{1}{\l_1^n}\right)+o\left(\frac{1}{\l_2^{n+1}}\right)
\nonumber\\
=&\gamma_2\frac{\ov c_1}{\l_1}\frac{\partial\e_{12}}{\partial
(\xi_1)_1}+o\left(\frac{1}{\l_2^{n+1}}\right).
\end{align}
Combining \eqref{usol}-\eqref{Aestimate}, \eqref{cetimate},
\eqref{Bdef}, \eqref{B'def}, \eqref{B1'}, \eqref{B'2}, \eqref{B3'}
and\eqref{scalar},  the proof of Lemma \ref{e12o} is completed.
\end{pf}\\
  Now, we will prove \eqref{h:14} and \eqref{Lambda}. From Lemma \ref{e12o},  we
  derive\begin{equation}\label{e12l2}
  \frac{1}{\l_2}\frac{\partial\e_{12}}{\partial
\xi_2}=o\biggr( \frac{1}{\lambda_2 ^{n-1}} \biggr).
  \end{equation}
Let us introduce the following change of variables
\begin{equation}\label{chv3}\left(\frac{\ln\l_2}{\lambda_2^{n-2}}\right)^{1/2}=
\Lambda\left(\frac{\Gamma_2}{3}\, \varepsilon\right)^{1/2}\, \hbox{
where }\Lambda>0.
\end{equation}
 Note that, multiplying \eqref{w5}, \eqref{w6} by respectively $\ln\l_1$ and $\ln\l_2$, and using
  \eqref{w10'} \eqref{lneq} and \eqref{chv3}, we obtain
\begin{equation}\label{w12}
 \frac{\Lambda^2}{k_0^\frac{n-2}{2}}-1=o(\Lambda^2);\quad
 H(\xi_2,\xi_2)\Lambda-\frac{4}{\Lambda}
=o(\Lambda).
\end{equation}
Taking into account \eqref{e12l2}, \eqref{chv3} and the fact that
$\e_{1,2}$ and $1/\l_2^{n-2}$ are of the same order, \eqref{w10}
implies
\begin{equation}\label{w13}\frac{\partial H(\xi_2,\xi_2)}{\partial a}\Lambda^2=
o\left(\Lambda^2\right). \end{equation}Recall that we have
$d(\xi_2,\partial \O)\geq c $. Thus, the functions $H$ and its
derivatives are bounded. Furthermore, the function $H$ satisfies
\eqref{minmajH}. Therefore, from \eqref{w12}, it is easy to see that
 $\Lambda$ is bounded from above and it is also bounded from below by some positive constant. Hence, $\Lambda$
converges to $\overline{\Lambda} >0$ (up to a subsequence). Let
$\xi_0$ be the limit of $\xi_2$ when $\e$ goes to $0$. Passing to
the limit in \eqref{w12} and \eqref{w13}, we get
\begin{equation}
k_0=\overline \Lambda^{\frac{4}{n-2}}; \quad
H(\xi_0,\xi_0)\overline{\Lambda} -4\overline{\Lambda}^{\,\,-1} =0;
\quad \frac{\partial H(\xi_0,\xi_0)}{\partial a} = 0,
\end{equation}
 which imply that
$\xi_0$ is a critical point of the Robin function $R$ and
$\overline{\Lambda}$ satisfies \eqref{Lambda}. Therefore, in view of
Lemma \ref{l:1a} and \eqref{w10'}, the conclusion in \eqref{h:14}
holds. The proof of Proposition \ref{t:4} is thereby completed.
\end{pfn}
\section{Description of the solution and its lower order term}
We denote by $i^*: L^\frac{2n}{n+2}(\O)\rightarrow H^1_0(\O)$ the
adjoint operator of the embedding $i: H^1_0(\O)\rightarrow
L^\frac{2n}{n-2}(\O)$, i.e. if $w\in L^\frac{2n}{n+2}(\O)$ then
$u=i^*(w)$ in $H^1_0(\O)$ is the unique solution of the equation
$-\Delta u=w$ in $\O$ and $u=0$ on $\partial\O$.\\
By the definition of the operator $i^*$, it is clear that $(P_\e)$
can be written as follows:
\begin{equation}\label{(2.2)}
u=i^*\left[f_\e(u)\right], \quad u\in  H^1_0(\O)
\end{equation}
where $f_\e$ is introduced in \eqref{nonpower}.\\
Next, we describe the shape of the solutions we are looking for.\\
Let $\xi$ be a point in $\O$ and,  given an integer number $k$, let
$\l_j$ for $j=1,...,k$, be positive parameters defined as multiple
of proper power of $\frac{\e}{|\ln \e|}$, namely
\begin{equation}\label{equ4}
\l_j=\biggr(\frac{\e}{|\ln
\e|}\biggr)^{-\frac{2(k-j)+1}{n-2}}\rho_j,\quad \mbox{with}\;
\rho_j>0.
\end{equation}
Let $\xi_j$, for $j=1,...,k$, also be $k$ points in $\O$ given by
\begin{equation}\label{equ5}
\xi_j=\xi+\l_j^{-1}\sigma_j,\quad
\mbox{with}\,\,\sigma_2,...,\sigma_k\in
\R^n\,\mbox{and}\,\,\sigma_1=0.
\end{equation}Let $\a_j$, for $j=1,...,k$, be positive parameters.
Fix a small $\eta> 0$  and assume that
\begin{equation}\label{equ6}
d(\xi,\partial \O)>\eta,\quad |\a_j-1|<\eta,\quad
\eta<\rho_j<\frac{1}{\eta}\quad\mbox{and}\quad
|\sigma_j|\leq\frac{1}{\eta}\quad\mbox{for}\,j=1,...,k.
\end{equation}
It is an immediate observation that
\begin{equation}\label{equ7}
\l_k=\biggr(\frac{\e}{|\ln
\e|}\biggr)^{-\frac{1}{n-2}}\rho_k,\quad\mbox{and}\,\,\frac{\l_{j}}{\l_{j+1}}=\biggr(\frac{\e}{|\ln
\e|}\biggr)^{-\frac{2}{n-2}}\frac{\rho_{j}}{\rho_{j+1}}.
\end{equation}
We look for a tower of sign-changing bubbles solution to $(P_\e)$,
of the form:
\begin{equation}\label{equ71}
u(x)=V(x)+v(x),\quad V(x)=V_{\ov\alpha,\ov \rho,\ov
\sigma,\xi}(x)=\sum_{j=1}^k (-1)^j\alpha_jP\d_{(\xi_j,\l_j)}(x),
\end{equation}
where we set $\ov\alpha=(\a_1,...,\a_k)\in \R_+^k$, $\ov
\rho=(\rho_1,...,\rho_k)\in \R_+^k$ and
$\ov\sigma=(\sigma_2,...,\sigma_k)\in \R^{(k-1)n}$.\\
The term $v$ has to be thought as a remainder term of lower order.
Let
\begin{align}\label{q3}&E:=
\Big\{v\in H^1_0(\O):\big<v,P\d_{i}\big>=\big<v,\frac{\partial
P\d_{i}}{\partial \l_i}\big>= \big<v,\frac{\partial
P\d_{i}}{\partial (\xi_i)_j}\big>=0\;\forall \; 1\leq j\leq
n,\,\forall \,1\leq i\leq k\Big\},
\end{align}
where $P\d_i=P\d_{(\xi_i,\l_i)}$ and $\xi_i^j$ is the $j^{th}$ component of $\xi_i$.\\
Let $\Pi: H^1_0(\O)\rightarrow E^\perp$ and $\Pi^\perp:
H^1_0(\O)\rightarrow E$ be the orthogonal projections. In order to
solve \eqref{(2.2)} we will solve the couple of equations:
\begin{equation}\label{(2.9)}
\Pi^\perp\left(V+v-i^*\left[f_\e(V+v)\right]\right)=0;
\end{equation}
\begin{equation}\label{(2.10)}
\Pi\left(V+v-i^*\left[f_\e(V+v)\right]\right)=0.
\end{equation}
Given $\xi$, $\ov \alpha$, $\ov \rho$ and $\ov \sigma$ satisfying
conditions \eqref{equ6}, one can solve uniquely \eqref{(2.9)} in
$v\in E$.
This solution, which will be denoted by $\overline{v}$ is the lower
order term in the description of the ansatz \eqref{equ71}. This is
the content of:
\begin{pro}\label{G21}
There exists $\e_0>0$ such that  for any $\xi\in \O$, $\ov \alpha\in
\R^k_+$, $\ov \rho\in \R^k_+$, $\ov \sigma\in \R^{(k-1)n}$
satisfying \eqref{equ4}-\eqref{equ6} and for any  $\e\in (0,\e_0)$
there exists a unique function $\ov v=\ov v(\e,\ov\alpha,\ov
\rho,\ov \sigma, \xi)\in E$ such
$$I_\e(V+\ov v(\e,\ov\alpha,\ov \rho,\ov \sigma, \xi))=\inf\{I_\e(V+v):v\in E\}. $$
Moreover $\ov v$ satisfies
\begin{align}\label{vest}
\|\overline{v}\| =\left\{\begin{array}{ll}O(\e\ln|\ln\e|
),&\hbox{ if }n\leq 6;\\
O\left(\left(\displaystyle\frac{\e}{|\ln\e|}\right)^\frac{p}{2}\right),
&\hbox{ if }n> 6.\end{array}\right.
\end{align}
In addition, there exists $(A,B,C)\in
\R^k\times\R^k\times(\R^{n})^k$ such that the following holds
\begin{equation}\label{v1}
-\Delta \ov v(\e,\ov\alpha,\ov \rho,\ov \sigma, \xi) =f_\e(V+\ov
v(\e,\ov\alpha,\ov \rho,\ov \sigma, \xi))-(-\Delta V)
+\sum_{i=1}^k(-\Delta)\biggr(A_iP\delta_i+B_i\frac{\partial
P\delta_i}{\partial \lambda_i}+\sum_{j=1}^n C_{ij}\frac{\partial
P\delta_i}{\partial (\xi_i)_j}\biggr).
\end{equation}
\end{pro}
The proof of such a result is contained in \cite{BFG} (see the proof
of \cite[Proposition 3.2]{BFG}). In the following, $\ov v$ is the
solution of
\eqref{(2.9)}.\\
 As in \cite{BFG}, we estimate the numbers  $A$, $B$
and $C$. Note that the functions $\frac{\partial P\delta_i}{\partial
\lambda_i}$ and $\frac{\partial P\delta_i}{\partial \xi_i^j}$, which
appear in \eqref{v1}, are multiplied respectively by $\l_i$ and
$\frac{1}{\l_i}$ in \cite{BFG}. Taking into account this change and
the $\e$-orders of the $\l_i's$ given in \eqref{equ4}, we get the
following.
\begin{pro}\label{pro42}For $i=1,...,k$ and $j=1,...,n$, we have
\begin{eqnarray*}\begin{cases}
 \displaystyle A_i=\left\{\begin{array}{ll}O(|1-\a_i|+\e\ln|\ln\e|
),&\hbox{ if }n\leq 6;\\
O\left(|1-\a_i|+\left(\displaystyle\frac{\e}{|\ln\e|}\right)^\frac{p}{2}\right),
&\hbox{ if }n> 6,\end{array}\right.\\
  \displaystyle      B_i=O\left(\frac{\lambda_i\varepsilon}{|\ln \e|}\right),    \\
      \displaystyle C_{ij}=O\left(\frac{|1-\a_i|\e}{\l_i|\ln \e|}+
      \frac{1}{\l_i}\left(\frac{\e}{|\ln \e|}\right)^\frac{n-1}{n-2}\right).
\end{cases}
\end{eqnarray*}
\end{pro}

As in \cite{R} and \cite{MP1}, since $|f_{\e}(u)|\leq |u|^p$, we
have the following estimate.
\begin{lem}\label{lem45}It holds
$$\int_{\partial \O}\biggr(\frac{\partial  \ov v}{\partial \nu}\biggr)^2dw=
\left\{\begin{array}{ll}O(\e^2\ln^2|\ln\e|
),&\hbox{ if }n\leq 6;\\
O\left(\left(\displaystyle\frac{\e}{|\ln\e|}\right)^p\right),
&\hbox{ if }n> 6.\end{array}\right.$$
\end{lem}
\begin{pf}We consider a smooth function $\zeta:\R^n\rightarrow \R$ such that
$$0\leq \zeta\leq 1 ,\quad \zeta(y)=0\quad\mbox{for}\quad|y-\xi|\leq \varrho,\quad \zeta(y)=1
 \quad\mbox{for}\quad|y-\xi|\geq 2\varrho$$
for some $\varrho>0$ (we choose $\varrho>0$ such that
$\varrho<d(\xi,\partial \O)/4$). It is elementary to see that $\zeta
\ov v$ is a solution of the following problem
\begin{eqnarray} \label{equaa1}\begin{cases}
 -\Delta(\zeta  \ov v)= g & \mbox{in} \, \, \O, \\
 \zeta  \ov v= 0     &  \mbox{on}\, \, \partial \O,
\end{cases}
  \end{eqnarray}
where
\begin{align}
g :=&\zeta\Big(f_\e(V+\ov
v)-\sum_{i=1}^k\a_i\d_{(\xi_i,\l_i)}^p\Big)+\zeta\sum_{i=1}^k\Big(A_i\d_{(\xi_i,\l_i)}^p+pB_i\d_{(\xi_i,\l_i)}^{p-1}\frac{\partial
\d_{(\xi_i,\l_i)}}{\partial \lambda_i}\notag\\&+p\sum_{j=1}^n
C_{ij}\d_{(\xi_i,\l_i)}^{p-1}\frac{\partial
\d_{(\xi_i,\l_i)}}{\partial (\xi_i)^j}\Big)-\Delta\zeta \ov
v-2\n\zeta\n \ov v.
\end{align}
Since $\zeta \ov v$ is a solution of\eqref{equaa1}, the following inequality holds:
\begin{equation}\label{equaa2}
\biggr|\frac{\partial}{\partial \nu}(\zeta  \ov v)\biggr|^2_{2,\partial \O}=\biggr|\frac{\partial  \ov v}{\partial \nu}\biggr|^2_{2,\partial \O}\leq C|g|^2_{\frac{2n}{n+1},\O}.
\end{equation}
Hence, we need to estimate $|g|^2_{2n/(n+1)}$. We observe that
\begin{align}\label{equaa3}
g=&O\biggr( \zeta|\ov v|^{p }+\sum_{i=1}^k\zeta\d_{(\xi_i,\l_i)}^p
+\sum_{i=1}^k\Big(|A_i|\zeta\d_{(\xi_i,\l_i)}^p+|B_i|\zeta\d_{(\xi_i,\l_i)}^{p-1}|\frac{\partial
\d_{(\xi_i,\l_i)}}{\partial \lambda_i}|\notag\\&+\sum_{j=1}^n
|C_{ij}|\zeta\d_{(\xi_i,\l_i)}^{p-1}|\frac{\partial
\d_{(\xi_i,\l_i)}}{\partial (\xi_i)^j}|\Big)+|\Delta\zeta \ov
v|+|\n\zeta\n  \ov v|\biggr).
\end{align}
Let us denote $q:=2n/(n+1)$. We have to estimate each term of
\eqref{equaa3} in the $L^q(\O)$-norm. It holds
\begin{align*}
 &\biggr|\sum_{i=1}^k\zeta\d_{(\xi_i,\l_i)}^{p}\biggr|_q^2=O\biggr(\biggr(\sum_{i=1}^k\int_{|x-\xi|\geq\varrho}\d_{(\xi_i,\l_i)}^{pq}\biggr)^{2/q}\biggr)
=O\biggr(\sum_{i=1}^k\frac{1}{\l_i^{n+2}}\biggr) \end{align*} and by
using Proposition \ref{pro42}, we get
\begin{align*}
&\biggr|\sum_{i=1}^k\Big(|A_i|\zeta\d_{(\xi_i,\l_i)}^p+|B_i||\zeta\d_{(\xi_i,\l_i)}^{p-1}\frac{\partial
\d_{(\xi_i,\l_i)}}{\partial \lambda_i}|+\sum_{j=1}^n
|C_{ij}||\zeta\d_{(\xi_i,\l_i)}^{p-1}\frac{\partial
\d_{(\xi_i,\l_i)}}{\partial
\xi_i^j}|\Big)\biggr|_q^2\\=&O\biggr(\sum_{i=1}^k|A_i|^2\biggr(\int_{|x-\xi|\geq\varrho}\d_{(\xi_i,\l_i)}^{pq}\biggr)^{2/q}
+\sum_{i=1}^k|B_i|^2\biggr(
\int_{|x-\xi|\geq\varrho}\Big|\d_{(\xi_i,\l_i)}^{p-1}\frac{\partial
\d_{(\xi_i,\l_i)}}{\partial
\lambda_i}\Big|^q\biggr)^{2/q}\\&+\sum_{i=1}^k\sum_{j=1}^n
|C_{ij}|^2\biggr(\int_{|x-\xi|\geq\varrho}|\d_{(\xi_i,\l_i)}^{p-1}\frac{\partial
\d_{(\xi_i,\l_i)}}{\partial
\xi_i^j}|^q\biggr)^{2/q}\biggr)\\=&O\biggr(\biggr(\sum_{i=1}^k\Big(|A_i|^2+\frac{|B_i|^2}{\l_i^2}+\l_i^2\sum_{j=1}^n|C_{ij}|^2\biggr)
\biggr(\int_{|x-\xi|\geq\varrho}\d_{(\xi_i,\l_i)}^{pq}\biggr)^{2/q}\biggr)=O\biggr(\sum_{i=1}^k\frac{1}{\l_i^{n+2}}\biggr).
\end{align*}
Furthermore, since $q<2<p+1$, by Holder's inequalities, we have
\begin{align*}
&\biggr(\int_{\O}|\Delta\zeta|^q|\ov v|^q\biggr)^{2/q}
=O\biggr(\biggr(\int_{\O}|  \ov v|^{p+1}\biggr)^{2/(p+1)}
\biggr(\int_{\O}|\Delta
\zeta|^{\frac{q(p+1)}{p+1-q}}\biggr)^{\frac{2(p+1-q)}{q(p+1)}}\biggr)
=O(\|\ov v\|^2),
\\&\biggr(\int_{\O}|\n\zeta|^q|\n\ov  v|^q\biggr)^{2/q}=O\biggr(\|\ov v\|^2
\biggr(\int_{\O}|\n \zeta|^{\frac{2q}{2-q}}\biggr)^{\frac{2 -q}{q
}}\biggr)=O(\|\ov v\|^2).
\end{align*}
The only remaining term is $\Big(\int_{\O}(\zeta|\ov
v|^{p})^q\Big)^{2/q}$ which is the most difficult to estimate
because $pq >p + 1$. We multiply equation \eqref{v1} by $\zeta
^{\frac{2(n-2)}{n+1}}|\ov v|^{\frac{2}{n+1}}\ov v$, and we integrate
on $\O$. Arguing as in \cite{R} (see also \cite{MP1}),integration by
parts and the Sobolev embedding theorem lead to the inequality
\begin{align*}
&-\int_\O\Delta\ov v \zeta ^{\frac{2(n-2)}{n+1}}|\ov v|^{\frac{2}{n+1}}\ov v\geq C\biggr(\int_\O(\zeta|\ov v|^p)^q\biggr)^{2/q}+O(\|\ov v\|^{2(n+2)/(n+1)}).
 \end{align*}
 On the other hand, we write
 \begin{align*}
 &  f_\e(u)-\sum_{i=1}^k(-1)^i\a_i\d_{(\xi_i,\l_i)}^p + \sum_{i=1}^k\Big(A_i\d_{(\xi_i,\l_i)}^p+pB_i\d_{(\xi_i,\l_i)}^{p-1}\frac{\partial \d_{(\xi_i,\l_i)}}{\partial
\lambda_i} p\sum_{j=1}^n C_{ij}\d_{(\xi_i,\l_i)}^{p-1}\frac{\partial \d_{(\xi_i,\l_i)}}{\partial
(\xi_i)^j}\Big) \\&=O\biggr(|\ov v|^p+\sum_{i=1}^k\delta_{(\xi_i,\l_i)}^p+\sum_{i=1}^k\Big(|A_i| \d_{(\xi_i,\l_i)}^p+|B_i|| \d_{(\xi_i,\l_i)}^{p-1}\frac{\partial \d_{(\xi_i,\l_i)}}{\partial
\lambda_i}| +\sum_{j=1}^n |C_{ij}|| \d_{(\xi_i,\l_i)}^{p-1}\frac{\partial \d_{(\xi_i,\l_i)}}{\partial
(\xi_i)^j}|\Big)\biggr).
\end{align*}
Arguing exactly as in \cite[Appendix C]{R}, we get:
$$\int_\O\zeta ^{\frac{2(n-2)}{n+1}}|\ov v|^{p+\frac{n+3}{n+1}}=O\biggr(\biggr(\int_{\O} |\ov v|^{\frac{(p-1)n}{2}} \biggr)^{\frac{2}{n}}\biggr(\int_{\O}(\zeta|\ov v|^{p})^q\biggr)^{\frac{2}{p+1}}\biggr)=o\biggr(\biggr(\int_{\O}(\zeta|\ov v|^{p})^q\biggr)^{\frac{2}{p+1}}\biggr)\biggr).$$
By using Holder's inequality, the definition of the function $\zeta$
and Proposition \eqref{pro42}, we also have the validity of the
following estimates:
\begin{align*}
\sum_{i=1}^k\int_{\O}\zeta ^{\frac{2(n-2)}{n+1}}|\ov v|^{
\frac{n+3}{n+1}}\d_{(\xi_i,\l_i)}^p &=O\biggr(\biggr(\int_{\O} |\ov
v|^{p+1}
\biggr)^{\frac{n+3}{(n+1)(p+1)}}\sum_{i=1}^k\biggr(\int_{|x-\xi|\geq
\varrho} \d_{( \xi_i,\l_i)}^{\frac{2n(n+1)p}{n^2+n+6}}
\biggr)^{\frac{n^2+n+6}{2n(n+1)}}\biggr)\\&= O\biggr(\|\ov
v\|^{\frac{n+3}{n+1}}\sum_{i=1}^k
\frac{1}{\l_i^{\frac{n+2}{2}}}\biggr),
\end{align*}
and
\begin{align*}
&\sum_{i=1}^k\int_{\O}\zeta ^{\frac{2(n-2)}{n+1}}|\ov v|^{\frac{n+3}{n+1}}\Big(|A_i| \d_{(\xi_i,\l_i)}^p+|B_i|| \d_{(\xi_i,\l_i)}^{p-1}\frac{\partial \d_{(\xi_i,\l_i)}}{\partial
\lambda_i}| +\sum_{j=1}^n |C_{ij}|| \d_{(\xi_i,\l_i)}^{p-1}\frac{\partial \d_{(\xi_i,\l_i)}}{\partial
(\xi_i)^j}|\Big)\biggr)\\& =O\biggr(\|\ov v\|^{\frac{n+3}{n+1}}\sum_{i=1}^k\biggr(\int_{|x-\xi|\geq \rho} \Big(|A_i|\d_{( \xi_i,\l_i)}^p+|B_i|| \d_{(\xi_i,\l_i)}^{p-1}\frac{\partial \d_{(\xi_i,\l_i)}}{\partial
\lambda_i}|+\sum_{j=1}^n |C_{ij}|| \d_{(\xi_i,\l_i)}^{p-1}\frac{\partial \d_{(\xi_i,\l_i)}}{\partial
(\xi_i)^j}|\Big)^{\frac{2n(n+1)}{n^2+n+6}} \biggr)^{\frac{n^2+n+6}{2n(n+1)}}\biggr)
\\& =O\biggr(\|\ov v\|^{\frac{n+3}{n+1}}\sum_{i=1}^k\l_i^{-\frac{n+2}{2}}\Big(|A_i|+\frac{|B_i|}{\l_i}+\l_i\sum_{j=1}^n |C_{ij}| \Big)\biggr)
=O\biggr(\frac{\|\ov v\|^{\frac{n+3}{n+1}}}{\l_k^{\frac{n+2}{2}}}\biggr),
\end{align*}
so that finally
\begin{align*}
\biggr(\int_{\O}(\zeta|\ov v|^{p})^q\biggr)^{2/q}=O\biggr(\|\ov v\|^{\frac{2(n+2)}{n+1}}+\frac{\|\ov v\|^{\frac{n+3}{n+1}}}{\l_k^{\frac{n+2}{2}}}\biggr).
\end{align*}
Taking account of \eqref{equ4} and \eqref{vest}, the desired result
follows.
\end{pf}
\section{The reduced energy}
 We are left now to solve \eqref{(2.10)}, more precisely to find the point $\xi$, the points
 $\sigma_2,...,\sigma_k\in \R^n$ and the
parameters $\rho_1,...,\rho_k$ and $\a_1,...,\a_k$ so that
\eqref{(2.10)} is satisfied. It happens that this problem has a
variational structure, in the sense
that solving \eqref{(2.10)} is reduced to find critical points to some given explicit finite dimensional functional.\\
Let $\widetilde{I_\e}:\R_+^k\times \R^k_+\times \R^{(k-1)n}\times
\O\rightarrow\R$ be defined by
\begin{equation}
\widetilde{I_\e}(\ov\a,\ov
\rho,\ov\sigma,\xi):=I_{\e}\Big(\sum_{i=1}^k(-1)^i\a_iP\d_{(\xi_i,\lambda_i)}
+\ov v\Big).
\end{equation}
As in Part $1$ of \cite[Proposition 2.2]{MP1}, we have the
following: if $(\ov\a,\ov \rho,\ov\sigma,\xi)$ is a critical point
of $\widetilde{I_\e}$, then $u:=
\sum_{i=1}^k(-1)^i\alpha_iP\delta_{(\xi_i,\lambda_i)}+\overline{v}$
is a critical point of $I_\varepsilon$. This result reduces the
existence of solutions for $(P_{\e})$ to the problem of finding
critical points of the reduced energy
functional $\widetilde{I_\e}$.\\
 Now, we introduce the asymptotic expansion for the reduced energy functional  $\widetilde{I_\e}$
 and its partial derivative in terms of the parameters $(\ov\a,\ov \rho,\ov
 \sigma,\xi)$.\\
 Using Propositions \ref{p22} and \ref{G21}, we have the
 following.
\begin{pro}\label{dalpha}For $h=1,...,k$, we have
\begin{align}\label{cr10}
&\frac{\partial \widetilde{I_\e}}{\partial \a_h}(\ov\a,\ov \rho,\ov
\sigma,\xi)=\a_h(1-\a_h^{p-1})S^{n/2}_n+\left\{\begin{array}{ll}O(\e\ln|\ln\e|
),&\hbox{ if }n\leq 6;\\
O\left(\left(\displaystyle\frac{\e}{|\ln\e|}\right)^\frac{p}{2}\right),
&\hbox{ if }n> 6,\end{array}\right.\\& \frac{\partial
\widetilde{I_\e}}{\partial \a_h}(\ov\a,\ov \rho,\ov
\sigma,\xi)=0\Leftrightarrow\a_h=1+\left\{\begin{array}{ll}O(\e\ln|\ln\e|
),&\hbox{ if }n\leq 6;\\
O\left(\left(\displaystyle\frac{\e}{|\ln\e|}\right)^\frac{p}{2}\right),
&\hbox{ if }n> 6,\end{array}\right.\label{cr11}
\end{align}
$C^0$-uniformly with respect to $\xi$ in compact sets of $\O$, $\ov
\sigma$ in compact sets of $\R^{(k-1)n}$, $\ov \a$ in compact sets
of $\R^k_+$ and $\ov \rho$ in compact sets of $\R^k_+$.
\end{pro}
The expansion of $\widetilde{I_\e}$ can be stated as follows.
\begin{pro}\label{pro52}We have
\begin{align*}
\widetilde{I_\e}(\ov\a,\ov
\rho,\ov\sigma,\xi)=&\sum_{i=1}^k\Big(\frac{\a_i^2}{2}-\frac{\a_i^{p+1}}{p+1}\Big)S_n^{\frac{n}{2}}+
k\frac{(n-2)S_n^{\frac{n}{2}}}{2n}\e\ln|\ln\e|\biggr(1+\frac{1}{|\ln\e|}\biggr)+\e\frac{(n-2)S_n^{\frac{n}{2}}}{2n}\sum_{i=1}^k
\ln\biggr(
\frac{2(k-i)+1}{2}\biggr)\notag\\&+\frac{\e}{|\ln\e|}\Psi(\ov
\rho,\ov\sigma,\xi)+o\biggr(\frac{\e}{|\ln \e|}\biggr),
\end{align*}
$C^0$-uniformly with respect to $\xi$ in compact sets of $\O$, $\ov
\sigma$ in compact sets of $\R^{(k-1)n}$, $\ov\a$ in compact sets of
$\R^k_+$ and $\ov \rho$ in compact sets of $\R^k_+$. Here
\begin{equation}\Psi(\ov \rho,\ov\sigma,\xi)=\frac{\ov c_1}{2}\frac{H(\xi,\xi)}{\rho_k^{n-2}}+
\frac{(n-2)S_n^{\frac{n}{2}}}{2n} \sum_{i=1}^k\frac{\ln
\rho_i}{2(k-i)+1} +\ov
c_1\sum_{i=2}^{k}\Big(\frac{\rho_{i}}{\rho_{i-1}}\Big)^{\frac{n-2}{2}}\frac{1}{\big(1
+|\sigma_{i}|^2\big)^{\frac{n-2}{2}}}.
\end{equation}
\end{pro}
\begin{pf}
We estimate $\e_{i(i+1)}$, using \eqref{eijdef}, \eqref{equ5} and
\eqref{equa7}, we have
\begin{align}\label{ee2}
\e_{i(i+1)}&=\frac{\Big(\frac{\l_{i+1}}{\l_i}\Big)^{(n-2)/2}}{\Big(1
+\Big(\frac{\l_{i+1}}{\l_i}\Big)^{2}+|\frac{\l_{i+1}}{\l_i}\sigma_i-\sigma_{i+1}|^2\Big)^{(n-2)/2}}
\notag\\&=\frac{\Big(\frac{\l_{i+1}}{\l_i}\Big)^{(n-2)/2}}{\Big(1
+|\sigma_{i+1}|^2+\Big(\frac{\l_{i+1}}{\l_i}\Big)^{2}+|\frac{\l_{i+1}}{\l_i}\sigma_i|^2-2\frac{\l_{i+1}}{\l_i}<\sigma_i,\sigma_{i+1}>\Big)^{(n-2)/2}}
\notag\\&=\Big(\frac{\l_{i+1}}{\l_i}\Big)^{(n-2)/2}\frac{1}{\big(1
+|\sigma_{i+1}|^2\big)^{(n-2)/2}}\Big(1+O\big(\frac{\l_{i+1}}{\l_i}\big)\Big)\notag\\
&=\frac{\e}{|\ln
\e|}\Big(\frac{\rho_{i+1}}{\rho_i}\Big)^{(n-2)/2}\frac{1}{\big(1
+|\sigma_{i+1}|^2\big)^{(n-2)/2}} +o\biggr(\frac{\e}{|\ln
\e|}\biggr).
\end{align}
Note that, from the previous computation we derive that
$\e_{i(i+1)}$ and $(\l_{i+1}/\l_i)^{(n-2)/2}$ are of the same
order.\\Using \eqref{ee2}, we get
\begin{align}\label{ee11}
\l_h\frac{\partial \e_{h(h+1)}}{\partial
\l_h}&=-\frac{n-2}{2}\biggr\{\frac{\l_h}{\l_{h+1}}-\frac{\l_{h+1}}{\l_h}+\l_h\l_{h+1}|\xi_h-\xi_{h+1}|^2\biggr\}\e_{h(h+1)}^{n/(n-2)}\notag\\&=-
\frac{n-2}{2}\frac{\l_{h}}{\l_{h+1}}\biggr\{1-\biggr(\frac{\l_{h+1}}{\l_h}\biggr)^2+\Big|\frac{\l_{h+1}}{\l_h}\sigma_h
-\sigma_{h+1}\Big|^2\biggr\}\e_{h(h+1)}^{n/(n-2)}\notag\\&=-\frac{n-2}{2}\frac{\e}{|\ln
\e|}\biggr(\frac{\rho_{h+1}}{\rho_h}\biggr)^{(n-2)/2}\frac{1}{\big(1
+|\sigma_{h+1}|^2\big)^{(n-2)/2}}+o\biggr(\frac{\e}{|\ln
\e|}\biggr),
\notag\\
\l_h\frac{\partial \e_{h(h-1)}}{\partial
\l_h}&=-\frac{n-2}{2}\biggr\{\frac{\l_h}{\l_{h-1}}-\frac{\l_{h-1}}{\l_h}+\l_h\l_{h-1}|\xi_h-\xi_{h-1}|^2\biggr\}\e_{h(h-1)}^{n/(n-2)}\notag\\&=-
\frac{n-2}{2}\frac{\l_{h-1}}{\l_{h}}\biggr\{\biggr(\frac{\l_h}{\l_{h-1}}\biggr)^2-1+\Big|\sigma_h
-\frac{\l_h}{\l_{h-1}}\sigma_{h-1}\Big|^2\biggr\}\e_{h(h-1)}^{n/(n-2)}\notag\\&
=\frac{n-2}{2}\frac{\e}{|\ln
\e|}\biggr(\frac{\rho_{h}}{\rho_{h-1}}\biggr)^{(n-2)/2}\frac{1-|\sigma_h|^2}{\big(1
+|\sigma_{h}|^2\big)^{n/2}}+o\biggr(\frac{\e}{|\ln \e|}\biggr),
\notag\\
\frac{1}{\l_h}\frac{\partial \e_{h(h+1)}}{\partial
\xi_h}&=-(n-2)\biggr(\frac{\l_{h+1}}{\l_h}\sigma_h-\sigma_{h+1}\biggr)\e_{h(h+1)}^{n/(n-2)}=o\biggr(\frac{\e}{|\ln
\e|}\biggr),\notag\\\frac{1}{\l_h}\frac{\partial
\e_{h(h-1)}}{\partial
(\xi_h)_j}&=-(n-2)\frac{\l_{h-1}}{\l_h}\biggr((\sigma_{h})_j-\frac{\l_h}{\l_{h-1}}(\sigma_{h-1})_j\biggr)\e_{h(h-1)}^{n/(n-2)}
\notag\\&=-(n-2)\frac{\e}{|\ln
\e|}\Big(\frac{\rho_{h}}{\rho_{h-1}}\Big)^{(n-2)/2}\frac{(\sigma_h)_j}{\big(1
+|\sigma_{h}|^2\big)^{n/2}}+o\biggr(\frac{\e}{|\ln \e|}\biggr),
\quad\mbox{for}\,j=1,...,n.
\end{align}
We have
\begin{align}\label{equa7}
\ln(\ln\l_i^\frac{n-2}{2})&=\ln\ln\left(\frac{n-2}{2}\right)+\ln|\ln\e|+\ln\biggr(\frac{2(k-i)+1}{n-2}\biggr)+\ln\biggr(1+
\frac{\ln|\ln\e|}{|\ln\e|}+\frac{n-2}{2(k-i)+1}\frac{\ln
\rho_i}{|\ln\e|}\biggr)\notag\\
&=\ln|\ln\e|\left(1+\frac{1}{|\ln\e|}\right)+\ln\biggr(\frac{2(k-i)+1}{2}\biggr)+\frac{n-2}{2(k-i)+1}\frac{\ln
\rho_i}{|\ln\e|}+O\biggr(\frac{(\ln|\ln\e|)^2}{(\ln\e)^2}\biggr).
\end{align}
Using Proposition \ref{p23'}, \eqref{equ4}, \eqref{equ5},
\eqref{equ6},\eqref{cr11}, \eqref{ee11} and \eqref{equa7}, we obtain
\begin{align}\label{ee1}
\widetilde{I_\e}(\ov\a,\ov
\rho,\ov\sigma,\xi)=&\sum_{i=1}^k\Big(\frac{\a_i^2}{2}-\frac{\a_i^{p+1}}{p+1}\Big)S_n^{\frac{n}{2}}
+k\frac{(n-2)S_n^{\frac{n}{2}}}{2n}\e\ln|\ln\e|\left(1+\frac{1}{|\ln\e|}\right)+\e\frac{(n-2)S_n^{\frac{n}{2}}}{2n}\sum_{i=1}^k
\ln\biggr(
\frac{2(k-i)+1}{2}\biggr)\notag\\&
-\ov
c_1\sum_{i=1}^{k-1}\gamma_i\gamma_{i+1}\Big(\frac{n+2}{2n}\a_i^{p}\a_{i+1}+\frac{n-2}{2n}\a_{i+1}^{p}\a_{i}-\frac{\a_i\a_{i+1}}{2}\Big)\e_{i(i+1)}\notag\\&-\ov
c_1\sum_{i=2}^{k}
\gamma_i\gamma_{i-1}\Big(\frac{n+2}{2n}\a_i^{p}\a_{i-1}+\frac{n-2}{2n}\a_{i-1}^{p}\a_{i}
-\frac{\a_i\a_{i-1}}{2}\Big)\e_{i(i-1)}\notag\\
&+\frac{\ov
c_1}{n}\sum_{i=1}^{k-1}\gamma_i\gamma_{i+1}\biggr(\a_i^p\a_{i+1}\l_i\frac{\partial\e_{i(i+1)}}{\partial\l_i}-\a_{i+1}^p\a_i\l_{i+1}\frac{\partial\e_{i(i+1)}}{\partial\l_{i+1}}\biggr)
\notag\\&+\frac{\ov
c_1}{n}\sum_{i=2}^{k}\gamma_i\gamma_{i-1}\biggr(\a_i^p\a_{i-1}\l_i\frac{\partial\e_{i(i-1)}}{\partial\l_i}-\a_{i-1}^p\a_i\l_{i-1}\frac{\partial\e_{i(i-1)}}{\partial\l_{i-1}}\biggr)
\notag\\& +\frac{\e}{|\ln \e|}\frac{\ov
c_1}{2}\frac{H(\xi,\xi)}{\rho_k^{n-2}}+\frac{\e}{|\ln
\e|}\frac{(n-2)^2S_n^{\frac{n}{2}}}{2n}
\sum_{i=1}^k\frac{1}{2(k-i)+1}\ln \rho_i+o\biggr(\frac{\e}{|\ln
\e|}\biggr).
\end{align}
Combining \eqref{ee1}, \eqref{ee2}, \eqref{ee11} and Proposition
\ref{dalpha}, the proof of Proposition \ref{pro52} follows.
\end{pf}
\begin{pro}\label{pro53}We have for  $h=2,...,k-1$ and $r=2,...,k,$
\begin{align*}
\frac{\partial \widetilde{I_\e}}{\partial \rho_h}(\ov\a,\ov \rho,\ov
\sigma,\xi)&= \frac{(n-2)\e}{|\ln
\e|}\biggr\{\frac{\Gamma_1}{(2(k-h)+1)\rho_h}-\frac{\ov
c_1}{2d_h}\biggr(\frac{\rho_{h+1}}{\rho_h}\biggr)^{(n-2)/2}
 +\frac{\ov c_1}{2\rho_h}\biggr(\frac{\rho_{h}}{\rho_{h-1}}\biggr)^{(n-2)/2} \biggr\}
+o\biggr(\frac{\e}{|\ln \e|}\biggr) ,\\ \frac{\partial
\widetilde{I_\e}}{\partial \rho_1}(\ov\a,\ov \rho,\ov \sigma,\xi)&
=\frac{(n-2)\e}{|\ln
\e|}\biggr\{\frac{\Gamma_1}{(2k-1)\rho_1}-\frac{\ov c_1}{2\rho_1}
\biggr(\frac{\rho_{2}}{\rho_1}\biggr)^{(n-2)/2}\biggr\}+o\biggr(\frac{\e}{|\ln \e|}\biggr),\\
\frac{\partial \widetilde{I_\e}}{\partial \rho_k}(\ov\a,\ov \rho,\ov
\sigma,\xi)&= \frac{(n-2)\e}{|\ln
\e|}\biggr\{\frac{\Gamma_1}{\rho_k}+\frac{\ov
c_1}{2\rho_k}\biggr(\frac{\rho_{k}}{\rho_{k-1}}\biggr)^{(n-2)/2}
-\frac{\ov c_1}{2\rho_k}\frac{H(\xi,\xi)}{\rho_k^{n-2}}\biggr\}+o\biggr(\frac{\e}{|\ln \e|}\biggr),\\
\frac{\partial \widetilde{I_\e}}{\partial (\sigma_r)_j}(\ov\a,\ov
\rho,\ov \sigma,\xi)&= -\frac{(n-2)\e}{|\ln
\e|}\Big(\frac{\rho_{r}}{\rho_{r-1}}\Big)^{\frac{n-2}{2}}
(\sigma_r)_j+o\biggr(\frac{\e}{|\ln \e|}\biggr),
\end{align*}
$C^0$-uniformly with respect to $\xi$ in compact sets of $\O$, $\ov
\sigma$ in compact sets of $\R^{(k-1)n}$, $\ov \a$ in compact sets
of $\R^k_+$ and $\ov \rho$ in compact sets of $\R^k_+$.
\end{pro}
\begin{pf}Let $\partial_s$ denote $\partial_{\rho_h}$ for $h = 1,\ldots, k$ and
$\partial_{(\sigma_r)_j}$ for $r = 2,\ldots, k $ and $j =
1,\ldots,n$. We have
\begin{equation}\label{531}
\partial_s \widetilde{I_\e}(\ov\a,\ov \rho,\ov \sigma,\ov\xi)=\n I_\e(V+\ov v)[\partial_{s}V+\partial _{s}\ov v]
\end{equation}
where $\partial_s$ is the partial derivative with respect to the variable $s$.\\
For $h=1,...,k$, in view of \eqref{equ4} and \eqref{equ5}, we have
\begin{align}\label{gradient}
\n I_\e(V+\ov v)[\partial_{\rho_h}V]=&\a_h\biggr(\frac{\e}{\ln
\e}\biggr)^{-\frac{2(k-h)+1}{n-2}}
\Big(\n I_\e(V+\ov v),(-1)^h\frac{\partial P\d_h}{\partial\l_h}\Big)\nonumber\\
&-\a_h\biggr(\frac{\e}{|\ln
\e|}\biggr)^{\frac{2(k-h)+1}{n-2}}\rho_h^{-2}\sum_{j=1}^n(\sigma_{h})_j\Big(\n
I_\e(V+\ov v),(-1)^h\frac{\partial P\d_h}{\partial(\xi_h)_j}\Big).
\end{align}
An easy computations show that
\begin{equation}\label{equa21}
\displaystyle\ln \l_h
=\frac{(2(k-h)+1)|\ln\e|}{n-2}\left(1+O\big(\frac{\ln|\ln\e|}{|\ln\e|}\big)\right)
\quad\hbox{ and }\quad\frac{1}{\ln \l_h}
=\frac{n-2}{(2(k-h)+1)|\ln\e|}+O\left(\frac{\ln|\ln\e|}{|\ln\e|^2}\right).
\end{equation}
Therefore, using Propositions \ref{p23}, \ref{p24},
\eqref{gradient}, \eqref{equa21} and \eqref{cr11}, we obtain
\begin{align}\label{532}
 \n I_\e(V+\ov v)[\partial_{\rho_h}V]&=(n-2)\frac{\e}{|\ln \e|}\biggr\{\frac{\Gamma_1}{(2(k-h)+1)\rho_h}
 -\frac{\ov c_1}{2\rho_h}\biggr(\frac{\rho_{h+1}}{\rho_h}\biggr)^{(n-2)/2}\nonumber
 \\&+\frac{\ov c_1}{2\rho_h}\biggr(\frac{\rho_{h}}{\rho_{h-1}}\biggr)^{(n-2)/2} \biggr\}
+o\biggr(\frac{\e}{|\ln \e|}\biggr) \quad \quad\mbox{for}\;
h=2,...,k-1,\nonumber\\ \n I_\e(V+\ov v)[\partial_{\rho_1}V]&=
(n-2)\frac{\e}{|\ln
\e|}\biggr\{\frac{\Gamma_1}{(2m-1)\rho_1}-\frac{c_1}{2\rho_1}
\biggr(\frac{\rho_{2}}{\rho_1}\biggr)^{(n-2)/2}\biggr\}+o\biggr(\frac{\e}{|\ln
\e|}\biggr),\nonumber\\
\n I_\e(V+\ov v)[\partial_{\rho_k}V]&=(n-2)\frac{\e}{|\ln
\e|}\biggr\{\frac{\Gamma_1}{\rho_k}+\frac{\ov
c_1}{2\rho_k}\biggr(\frac{\rho_{k}}{\rho_{k-1}}\biggr)^{(n-2)/2}
-\frac{\ov
c_1}{2\rho_k}\frac{H(\xi,\xi)}{\rho_k^{n-2}}\biggr\}+o\biggr(\frac{\e}{|\ln
\e|}\biggr).
\end{align}
In similar way, from Proposition \ref{p24} we derive, for
$r=2,...,k$ and $j=1,...,n$ that
\begin{align}\label{533}
\n I_\e(V+\ov v)[\partial_{(\sigma_r)_j}V]&= (-1)^r\a_r\Big(\n
I_\e(V+\ov v),\frac{1}{\l_r}\frac{\partial
P\d_r}{\partial(\xi_r)_j}\Big)\nonumber\\&= -\frac{\e}{|\ln
\e|}(n-2)\Big(\frac{\rho_{r}}{\rho_{r-1}}\Big)^{(n-2)/2}
(\sigma_r)_j+o\biggr(\frac{\e}{|\ln \e|}\biggr).
\end{align}
We shall prove now
\begin{equation}\label{534}
\n I_\e(V+\ov v)[\partial _{s}\ov v]=o\biggr(\frac{\e}{|\ln
\e|}\biggr).
\end{equation}
Recall that
\begin{align}\label{516}
\n I_\e(V+\ov v)[\partial _{s}\ov
v]&=<\sum_{i=1}^k\Big(A_iP\delta_i+B_i\frac{\partial
P\delta_i}{\partial \lambda_i}+\sum_{j=1}^n C_{ij}\frac{\partial
P\delta_i}{\partial (\xi_i)_j}\Big),\partial_{s} \ov
v>\notag\\&=\sum_{i=1}^kA_i<P\delta_i,\partial_{s} \ov
v>+\sum_{i=1}^kB_i<\frac{\partial P\delta_i}{\partial
\lambda_i},\partial_{s} \ov v>+\sum_{i=1}^k\sum_{j=1}^n
C_{ij}<\frac{\partial P\delta_i}{\partial (\xi_i)_j},\partial_{s}\ov
v> \notag\\&=\sum_{i=1}^kA_i<\partial_{s}P\delta_i, \ov
v>+\sum_{i=1}^kB_i<\partial_{s}\frac{\partial P\delta_i}{\partial
\lambda_i}, \ov v>+\sum_{i=1}^k\sum_{j=1}^n
C_{ij}<\partial_{s}\frac{\partial P\delta_i}{\partial (\xi_i)_j},\ov
v>
\end{align}
since $\ov v\in E$. Moreover, in view of \eqref{equ4} and
\eqref{equ5}, it is easy to check that
\begin{align}\label{cr8}
& \|\partial_{\rho_h}P\d_{h}\|=O(1),\quad
\|\partial_{(\sigma_{h})_j}P\d_{h}\|=O(1),\quad
\|\partial_{(\xi_h)_j}P\d_{h}\|=O(\l_h),\notag\\&\Big\|\partial_{\rho_h}\frac{\partial
P\d_h}{\partial \l_h}\Big\|=O(\frac{1}{\l_h}),\quad
\Big\|\partial_{(\sigma_{h})_j}\frac{\partial P\d_h}{\partial
\l_h}\Big\|=O(\frac{1}{\l_h}),\quad
\Big\|\partial_{(\xi_h)_j}\frac{\partial P\d_h}{\partial
\l_h}\Big\|=O(1),\notag\\& \Big\|\partial_{\rho_h}\frac{\partial
P\d_h}{\partial \xi_h}\Big\|=O(\l_h),\quad
\Big\|\partial_{(\sigma_{h})_j}\frac{\partial P\d_h}{\partial
\xi_h}\Big\|=O(\l_h),\quad\Big\|\partial_{(\xi_h)_j}\frac{\partial
P\d_h}{\partial \xi_h}\Big\|=O(\l_h^2), \quad\mbox{for}\;
j=1,\ldots,n.
\end{align}
 Through \eqref{516}, \eqref{cr8}
and Proposition \ref{pro42}, for $s=\rho_h$, we get
\begin{align*}
\n I_\e(V+\ov v)[\partial _{\rho_h}\ov
v]&=O\biggr(\Big(|A_h|\|\partial_{\rho_h}P\d_h\|+|B_h|\Big\|\partial_{\rho_h}\Big(\frac{\partial
P\d_h}{\partial\l_h}\Big)\Big\| +\sum_{j=1}^n
|C_{hj}|\Big\|\partial_{\rho_h}\Big(\frac{\partial
P\d_h}{\partial(\xi_h)_j}\Big)\Big\|\Big)\|\ov
v\|\biggr)=o\biggr(\frac{\e}{|\ln \e|}\biggr),
\end{align*}
and for $s=(\sigma_r)_j$, we obtain
\begin{align*}
\n I_\e(V+\ov v)[\partial _{(\sigma_r)_j}\ov
v]&=O\biggr(\Big(|A_r|\|\partial_{(\sigma_r)_j}P\d_r\|+|B_r|\Big\|\partial_{(\sigma_r)_j}\Big(\frac{\partial
P\d_r}{\partial\l_r}\Big)\Big\| +\sum_{\ell=1}^n
|C_{r\ell}|\Big\|\partial_{(\sigma_r)_j}\Big(\frac{\partial
P\d_r}{\partial(\xi_r)_\ell}\Big)\Big\|\Big)\|\ov
v\|\biggr)\\
&=o\biggr(\frac{\e}{|\ln \e|}\biggr).
\end{align*}
Therefore, \eqref{534} follows. Finally, combining \eqref{531} and
\eqref{532}-\eqref{534}, the proof of is Proposition \ref{pro53}
completed.
\end{pf}
\begin{pro}\label{pro54}For $i=1,...,n$, we have
\begin{equation*}
\frac{\partial \widetilde{I_\e}}{\partial (\xi)_i}(\ov\a,\ov
\rho,\ov \sigma,\xi)= \frac{\ov c_1}{2\rho_k^{n-2}}\frac{\e}{|\ln
\e|}\partial_{(x)_i}R(\xi)+o\biggr(\frac{\e}{|\ln \e|}\biggr).
\end{equation*}
$C^0$-uniformly with respect to $\xi$ in compact sets of $\O$, $\ov
\sigma$ in compact sets of $\R^{(k-1)n}$, $\ov \a$ in compact sets
of $\R^k_+$ and $\ov \rho$ in compact sets of $\R^k_+$.
\end{pro}
\begin{pf}
Let $ u=V+\ov v$ and set
\begin{equation}\label{v12}
S( u):=-\Delta u-\frac{| u|^{p-1}
u}{[\ln(e+|u|)]^{\e}}=-\sum_{m=1}^k\Delta\Big(A_mP\delta_m+B_m\frac{\partial
P\delta_m}{\partial \lambda_m}+\sum_{j=1}^n C_{mj}\frac{\partial
P\delta_m}{\partial (\xi_m)_j}\Big)
\end{equation}
because of \eqref{v1}.\\
Using formulas $(2.7)$ and $(2.9)$ in \cite{R}, for $m=1\ldots k$
and $i=1\ldots,n$, we have
\begin{align}\label{r1}
\int_{\partial \O}|\n P\d_{(\xi_m,\l_m)}|^2\nu_idw&=\frac{\ov
c_1}{\l^{n-2}}
\partial_{(x)_i}R(\xi_m)+O\biggr(\frac{1}{ \l_m  ^{n-1}}\biggr),\\
|\n P\d_{(\xi_m,\l_m)}|_{L^2(\partial
\O)}&=O(\l_m^{-\frac{n-2}{2}})\label{r2},
\end{align}
since the concentration points $\xi_1,\ldots,\xi_k$ are close to the point $\xi$ which is far away from the boundary of $\O$.\\
 Let $i\in \{1,\ldots,n\}$, we claim that
\begin{equation}\label{claim}
\frac{\partial \widetilde{I_\e}}{\partial (\xi)_i}(\ov\a,\ov
\rho,\ov \sigma,\ov\xi)=\n I_\e(V+\ov
v)[\partial_{(\xi)_i}V+\partial _{(\xi)_i}\ov
v]=\int_{\O}S(u)\partial_{(\xi)_i} u=-\int_{\O}S(u)\partial_{(x)_i}
u+o\biggr(\frac{\e}{|\ln \e|}\biggr).
\end{equation}
In fact, by integrating by parts and using \eqref{v12}, \eqref{q1},
\eqref{equ5}, \eqref{r2} and Proposition \ref{pro42}, we have
\begin{align*}
\int_{\O}S( u)\partial_{(\xi)_i}
V&=-\sum_{m=1}^k\int_{\O}\Delta\Big(A_mP\delta_m+B_m\frac{\partial
P\delta_m}{\partial \lambda_m}+\sum_{j=1}^n C_{mj}\frac{\partial
P\delta_m}{\partial (\xi_m)_j}\Big)\partial_{(\xi)_i}
V\\&=\sum_{m=1}^k\int_{\O}\Big(A_mP\delta_m+B_m\frac{\partial
P\delta_m}{\partial \lambda_m}+\sum_{j=1}^n C_{mj}\frac{\partial
P\delta_m}{\partial (\xi_m)_j}\Big)\partial_{(\xi)_i} (-\Delta
V)\\&=\sum_{m=1}^k\int_{\O}\Big(A_mP\delta_m+B_m\frac{\partial
P\delta_m}{\partial \lambda_m}+\sum_{j=1}^n C_{mj}\frac{\partial
P\delta_m}{\partial (\xi_m)_j}\Big)\partial_{(x)_i} (\Delta
V)\\&=\sum_{m=1}^k\int_{\O}\Delta\Big(A_mP\delta_m+B_m\frac{\partial
P\delta_m}{\partial \lambda_m}+\sum_{j=1}^n C_{mj}\frac{\partial
P\delta_m}{\partial (\xi_m)_j}\Big)\partial_{(x)_i}
V\\&-\sum_{m=1}^k\int_{\partial\O}\frac{\partial}{\partial
\nu}\Big(A_mP\delta_m+B_m\frac{\partial P\delta_m}{\partial
\lambda_m}+\sum_{j=1}^n C_{mj}\frac{\partial P\delta_m}{\partial
(\xi_m)_j}\Big)\partial_{(x)_i}  V\\&=-\int_{\O}S(
u)\partial_{(x)_i} V+o\biggr(\frac{\e}{|\ln \e|}\biggr).
\end{align*}
Using \eqref{v12}, \eqref{q1}, $\ov v\in E$, \eqref{equ5} and
integration by parts, we also have
\begin{align*}
\int_{\O}S(u)\partial_{(\xi)_i} \ov
v&=\sum_{m=1}^k\bigg<\Big(A_mP\delta_m+B_m\frac{\partial
P\delta_m}{\partial \lambda_m}+\sum_{j=1}^n C_{mj}\frac{\partial
P\delta_m}{\partial (\xi_m)_j}\Big),\partial_{(\xi)_i} \ov v\bigg>\\
&=-\sum_{m=1}^k\bigg<\partial_{(\xi)_i}\Big(A_mP\delta_m+B_m\frac{\partial
P\delta_m}{\partial \lambda_m}+\sum_{j=1}^n C_{mj}\frac{\partial
P\delta_m}{\partial (\xi_m)_j}\Big), \ov v\bigg>
\\&=-\sum_{m=1}^k\int_{\O}\Big(A_m\partial_{(\xi)_i}(\delta_m^p)+pB_m\partial_{(\xi)_i}(\delta_m^{p-1}\frac{\partial
\delta_m}{\partial \lambda_m})+p\sum_{j=1}^n
C_{mj}\partial_{(\xi)_i}(\delta_m^{p-1}\frac{\partial
\delta_m}{\partial (\xi_m)_j}\big)\Big) \ov
v\\&=\sum_{m=1}^k\int_{\O}\Big(A_m\partial_{(x)_i}(\delta_m^p)+pB_m\partial_{(x)_i}(\delta_m^{p-1}\frac{\partial
\delta_m}{\partial \lambda_m})+p\sum_{j=1}^n
C_{mj}\partial_{(x)_i}(\delta_m^{p-1}\frac{\partial
\delta_m}{\partial (\xi_m)_j}\big)\Big) \ov
v\\&=-\sum_{m=1}^k\int_{\O}\big(A_m\delta_m^p+pB_m\delta_m^{p-1}\frac{\partial
\delta_m}{\partial \lambda_m}+p\sum_{j=1}^n
C_{mj}\delta_m^{p-1}\frac{\partial \delta_m}{\partial
(\xi_m)_j}\big)\partial_{x_i} \ov v\\&=-\int_{\O}S(
u)\partial_{(x)_i} \ov v.
\end{align*}
Our claim \eqref{claim} is proved.\\
 Moreover, using integration by
parts, we have
\begin{align*}
-\int_{ \O}\Delta  u\partial_{(x)_i} u =-\int_{ \O}
u\Delta(\partial_{(x)_i}u)-\int_{\partial \O}\biggr(\frac{\partial
u}{\partial\nu}\biggr)^2\nu_idw =\int_{ \O}
\partial_{(x)_i} u\,\Delta u-\int_{\partial \O}\biggr(\frac{\partial
u}{\partial\nu}\biggr)^2\nu_idw,
\end{align*}
thus
\begin{equation}\label{541}
\int_{\O}S(u)\partial_{(x)_i} u=-\int_{ \O}\Delta u\partial_{(x)_i}
u-\int_{\O}\frac{| u|^{p-1}u}{[\ln(e+|u|)]^{\e}}\partial_{(x)_i}
u=-\frac{1}{2}\int_{\partial \O}|\n u|^2\nu_idw
\end{equation}
since $\displaystyle\int_{\O}\frac{| u|^{p-1}u}{[\ln(e+|u|)]^{\e}}\partial_{(x)_i} u=\int_{\O}\partial_{(x)_i}
(F_{\e}(u(x)))=0$. \\
Now, by using \eqref{equ4}, \eqref{r1}, \eqref{r2} and Lemma
\ref{lem45}, we obtain
\begin{align}\label{542}
\int_{\partial \O}|\n u|^2\nu_idw
=&\int_{\partial \O}|\n V|^2\nu_idw+\int_{\partial \O}|\n \ov
v|^2\nu_idw+2\int_{\partial \O}\n V \n \ov v\nu_idw\notag\\
=&\sum_{m=1}^k\int_{\partial \O}|\n
P\d_{m}|^2\nu_idw+\sum_{l,h=1,l\neq h}^k\int_{\partial \O}\n
P\d_{h}\n P\d_{l}\nu_idw +\int_{\partial \O}|\n \ov
v|^2\nu_idw\notag\\&+O\biggr(\sum_{m=1}^k|\n
P\d_m|_{2,\partial\O}|\n \ov v|_{2,\partial \O}\biggr)\notag\\=&
\int_{\partial \O}|\n P\d_{k}|^2\nu_idw+o\biggr(\frac{\e}{|\ln \e|}\biggr)\notag\\
&=\frac{\ov
c_1}{\rho_k^{n-2}}\frac{\e}{|\ln\e|}\partial_{(x)_i}R(\xi_k)+o\biggr(\frac{\e}{|\ln
\e|}\biggr).
\end{align}
The desired result follows by \eqref{claim}, \eqref{541},
\eqref{542} and \eqref{equ5}.
\end{pf}
\section{Proof of Theorem \ref{t:1}}
We recall that, in the previous section, we have mentioned that our
aim is prove that the reduced energy has a critical point. This fact
is equivalent to the existence of the requested solution.\\
 Let us perform the change of variables:
$$s_i=\rho_{i+1}/\rho_{i},\quad\mbox{for}\;1\leq i\leq k-1\quad\mbox{and}\;\;s_k=\rho_k.$$
The function $\Psi$ in the new variables $\ov s=(s_1,...,s_k)$ reads
as
\begin{align*}\widehat{\Psi}(\ov s,\ov \sigma,\xi)&:=-\frac{(n-2)^2S_n^{\frac{n}{2}}}{2n}
\sum_{i=1}^{k-1}\sum_{l=k-(i-1)}^k\frac{1}{2l-1}\ln s_i +\ov
c_1\sum_{i=1}^{k-1}s_i^{\frac{n-2}{2}}\frac{1}{\big(1
+|\sigma_{i+1}|^2\big)^{(n-2)/2}}\\&+\frac{(n-2)^2S_n^{\frac{n}{2}}}{2n}\ln
s_k\sum_{l=1}^k\frac{1}{2l-1}+\frac{\ov
c_1}{2}\frac{H(\xi,\xi)}{s_k^{n-2}}.
\end{align*}
Let $\xi_0$ be a stable critical point of the Robin function.
 For any
$\ov\sigma\in\R^{(k-1)n}$ and for any $\xi \in B(\xi_0,\varrho)$ the
function $\ov s\rightarrow \widehat{\Psi}(\ov s,\ov \sigma,\xi)$ has
exactly one critical point $\ov s^{(1)}(\ov \sigma,\xi)$. More
precisely
\begin{align*}&s^{(1)}_q=\Biggr[\frac{(n-2)S_n^{\frac{n}{2}}}{n\ov c_1}\sum_{l=k-(q-1)}^k\frac{1}{2l-1}(1+|\sigma_{q+1}|^2)^{\frac{n-2}{2}}\Biggr]^{\frac{2}{n-2}}, \quad\mbox{for}\;q=1,...,k-1,\\&
s^{(1)}_k=\Biggr[\frac{n\ov
c_1}{(n-2)S_n^{\frac{n}{2}}}\frac{H(\xi,\xi)}{\sum_{l=1}^k\frac{1}{2l-1}}\Biggr]^{\frac{1}{n-2}}.
\end{align*}
It is useful to point out that $s_k^{(1)}\rightarrow \infty$ as $\xi\rightarrow \partial \O$. Moreover, it is easy to check that it is non-degenerate, i.e.
the matrix $(\partial^2_{s_hs_l}\widehat{\Psi}(\ov s^{(1)},\ov \sigma,\xi))_{h,l=1,...,k}$ is invertible.

Now, let us consider the reduced function:
\begin{align*}
(\ov\sigma,\xi)\rightarrow \widehat{\Psi}(\ov s^{(1)},&\ov
\sigma,\xi) =\frac{(n-2)^2S_n^{\frac{n}{2}}}{2n}\ln
s_k^{(1)}\sum_{l=1}^k\frac{1}{2l-1}
+\frac{(n-2)S_n^{\frac{n}{2}}}{2n}\sum_{l=1}^k\frac{1}{2l-1}+\frac{(n-2)S_n^{\frac{n}{2}}}{n}\sum_{i=1}^{k-1}
\sum_{l=k-(i-1)}^k\frac{1}{2l-1}\\&
-\frac{2}{n-2}\ln\biggr(\frac{(n-2)S_n^{\frac{n}{2}}}{n\ov c_1}
\sum_{i=1}^{k-1}\sum_{l=k-(i-1)}^k\frac{1}{2l-1}\biggr)-\frac{(n-2)^2S_n^{\frac{n}{2}}}{2n}
\sum_{i=1}^{k-1}\sum_{l=k-(i-1)}^k\frac{1}{2l-1}\ln
(1+|\sigma_{i+1}|^2).
\end{align*}
Since $\xi_0$ is a stable critical point of  $\xi\rightarrow \ln
s_k^{(1)}$ and each $\sigma_i = 0$ is a strict maximum point of
$\sigma_i\rightarrow-\ln(1+|\sigma_i|^2)$, by degree theory we
easily deduce that the point $(0,\xi_0)$ is a $C^1$-stable critical
point of the reduced function, in the sense that small
$C^1$-perturbations of the reduced function still have a critical
point, close to $(0,\xi_0)$.\\ 
Concerning the variable $\ov\alpha$, let us denote the functions
$$\phi(\alpha):=\frac{\alpha^2}{2}-\frac{\alpha^{p+1}}{p+1}\, \hbox{ and }\, \Phi(\ov\alpha):=S_n^\frac{n}{2}
\sum_{i=1}^k\phi(\a_i)$$ defined on $\mathbb{R}_+$ and
$(\mathbb{R}_+)^k$ respectively. $\phi$ is a $C^2$-function
satisfying $\phi'(1)=0$ and $\phi''(1)=1-p<0$. Thus $1$ is a
non-degenerate maximum point of $\phi$ and we derive that
$\ov\alpha^{(1)}:=(1,\ldots,1)$ is $C^1$-stable critical point of
$\Phi$.\\
Therefore, we can conclude that,  the point $(\ov\alpha^{(1)},\ov
s^{(1)}, 0,\xi_0)$ is a $C^1$-stable critical point for
$\Phi+\frac{\e}{|\ln\e|}\widehat{\Psi}$, in the sense that small
$C^1$-perturbations of the function
$\Phi+\frac{\e}{|\ln\e|}\widehat{\Psi}$ still have a critical point,
close to $(\ov\alpha^{(1)},\ov s^{(1)}, 0,\xi_0)$. Propositions
\ref{pro52}, \ref{pro53} and \ref{pro54} conclude the proof of
Theorem \ref{t:1}.
\section{Appendix}
 In this section, we collect some technical Lemmas used in this
 paper. The first one concerns the
 nonlinearity $f_\e$ and its derivatives. It is contained in \cite{Pistoianew} and \cite{BFG}.
\begin{lem}\label{lemmaB1}
\begin{enumerate}
\item For any $\e>0$, and any $U\in \mathbb{R}$, we have $|f_\e(U)-f_0(U)|\leq \e |U|^{p}\ln\ln(e+|U|)$.
\item There exists $c>0$ such that, for $\e$ small enough
and any $U,V\in \mathbb{R}$
\begin{equation}\label{B'}
|f_\e(U+V)-f_\e(U)|\leq c( |U|^{p-1}+|V|^{p-1})|V|, \quad \forall
n\geq 3.
\end{equation}
\item For $\e$ small enough, and any $U\in \mathbb{R}$,
\begin{equation}\label{B.1} |f'_\e(U)|\leq c|U|^{p-1},\end{equation} and
\begin{equation}\label{f'}|f'_\e(U)-f'_0(U)|\leq \e
|U|^{p-1}\left(p\ln\ln(e+|U|)+\frac{1}{\ln(e+|U|)}\right).
\end{equation}
\item There exists $c>0$ such that, for $\e$ small enough and any $U,V\in
\mathbb{R}$,
\begin{equation}\label{B.2} |f'_\e(U+V)-f'_\e(U)|\leq \left\{\begin{array}{ll}c(
|U|^{p-2}+|V|^{p-2})|V|\hbox{ if }n\leq 6,\\
c( |V|^{p-1}+\e|U|^{p-1})\quad \hbox{ if }n> 6.\end{array}\right.
\end{equation}
\item For $\e$ small enough, and any $U\in \mathbb{R}$,
\begin{equation}\label{B''}
|f''_\e(U)|\leq c|U|^{p-2}, \quad \forall n\geq 3.
\end{equation}
\end{enumerate}
\end{lem}
We state the preliminary result given by  \cite[Proposition 1]{R} as
follows.
\begin{pro}\label{p21} Let $a\in
\O$ and $\l>0$ such that $\l d:=\l d(a, \partial\O)$ is large
enough. For $\varphi_{(a,\l)}=\d_{(a,\l)}-P\d_{(a,\l)}$, we have the
following estimates
\begin{eqnarray*}(a)\;\; 0\leq \varphi_{(a,\l)}\leq \d_{(a,\l)},\quad
(b)\quad \varphi_{(a,\l)}=c_0\frac{H(a,.)}{\l^{(n-2)/2}}+f_{(a,\l)},
\end{eqnarray*}
where $c_0$ is defined in \eqref{d} and $f_{(a,\l)}$ satisfies
\begin{align*}&f_{(a,\l)}=O\biggr(\frac{1}{\l^{\frac{n+2}{2}}d^{n}}\biggr),\quad
\l\frac{\partial f_{(a,\l)}}{\partial
\l}=O\biggr(\frac{1}{\l^{\frac{n+2}{2}}d^{n}}\biggr),
\\&\frac{1}{\l}\frac{\partial f_{(a,\l)}}{\partial
a}=O\biggr(\frac{1}{\l^{\frac{n+4}{2}}d^{n+1}}\biggr).
\end{align*}
\begin{align*}(c)\;&|\varphi_{(a,\l)}|_{2n/(n-2)}=O\biggr(\frac{1}{(\l
d)^{(n-2)/2}}\biggr),\quad
\biggr|\l\frac{\partial\varphi_{(a,\l)}}{\partial
\l}\biggr|_{2n/(n-2)}=O\biggr(\frac{1}{(\l
d)^{(n-2)/2}}\biggr),\\&\|\varphi_{(a,\l)}\|=O\biggr(\frac{1}{(\l
d)^{(n-2)/2}}\biggr),\quad
\biggr|\frac{1}{\l}\frac{\partial\varphi_{(a,\l)}}{\partial
a}\biggr|_{2n/(n-2)}=O\biggr(\frac{1}{(\l d)^{n/2}}\biggr),
\end{align*}
where $|.|_q$ denotes the usual norm in $L^q(\Omega)$ for each
$1\leq q\leq \infty$.
\end{pro}
We also need the following lemma.
\begin{lem}\label{lemme1.55}
\begin{enumerate}[(a)]
\item For every $y\in \Omega,$ we have:
$$  (y-a).\n_y P\delta_{(a,\lambda)} (y)=O(\delta_{(a,\lambda)}).$$
\item For every $y\in \Omega,$ we have:
$$   \n_y P\delta_{(a,\lambda)} (y)=O(\lambda\delta_{(a,\lambda)}).$$
\end{enumerate}
\end{lem}
\begin{pf}
We have
$$|\n_y P\d_{(a,\lambda)}(y)|=
\biggr|\int_{\Omega}\n_y
G(x,y)\d_{(a,\lambda)}^{\frac{n+2}{n-2}}(x)dx\biggr|\leq c
\int_{\Omega}
\frac{1}{|x-y|^{n-1}}\d_{(a,\lambda)}^{\frac{n+2}{n-2}}(x)dx,$$
where we have used that $|\nabla_y G(x,y)|\leq c/|x-y|^{n-1}$.
\begin{enumerate}[(a)]
\item To prove $(a)$ we distinguish two cases.\\
{\it First case:} if $\lambda|y-a|\geq 1$. In this case we get
\begin{equation}\label{*}
1/(\l^{(n-2)/2}|y-a|^{n-2})\leq c \d_{(a,\lambda)}(y).
\end{equation}
We have
\begin{align*}
I:
=|y-a|\int_{\Omega}\frac{1}{|x-y|^{n-1}}\d_{(a,\lambda)}^{\frac{n+2}{n-2}}(x)dx
 =|y-a|\int_{|x-y|\leq  |y-a|/ 2 }...+|y-a|\int_{|x-y|\geq |y-a|/ 2}...
=:I_1+I_2.
\end{align*}
Concerning  $I_1$, we have $\displaystyle
\delta_{(a,\lambda)}(x)\leq c\delta_{(a,\lambda)}(y)$ since
$|x-a|\geq |y-a|/2$. This implies that:
$$I_1\leq c\delta_{(a,\lambda)}^{\frac{n+2}{n-2}}(y)|y-a|\int_{|x-y|\leq  |y-a|/2}\frac{1}{|x-y|^{n-1}}dx\leq c\delta_{(a,\lambda)}\frac{\l^2|y-a|^2}{(1+\l^2|y-a|^2)^2}\leq c\delta_{(a,\lambda)}.$$
For $I_2$, since $\displaystyle |x-y|\geq \frac{1}{2 }|y-a|$, then
\begin{align*}
I_2\leq\frac{1}{|y-a|^{n-2}}\int_{\mathbb{R}^n}
\d_{(a,\lambda)}^{\frac{n+2}{n-2}} \leq
\frac{c}{\l^{(n-2)/2}|y-a|^{n-2}}
 \leq c\delta_{(a,\lambda)}(y)
\end{align*}
by using \eqref{*}.
Hence we get the desired result in this case.\\
{\it Second case:}
  If $\lambda|y-a|\leq 1.$
 In this case we have:
\begin{equation}\label{**}\displaystyle
\delta_{(a,\lambda)}(y)\geq c\l^{(n-2)/2}.
\end{equation}
We split the integral $I$ as follows
\begin{align*}
I:=|y-a|\int_{\Omega}\frac{1}{|x-y|^{n-1}}\d_{(a,\lambda)}^{\frac{n+2}{n-2}}(x)dx=|y-a|\int_{|x-a|\leq
\frac{M}{\lambda}}...+|y-a|\int_{|x-a|\geq \frac{M}{\lambda}}...
=:J_1+J_2,
\end{align*}
where $M$ is a large positive constant.\\
Concerning $J_2,$ since we have:  $x$ satisfies
$\displaystyle\lambda|x-a|\geq M,$ we obtain $\l|x-y|\geq M-1$ and
$\displaystyle  |x-y|\geq c|a-x|.$ Thus, in view of \eqref{**}, we
get
\begin{align}\label{1.553}
J_2&\leq|y-a|\int_{|x-a|\geq
\frac{M}{\lambda}}\frac{\l^{(n+2)/2}}{|a-x|^{n-1}(1+\l^2|x-a|^2)^{(n+2)/2}}\notag\\&\leq
c\l^{n/2}|y-a|\int_M^{\infty}\frac{1}{(1+r^2)^{(n+2)/2}} \leq
c\l^{(n-2)/2}\leq c\d_{(a,\lambda)}(y).
\end{align}
 For $J_1$, we have:
\begin{equation}\label{1.554}
J_1\leq c|y-a|\l^{(n+2)/2}\int_{|x-y|\leq
\frac{M+1}{\lambda}}\frac{dx}{|x-y|^{n-1}}\leq c|y-a|\l^{n/2}\leq
c\l^{(n-2)/2}\leq c\d_{(a,\lambda)}(y)
\end{equation}
where the last inequality follow from \eqref{**}.\\
\eqref{1.553} and \eqref{1.554} prove part $(a)$ of the lemma in this case.\\
\item The proof of $(b)$ is similar to that of $(a)$ taking into account power changes.
\end{enumerate}
\end{pf}

\end{document}